\documentclass[11pt]{article}
\pdfoutput=1

\usepackage[margin=1in]{geometry}
\usepackage{amssymb,amsbsy,exscale,amsmath,amsfonts,amssymb,amscd}
\usepackage{amsthm}
\usepackage{epsfig}
\usepackage{graphicx}
\usepackage{makeidx}
\usepackage{multicol}
\usepackage{bm} 
\usepackage{subfigure}
\usepackage{accents}
\usepackage{color}
\usepackage{placeins}

\newcommand{\rg}{ \mathrm {\bf {sp}}}

\begin{document}
\title{Band Gap Optimization of Two-Dimensional Photonic Crystals Using Semidefinite Programming and Subspace Methods\footnote{This
research has been supported through AFOSR grant FA9550-08-1-0350 and the Singapore-MIT Alliance.}}

\author{H. Men\footnote{National University of Singapore, Center for Singapore-MIT Alliance, Singapore 117576, email: \texttt{men@nus.edu.sg}},
N.~C. Nguyen\footnote{MIT Department of Aeronautics and Astronautics, 77 Massachusetts Ave., Cambridge, MA 02139, USA, email: \texttt{cuongng@mit.edu}},
R.~M. Freund\footnote{MIT Sloan School of Management, 50 Memorial Drive, Cambridge, MA 02142, USA, email: \texttt{rfreund@mit.edu}},
P.~A. Parrilo\footnote{MIT Department of Electrical Engineering and Computer Science, 77 Massachusetts Ave., Cambridge, MA 02139, USA, email: \texttt{parrilo@mit.edu}},
and J. Peraire\footnote{MIT Department of Aeronautics and Astronautics, 77 Massachusetts Ave., Cambridge, MA 02139, USA, email: \texttt{peraire@mit.edu}}}

\maketitle


\begin{abstract}
In this paper, we consider the optimal design of photonic crystal band
structures for two-dimensional square lattices. The mathematical
formulation of the band gap optimization problem leads to an
infinite-dimensional Hermitian eigenvalue optimization problem
parametrized by the dielectric material and the wave vector. To make
the problem tractable, the original eigenvalue problem is discretized
using the finite element method into a series of finite-dimensional
eigenvalue problems for multiple values of the wave vector
parameter. The resulting optimization problem is large-scale and
non-convex, with low regularity and non-differentiable objective. By
restricting to appropriate eigenspaces, we reduce the large-scale
non-convex optimization problem via reparametrization to a sequence of
small-scale convex semidefinite programs (SDPs) for which modern SDP
solvers can be efficiently applied. Numerical results are presented
for both transverse magnetic (TM) and transverse electric (TE)
polarizations at several frequency bands. The optimized structures
exhibit patterns which go far beyond typical physical intuition on
periodic media design.
\end{abstract}

\section{Introduction}

The propagation of waves in periodic media has attracted considerable
interest in recent years. This interest stems from the possibility of
creating periodic structures that exhibit band gaps in their spectrum,
i.e., frequency regions in which the wave propagation is
prohibited. Band gaps occur in many wave propagation phenomena
including electromagnetic, acoustic and elastic waves. Periodic
structures exhibiting electromagnetic wave band gaps, or photonic
crystals, have proven very important as device components for
integrated optics including frequency filters \cite{fan1998cdf},
waveguides \cite{fan1995gad}, switches \cite{SoljacicIb02}, and
optical buffers \cite{yanik2005sas}.

The optimal conditions for the appearance of gaps were first studied
for one-dimensional crystals by Lord Rayleigh in $1887$
\cite{rayleigh1887mvf}. In a one-dimensional periodic structure, one can
widen the band gap by increasing the contrast in the refractive index
and difference in width between the materials. Furthermore, it is
possible to create band gaps for any particular frequency by changing
the periodicity length of the crystal. Unfortunately, however, in two
or three dimensions one can only suggest rules of thumb for the
existence of a band gap in a periodic structure, since no rigorous
criteria have yet been determined. This made the design of two- or
three-dimensional crystals a trial and error process, being far from
optimal. Indeed, the possibility of two- and three-dimensionally
periodic crystals with corresponding two- and three-dimensional band
gaps was not suggested until 100 years after Rayleigh's discovery of
photonic band gap in one dimension, by
Yablonovitch~\cite{yablonovitch1987ise} and John~\cite{john1987slp} in 1987.

From a mathematical viewpoint, the calculation of the band gap reduces
to the solution of an infinite-dimensional Hermitian eigenvalue
problem which is parametrized by the dielectric function and the wave
vector. In the design setting, however, one wishes to know the answer
to the question: which periodic structures, composed of arbitrary
arrangements of two or more different materials, produce the largest
band gaps around a certain frequency?  This question can be rigorously
addressed by formulating an optimization problem for the parameters
that represent the material properties and geometry of the periodic
structure.  The resulting problem is infinite-dimensional with an
infinite number of constraints.  After appropriate discretization in
space and consideration of a finite set of wave vectors, one obtains a
large-scale finite-dimensional eigenvalue problem which is non-convex
and is known to be non-differentiable when eigenvalue multiplicities
exist.  The current state-of-the-art work done on this problem falls
into two broad categories.  The first kind tries to find the
``optimal'' band structure by parameter studies -- based on prescribed
inclusion shapes (e.g., circular or hexagonal inclusions)
\cite{doosje2000pbo} or fixed topology \cite{yang2008obg}. The second kind
attempts to use formal topology optimization
techniques~\cite{sigmund2003sdp,
cox2000bso, burger87ipt, kao2005mbg}. Both approaches
typically use gradient-based optimization methods.  While these
methods are attractive and have been quite successful in practice, the
optimization processes employed explicitly compute the sensitivities
of eigenvalues with respect to the dielectric function, which are
local subgradients for such non-differentiable problem.  As a result,
gradient-based solution methods often suffer from the lack of
regularity of the underlying problem when eigenvalue multiplicities
are present, as they typically are at or near the solution.

In this paper we propose a new approach based on semidefinite
programming (SDP) and subspace methods for the optimal design of
photonic band structure. In the last two decades, SDP has emerged as
the most important class of models in convex optimization; see
\cite{alizadeh1995ipm,alizadeh1998pdi,nesterov1994ipp,vandenberghe1996sp,wolkowicz2000hsp}.  SDP encompasses a huge array of convex
problems as special cases, and is computationally tractable (usually
comparable to least-square problems of comparable dimensions). There
are three distinct properties that make SDP very suitable for the band
gap optimization problem. First, the underlying differential operator
is Hermitian and positive semidefinite. Second, the objective and
associated constraints involve bounds on eigenvalues of matrices. And
third, as explained below, we can approximate the original non-convex
optimization problem by a semidefinite program for which SDP can be
well applied, thanks to its efficiency and robustness of handling this
type of spectral objective and constraints.

In our approach, we first reformulate the original problem of
maximizing the band gap between two consecutive eigenvalues as an
optimization problem in which we optimize the gap in eigenvalues
between two orthogonal subspaces.  The first eigenspace consists of
eigenfunctions corresponding to eigenvalues below the band gap,
whereas the second eigenspace consists of eigenfunctions whose
eigenvalues are above the band gap. In this way, the eigenvalues are
no longer present in our formulation; however, like the original
problem, the exactly reformulated optimization problem is
large-scale. To reduce the problem size, we truncate the
high-dimensional subspaces to only a few eigenfunctions below and
above the band gap~\cite{cances2007fac,pau:046704}, thereby
obtaining a new small-scale yet non-convex optimization problem.
Finally, we keep the subspaces fixed at a given decision parameter
vector and use a reparametrization of the decision variables to obtain
a convex semidefinite optimization problem for which SDP solution
methods can be effectively applied.  We apply this approach to
optimize band gaps in two-dimensional photonic crystals for either the
transverse magnetic (TM) or the transverse electric (TE)
polarizations.

The rest of the paper is organized as follows. In
Section~\ref{sec_BSC} we introduce the governing differential
equations and the mathematical formulation of the band gap
optimization problem. We then discuss the discretization process and
present the subspace restriction approach. In Section~\ref{sec_BSO} we
introduce the semidefinite programming formulation of the band
structure optimization, and lay out the optimization steps involved in
solving the problem. Numerical results are presented in
Section~\ref{sec_OR} for both the TE and TM polarizations in square
lattices. Finally, in Section~\ref{conclusion} we conclude with
several remarks on anticipated future research directions.

\section{The Band Gap Optimization Problem}
\label{sec_BSC}
\subsection{Governing Equations}

Our primary concern is the propagation of electromagnetic linear waves
in periodic media, and the design of such periodic structures, or
photonic crystals, to create optimal band gaps in their spectrum. The
propagation of electromagnetic waves in photonic crystals is governed
by Maxwell's equations. The solutions to these equations are in
general very complex functions of space and time. Due to linearity
however, it is possible to separate the time dependence from the
spatial dependence by expanding the solution in terms of harmonic
modes -- any time-varying solution can always be reconstructed by a
linear combination of these harmonic modes using Fourier analysis. By
considering only harmonic solutions, the problem is considerably
simplified since it reduces to a series of eigenvalue problems for the
spatially varying part of the solutions (eigenfunctions) and the
corresponding frequencies (eigenvalues).

In the absence of sources and assuming a monochromatic wave, i.e.,
with magnetic field $\bm{H}(\bm{r},t) = \bm{H}(\bm{r})e^{-i\omega t}$,
and electric field $\bm{E}(\bm{r},t) = \bm{E}(\bm{r})e^{-i\omega t}$,
Maxwell's equations can be written in the following form:
\begin{align*}
\bm{\nabla} \times \left(\frac{1}{\varepsilon(\bm{r})} \bm{\nabla} \times \bm{H}(\bm{r})\right)
&= \left( \frac{\omega}{c}\right)^2 \bm{H}(\bm{r}), \qquad \mbox{in } \mathbb{R}^3,\\
\frac{1}{\varepsilon(\bm{r})}  \bm{\nabla} \times \left(\bm{\nabla} \times \bm{E}(\bm{r})\right)
&= \left( \frac{\omega}{c}\right)^2 \bm{E}(\bm{r}), \qquad \mbox{in } \mathbb{R}^3,
\end{align*}
where $c$ is the speed of light, and $\varepsilon(\bm{r})$ is the
dielectric function. {In two dimensions, there are two possible
polarizations of the magnetic and electric fields. In TE (transverse
electric) polarization, the electric field is confined to the plane of
wave propagation and the magnetic field $\bm{H} = (0,0,H)$ is
perpendicular to this plane. In contrast, in TM (transverse magnetic)
polarization, the magnetic field is confined to the plane of wave
propagation and the electric field $\bm{E} = (0, 0, E)$ is
perpendicular to this plane. In such cases, the Maxwell's equations
can be reduced to scalar eigenvalue problems
\begin{equation}\label{eqBSC_2}
\mathrm{TE:} \qquad
   - \bm{\nabla} \cdot \left( \frac{1}{\varepsilon(\bm{r})} \, \bm{\nabla} {H}(\bm{r})\right) = \left( \frac{\omega}{c}\right)^2 {H}(\bm{r}), \qquad \mbox{in } \mathbb{R}^2,
\end{equation}
\begin{equation}\label{eqBSC_2B}
\mathrm{TM:} \qquad
- \bm{\nabla} \cdot \left(\bm{\nabla} {E}(\bm{r})\right) = \left( \frac{\omega}{c}\right)^2 \varepsilon(\bm{r})  {E}(\bm{r}), \qquad \mbox{in } \mathbb{R}^2 \ .
\end{equation}
Note that the reciprocal of the dielectric function is present in the
differential operator for the TE case, whereas the dielectric function
is present in the right-hand side for the TM case.}

\begin{figure}
\centering
\includegraphics[scale=0.5]{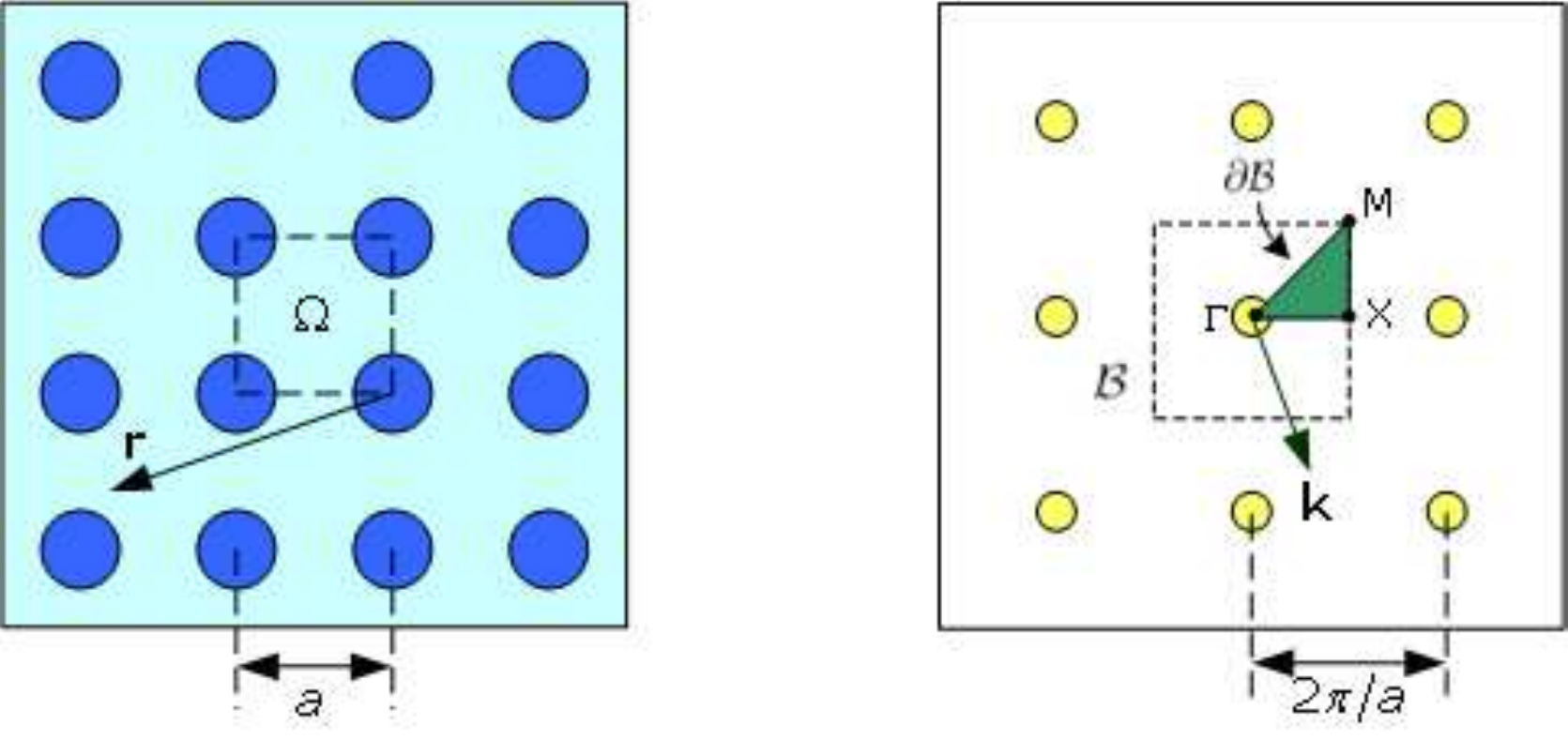}
\caption{Left: A photonic crystal on a square lattice. The dashed box
represents the primitive unit cell ($\Omega$), where $\bm a$ is the
periodicity length of the lattice. Right: The
reciprocal lattice, and the dashed box represents the first Brillouin zone ($\mathcal{B}$). The irreducible zone is the green triangular
wedge, and its boundary is denoted by $\partial \mathcal{B}$.}
\label{SquLat_BriZone}
\end{figure}

For two-dimensional square lattices the dielectric function satisfies
$\varepsilon(\bm{r}) = \varepsilon(\bm{r} + \bm{R})$, where $\bm{R}$
are the crystal lattice vectors\footnote{For a square lattice,
$\bm{R}$ denotes the vectors spanned by $\{a \bm{e}_1, a \bm{e}_2 \}$,
where $\bm{e}_1$ and $\bm{e}_2$ are the unit basis vectors and $a$ is
the periodicity length of the crystal \cite{joannopoulos2008pcm}.}. By
applying the Bloch-Floquet theory~\cite{bloch1929qek,floquet1883edl} for periodic
eigenvalue problems we obtain that
$$H(\bm{r}) = e^{i\bm{k}\cdot \bm{r}}H_{\bm{k}}(\bm{r}), \quad \mbox{and} \quad E(\bm{r}) = e^{i\bm{k}\cdot \bm{r}}E_{\bm{k}}(\bm{r}) , $$
where $H_{\bm{k}}(\bm{r})$ and $E_{\bm{k}}(\bm{r})$ satisfy
\begin{equation}\label{eqBSC_3}
\mathrm{TE:} \qquad
    (\bm{\nabla}+i\bm k) \cdot \left(\frac{1}{\varepsilon(\bm{r})} (\bm{\nabla}+i\bm k) H_{\bm k}(\bm{r})\right) =
\left( \frac{\omega}{c}\right)^2 H_{\bm k}(\bm{r}), \qquad \mbox{in } \Omega,
\end{equation}
\begin{equation}\label{eqBSC_3B}
\mathrm{TM:} \qquad
  (\bm{\nabla}+i\bm k) \cdot \left( (\bm{\nabla}+i\bm k) E_{\bm k}(\bm{r}) \right) =
\left( \frac{\omega}{c}\right)^2 \varepsilon(\bm{r}) E_{\bm k}(\bm{r}),  \qquad \mbox{in } \Omega,
\end{equation}
respectively. Thus, the effect of considering periodicity is reduced
to replacing the indefinite periodic domain by the unit cell $\Omega$
and $\bm{\nabla}$ by $\bm{\nabla} + i\bm k$ in the original equation,
where $\bm{k}$ is a wave vector in the first Brillouin zone
$\mathcal{B}$. Note that the unit cell $\Omega$ and the Brillouin zone
$\mathcal{B}$ depend on the lattice type (e.g., square or triangular
lattices) as well as the crystal lattice vectors $\bm{R}$. If we
further take into consideration the symmetry group of the square
lattice~\cite{weyl1952sp}, we only need to consider all possible
wavevectors $\bm k$ on the irreducible Brillouin zone, or (under
certain conditions) its boundary
\cite{joannopoulos2008pcm}. Figure~\ref{SquLat_BriZone} shows an
example of the unit cell and the Brillouin zone for a square lattice.


For notational convenience, we write the above equations in the following operator form
\begin{equation}
\label{eigen_cont}
\mathcal{A} u = \lambda \mathcal{M} u, \qquad \mbox{in } \Omega ,
\end{equation}
where, for the TE case, $u \equiv  H_{{\bm k}}(\bm{r})$, $\lambda \equiv \omega^2_{\text{TE}}/c^2$, and
\begin{equation}
\label{eqBSC_4}
\mathcal{A}(\varepsilon, \bm{k})  \equiv  - (\bm{\nabla}+i\bm k) \cdot \left( \frac{1}{\varepsilon(\bm{r})} (\bm{\nabla}+i\bm k) \right),
\qquad \mathcal{M} \equiv I ;
\end{equation}
whereas, for the TM case, $u \equiv  E_{{\bm k}}(\bm{r})$, $\lambda \equiv \omega^2_{\text{TM}}/c^2$, and
\begin{equation}
\label{eqBSC_4B}
\mathcal{A}(\bm{k})  \equiv - (\bm{\nabla}+i\bm k) \cdot (\bm{\nabla}+i\bm k),
\qquad \mathcal{M}(\varepsilon) \equiv \varepsilon(\bm{r}) I  .
\end{equation}
Here $I$ denotes the identity operator. We denote by $(u^m,
\lambda^m)$ the $m$-th pair of eigenfunction and eigenvalue of
(\ref{eigen_cont}) and assume that these eigenpairs are numbered in
ascending order: $ 0 < \lambda^1 \le \lambda^2 \le \dots \le
\lambda^{\infty}$.

\subsection{The Optimization Problem}
\label{sec:theoptproblem}
The objective in photonic crystal design is to maximize the band gap
between two consecutive frequency modes. Due to the lack of
fundamental length scale in Maxwell's equations, it can be shown that
the magnitude of the band gap scales by a factor of $s$ when the
crystal is expanded by a factor of $1/s$. Therefore, it is more
meaningful to maximize the \textbf{gap-midgap} ratio instead of the
absolute band gap~\cite{joannopoulos2008pcm}. The gap-midgap ratio
between $\lambda^m$ and $\lambda^{m+1}$ is defined as
\begin{equation*}
J(\varepsilon(\bm{r})) = \frac{\inf_{\bm{k} \in \partial \mathcal{B}} \lambda^{m+1}(\varepsilon(\bm{r}),\bm{k}) - \sup_{\bm{k} \in \partial \mathcal{B}} \lambda^{m}(\varepsilon(\bm{r}),\bm{k})}{\inf_{\bm{k} \in \partial \mathcal{B}} \lambda^{m+1}(\varepsilon(\bm{r}),\bm{k}) + \sup_{\bm{k} \in \partial \mathcal{B}} \lambda^{m}(\varepsilon(\bm{r}),\bm{k})} ,
\end{equation*}
where $\partial \mathcal{B}$ represents the irreducible Brillouin zone
boundary; see Figure~\ref{SquLat_BriZone} for example.

A typical characterization of the dielectric function
$\varepsilon(\bm{r})$ is the distribution of two different
materials. Suppose that we are given two distinct materials with
dielectric constants $\varepsilon_{\min}$ and $\varepsilon_{\max}$
where $\varepsilon_{\min} < \varepsilon_{\max}$. We wish to find
arrangements of the materials within the unit cell $\Omega$ which
result in maximal gap-midgap ratio. To this end, we decompose the unit
cell $\Omega$ into $N_{\varepsilon}$ disjoint subcells $K_i, 1 \le i
\le N_{\varepsilon}$, such that $\Omega = \cup_{i=1}^{N_{\varepsilon}}
{K_i}$ and $K_i \cap K_j = \emptyset$ for $i \neq j$. Here we take
this subcell grid to be the same as the finite element triangulation
of the unit cell as we are going to discretize the continuous
eigenvalue problem by the finite element method. Our dielectric
function ${\varepsilon}(\bm{r})$ takes a unique value between
$\varepsilon_{\min}$ and $\varepsilon_{\max}$ on each subcell, namely,
${\varepsilon}(\bm{r}) = \varepsilon_i \in \mathbb{R}$ on $K_i$ and
$\varepsilon_{\min} \le \varepsilon_i \le
\varepsilon_{\max}$. However, due to the symmetry of square lattice,
we only need to define the dielectric function ${\varepsilon}(\bm{r})$
over part of the unit cell ($1/8$ of the unit cell). Hence, in
general, the dielectric function $\varepsilon(\bm{r})$ is discretized
into a finite dimensional vector $\bm{\varepsilon} = (\varepsilon_1,
\ldots, \varepsilon_{n_{\varepsilon}}) \in
\mathbb{R}^{n_{\varepsilon}}$ (with $n_{\epsilon} \le N_{\epsilon}$)
which resides in the following admissible region:
\begin{equation*}
\mathcal{Q}_{ad} \equiv \{ \bm{\varepsilon} = (\varepsilon_1, \ldots, \varepsilon_{n_{\varepsilon}}) \in \mathbb{R}^{n_{\varepsilon}} \ : \  \varepsilon_{\min} \le \varepsilon_i \le \varepsilon_{\max}, \; 1 \le i \le n_{\varepsilon}\} .
\end{equation*}
This region consists of piecewise-constant functions whose value on
every subcell varies between $\varepsilon_{\min}$ and
$\varepsilon_{\max}$.  Moreover, to render this problem computationally
tractable, we replace the irreducible Brillouin zone boundary
$\partial \mathcal{B}$ by a finite subset
\begin{equation*}
\mathcal{S}_{n_k} = \{\bm{k}_t \in \partial \mathcal{B}, \ 1 \le t \le n_k\},
\end{equation*}
where $\bm{k}_t, 1 \le t \le n_k,$ are wave vectors chosen along the
irreducible Brillouin zone boundary.  As a result, the band gap
optimization problem that maximizes the gap-midgap ratio between
$\lambda^m$ and $\lambda^{m+1}$ can be stated as follows:
\begin{equation}
\label{eqBSO_1}
\begin{split}
 & \max_{\bm{\varepsilon}} \; J^\ast(\bm{\varepsilon}) = \frac{\min_{\bm{k} \in  \mathcal{S}_{n_k}} \lambda^{m+1}(\bm{\varepsilon},\bm{k}) - \max_{\bm{k} \in  \mathcal{S}_{n_k}} \lambda^{m}(\bm{\varepsilon},\bm{k})}{\min_{\bm{k} \in  \mathcal{S}_{n_k}} \lambda^{m+1}(\bm{\varepsilon},\bm{k}) + \max_{\bm{k} \in  \mathcal{S}_{n_k}} \lambda^{m}(\bm{\varepsilon},\bm{k})} \\[1ex]
& \mbox{s.t.}\ \ \ \  \mathcal{A}(\bm \varepsilon,\bm{k}) {u}^j = \lambda^j \mathcal{M}(\bm \varepsilon) {u}^j, \quad j = m, m+1, \  \bm{k} \in  \mathcal{S}_{n_k} ,\\
& \qquad \ \ \varepsilon_{\min} \le \varepsilon_i \le \varepsilon_{\max}, \qquad \quad \ \ 1 \le i \le n_{\varepsilon} .
\end{split}
\end{equation}
In this problem a subtle difference between TE and TM polarizations
lies in the operators of the eigenvalue problem: $\mathcal{A}$ and
$\mathcal{M}$ take the form of either~(\ref{eqBSC_4}) for the TE case
or (\ref{eqBSC_4B}) for the TM case. In either case, note that the
eigenvalue problems embedded in (\ref{eqBSO_1}) must be addressed as
part of any computational strategy for the overall solution of
(\ref{eqBSO_1}).


\subsection{Discretization of the Eigenvalue Problem}

We consider here the finite element method to discretize the
continuous eigenvalue problem~(\ref{eigen_cont}). This produces the
following discrete eigenvalue problem
\begin{equation}
\label{eqBSC_1}
A_h(\bm{\varepsilon},\bm{k})u^j_{h} = \lambda^j_{h} M_h(\bm{\varepsilon})u^j_{h},
\qquad j = 1, \ldots, \mathcal{N}, \quad \bm{k} \in  \mathcal{S}_{n_k},
\end{equation}
where $A_h(\bm{\varepsilon},\bm{k}) \in \mathbb{C}^{\mathcal{N} \times
\mathcal{N}}$ is a Hermitian stiffness matrix and
$M_h(\bm{\varepsilon}) \in \mathbb{R}^{\mathcal{N} \times
\mathcal{N}}$ is a symmetric positive definite mass matrix. These
matrices are sparse and typically very large ($\mathcal{N} \gg 1$). We
consider the approximate eigenvalues in ascending order:
$\lambda^1_{h} \le \lambda^2_{h} \le \dots \le
\lambda^{\mathcal{N}}_{h}$.

It is important to note that the dependence of the above matrices on
the design parameter vector $\bm{\varepsilon}$ is different for the TE
and TM polarizations.  In the TE case, $A_h^{\rm TE}$ depends on
$\bm{\varepsilon}$ and $M_h^{\rm TE}$ does not, whereas in the TM case
$M_h^{\rm TM}$ depends on $\bm{\varepsilon}$ and $A_h^{\rm TM}$ does
not. More specifically, since $\varepsilon(\bm{r})$ is a
piecewise-constant function on $\Omega$, the
$\bm{\varepsilon}$-dependent matrices can be expressed as
\begin{equation}
\label{Affine}
A_h^{\rm TE}(\bm{\varepsilon},\bm{k}) = \sum_{i=1}^{n_{\varepsilon}} \frac{1}{{\varepsilon}_i} A^{\rm TE}_{h,i}(\bm{k}), \qquad
M_h^{\rm TM}(\bm{\varepsilon}) = \sum_{i=1}^{n_{\varepsilon}} {\varepsilon}_i M^{\rm TM}_{h,i} ,
\end{equation}
where the matrices $A^{\rm TE}_{h,i}(\bm{k})$ and $M^{\rm TM}_{h,i}, 1
\le i \le n_{\varepsilon}$ are independent of $\bm{\varepsilon}$. We
note that $A_h^{\rm TE}(\bm{\varepsilon},\bm{k})$ is linear with
respect to $1/{\varepsilon}_i$, $1 \le i \le n_{\varepsilon}$, while
$M_h^{\rm TM}(\bm{\varepsilon})$ is linear with respect to
${\varepsilon}_i$, $1 \le i \le n_{\varepsilon}$. The affine
expansion~(\ref{Affine}) is a direct consequence of the fact that we
use piecewise-constant approximation for the dielectric function
$\varepsilon(\bm{r})$.  (In the TE case, we will shortly change our
decision variables to $y_i = 1/{\varepsilon}_i$, $1 \le i \le
n_{\varepsilon}$, so as to render $A_h^{\rm TE}$ affine in the
variables $y_1, \ldots, y_{n_\varepsilon}$.)

After discretizing the eigenvalue problem (\ref{eigen_cont}) by the
finite element method, we obtain the following band gap optimization
problem:
\begin{equation}
\label{eqBSO_1B}
\begin{split}
 & \max_{\bm{\varepsilon}} J_h(\bm{\varepsilon}) = \frac{\min_{\bm{k} \in  \mathcal{S}_{n_k}} \lambda_h^{m+1}(\bm{\varepsilon},\bm{k}) - \max_{\bm{k} \in  \mathcal{S}_{n_k}} \lambda_h^{m}(\bm{\varepsilon},\bm{k})}{\max_{\bm{k} \in  \mathcal{S}_{n_k}} \lambda_h^{m+1}(\bm{\varepsilon},\bm{k}) + \max_{\bm{k} \in  \mathcal{S}_{n_k}} \lambda_h^{m}(\bm{\varepsilon},\bm{k})} \\[1ex]
& \mbox{s.t.}\ \ \ \  A_h(\bm \varepsilon,\bm{k}) {u}^j_h = \lambda^j_h \mathcal{M}_h(\bm \varepsilon) {u}^j_h, \quad j = m, m+1, \  \bm{k} \in  \mathcal{S}_{n_k} ,\\
& \qquad \ \ \varepsilon_{\min} \le \varepsilon_i \le \varepsilon_{\max}, \qquad \quad \quad \ 1 \le i \le n_{\varepsilon} .
\end{split}
\end{equation}
Unfortunately, this optimization problem is non-convex; furthermore it
suffers from lack of regularity at the optimum.  The reason for this
is that the eigenvalues $\lambda^m_h$ and $\lambda^{m+1}_h$ are
typically not smooth functions of $\bm{\varepsilon}$ at points of
multiplicity, and multiple eigenvalues at the optimum are typical of
structures with symmetry.  As a consequence, the gradient of the
objective function $J(\bm{\varepsilon})$ with respect to
$\bm{\varepsilon}$ is not well-defined at points of eigenvalue
multiplicity, and thus gradient-based descent methods often run into
serious numerical difficulties and convergence problems.

\section{Band Structure Optimization}
\label{sec_BSO}
In this section we describe our approach to solve the band gap
optimization problem based on a subspace method and semidefinite
programming (SDP). In our approach, we first reformulate the original
problem as an optimization problem in which we aim to maximize the
band gap obtained by restriction of the operator to two orthogonal
subspaces. The first subspace consists of eigenfunctions associated to
eigenvalues below the band gap, and the second subspace consists of
eigenfunctions whose eigenvalues are above the band gap.  In this way,
the eigenvalues are no longer explicitly present in the formulation,
and eigenvalue multiplicity no longer leads to lack of regularity. The
reformulated optimization problem is exact but non-convex and
large-scale. To reduce the problem size, we truncate the
high-dimensional subspaces to only a few eigenfunctions below and
above the band gap~\cite{cances2007fac,pau:046704}, thereby
obtaining a new small-scale yet non-convex optimization
problem. Finally, we keep the subspaces fixed at a given decision
parameter vector to obtain a convex semidefinite optimization problem
for which SDP solution methods can be efficiently applied.

\subsection{Reformulation of the Band Gap Optimization Problem using Subspaces}

We first define two additional decision variables:
\begin{equation*}
\lambda_h^u := \min_{\bm{k} \in  \mathcal{S}_{n_k}} \lambda_h^{m+1}(\bm{\varepsilon},\bm{k}) \ , \qquad \lambda_h^{\ell} := \max_{\bm{k} \in  \mathcal{S}_{n_k}} \lambda_h^{m}(\bm{\varepsilon},\bm{k}) \ ,
\end{equation*}
and then rewrite the original problem (\ref{eqBSO_1B}) as
\begin{equation}
\label{FormP0}
\begin{array}{llll}
P_0: &  \underset{\bm \varepsilon, \lambda_h^u, \lambda_h^{\ell}}{\rm max}  & \ \ \ \ \ \ \displaystyle  \frac{\lambda_h^u - \lambda_h^{\ell}}{\lambda_h^u + \lambda_h^{\ell}} \\ [1ex]   \\[1.5ex]
& \mbox{ s.t. }  & \lambda_{h}^m(\bm{\varepsilon},\bm{k}) \le \lambda_h^{\ell} \ , \ \lambda_h^u  \le
\lambda^{m+1}_h(\bm{\varepsilon},\bm{k}), \quad &  \forall \bm{k} \in \mathcal{S}_{n_k}, \\[1ex]
& & A_h(\bm{\varepsilon},\bm{k})  {u}_h^m   =   \lambda^{m}_h  M_h(\bm{\varepsilon}) {u}^{m}_h, \quad & \forall \bm{k} \in \mathcal{S}_{n_k} , \\[1ex]
& & A_h(\bm{\varepsilon},\bm{k})  {u}_h^{m+1}   =   \lambda^{m+1}_h  M_h(\bm{\varepsilon}) {u}^{m+1}_h, \quad & \forall \bm{k} \in \mathcal{S}_{n_k} , \\[1ex]
& & \displaystyle \varepsilon_{\min} \le \varepsilon_i \le \varepsilon_{\max},   \quad & i = 1,\ldots, n_\varepsilon, \\[1ex]
& & \lambda_h^u \ , \ \lambda_h^{\ell}  > 0 . &
\end{array}
\end{equation}

Next, we introduce the following matrices:
\begin{equation*}
\Phi^{\bm{\varepsilon}}(\bm{k}) := [\Phi^{\bm{\varepsilon}}_{\ell}(\bm{k}) \ | \ \Phi^{\bm{\varepsilon}}_{u}(\bm{k})] := [u_h^1(\bm{\varepsilon},\bm{k}) \, \ldots  \, u_h^m(\bm{\varepsilon},\bm{k}) \ | \ u_h^{m+1}(\bm{\varepsilon},\bm{k}) \, \ldots \, u_h^{\mathcal{N}}(\bm{\varepsilon},\bm{k})] ,
\end{equation*}
where $\Phi^{\bm{\varepsilon}}_{\ell}(\bm{k})$ and
$\Phi^{\bm{\varepsilon}}_{u}(\bm{k})$ consist of the first $m$ eigenvectors and
the remaining $\mathcal{N}-m$ eigenvectors, respectively, of the eigenvalue problem:
\begin{equation*}
A_h(\bm{\varepsilon},\bm{k})  {u}_h^j   =   \lambda^{j}_h  M_h(\bm{\varepsilon}) {u}^{j}_h, \qquad 1 \le j \le \mathcal{N}.
\end{equation*}
We will also denote the subspaces spanned by the eigenvectors of $\Phi^{\bm{\varepsilon}}_{\ell}(\bm{k})$ and $\Phi^{\bm{\varepsilon}}_{u}(\bm{k})$ as $\rg(\Phi^{\bm{\varepsilon}}_{\ell}(\bm{k}))$ and $\rg(\Phi^{\bm{\varepsilon}}_{u}(\bm{k}))$, respectively.


The first three sets of constraints in (\ref{FormP0}) can be represented exactly as
\begin{equation*}
\begin{split}
\Phi^{\bm{\varepsilon}\ast}_{\ell}(\bm{k})[A_h(\bm{\varepsilon},\bm{k}) - \lambda_h^{\ell}  M_h(\bm{\varepsilon}) ] \Phi^{\bm{\varepsilon}}_{\ell}(\bm{k}) &\preceq 0, \qquad \forall \bm{k} \in \mathcal{S}_{n_k} \\
\Phi^{\bm{\varepsilon}\ast}_{u}(\bm{k}) [A_h(\bm{\varepsilon},\bm{k})- \lambda_h^{u} M_h(\bm{\varepsilon})] \Phi^{\bm{\varepsilon}}_{u}(\bm{k}) &\succeq  0, \qquad  \forall \bm{k} \in \mathcal{S}_{n_k},
\end{split}
\end{equation*}
where ``$\succeq$'' is the L\"{o}wner partial ordering on symmetric
matrices, i.e., $A \succeq B$ if and only if $A-B$ is positive
semidefinite.  We therefore obtain the following equivalent
optimization problem:
\begin{equation}
\label{FormP1}
\begin{array}{llll}
P_1: &  \underset{\bm \varepsilon, \lambda_h^u, \lambda_h^{\ell}}{\rm max}  & \ \ \ \ \ \ \displaystyle  \frac{\lambda_h^u - \lambda_h^{\ell}}{\lambda_h^u + \lambda_h^{\ell}} \\ [1ex]   \\[1.5ex]
& \mbox{ s.t. }  & \Phi^{\bm{\varepsilon}\ast}_{\ell}(\bm{k})[A_h(\bm{\varepsilon},\bm{k}) - \lambda_h^{\ell}  M_h(\bm{\varepsilon}) ] \Phi^{\bm{\varepsilon}}_{\ell}(\bm{k}) \preceq 0, \quad & \forall \bm{k} \in \mathcal{S}_{n_k} , \\[1ex]
& & \Phi^{\bm{\varepsilon}\ast}_{u}(\bm{k}) [A_h(\bm{\varepsilon},\bm{k})- \lambda_h^{u} M_h(\bm{\varepsilon})] \Phi^{\bm{\varepsilon}}_{u}(\bm{k}) \succeq 0, \quad & \forall \bm{k} \in \mathcal{S}_{n_k} , \\[1ex]
& & \displaystyle \varepsilon_{\min} \le \varepsilon_i \le \varepsilon_{\max},   \quad & i = 1,\ldots, n_\varepsilon, \\[1ex]
& & \lambda_h^u \ , \ \lambda_h^{\ell}  > 0 . &
\end{array}
\end{equation}
Although the reformulation $P_1$ is exact, there is however a subtle
difference in the interpretation of $P_0$ and $P_1$: $P_0$ can be
viewed as maximizing the gap-midgap ratio between the two eigenvalues
$\lambda_h^{m}$ and $\lambda_h^{m+1}$; whereas $P_1$ can be viewed as
maximizing the gap-midgap ratio between the two subspaces
$\rg(\Phi^{\bm{\varepsilon}}_{\ell}(\bm{k}))$ and
$\rg(\Phi^{\bm{\varepsilon}}_{u}(\bm{k}))$.  The latter viewpoint allows us
to develop an efficient subspace approximation method for solving the
band gap optimization problem as discussed below.

\subsection{Subspace Approximation and Reduction}

Let us assume that we are \emph{given} a parameter vector
$\hat{\bm{\varepsilon}}$.  We then introduce the associated
matrices
\begin{equation*}
\Phi^{\hat{\bm{\varepsilon}}}(\bm{k}) := [\Phi^{\hat{\bm{\varepsilon}}}_{\ell}(\bm{k}) \ | \ \Phi^{\hat{\bm{\varepsilon}}}_{u}(\bm{k})] = [u_h^1(\hat{\bm{\varepsilon}},\bm{k}) \, \ldots  \, u_h^m(\hat{\bm{\varepsilon}},\bm{k}) \ | \ u_h^{m+1}(\hat{\bm{\varepsilon}},\bm{k}) \, \ldots \, u_h^{\mathcal{N}}(\hat{\bm{\varepsilon}},\bm{k})] \ ,
\end{equation*}
where $\Phi^{\hat{\bm{\varepsilon}}}_{\ell}(\bm{k})$ and $\Phi^{\hat{\bm{\varepsilon}}}_{u}(\bm{k})$ consist of the first $m$ eigenvectors and the remaining $\mathcal{N}-m$ eigenvectors, respectively, of the eigenvalue problem
\begin{equation*}
A_h(\hat{\bm{\varepsilon}},\bm{k})  {u}_h^j   =
\lambda^{j}_h  M_h(\hat{\bm{\varepsilon}}) {u}^{j}_h, \qquad 1 \le j \le \mathcal{N} .
\end{equation*}
Under the presumption that
$\rg(\Phi_{\ell}^{\hat{\bm{\varepsilon}}}(\bm{k}))$ and
$\rg(\Phi_{u}^{\hat{\bm{\varepsilon}}}(\bm{k}))$ are reasonable
approximations of $\rg(\Phi_{\ell}^{{\bm{\varepsilon}}}(\bm{k}))$ and
$\rg(\Phi_{u}^{{\bm{\varepsilon}}}(\bm{k}))$ for $\bm{\varepsilon}$ near
$\hat{\bm{\varepsilon}}$, we replace
$\Phi^{\bm{\varepsilon}}_{\ell}(\bm{k})$ with
$\Phi^{\hat{\bm{\varepsilon}}}_{\ell}(\bm{k})$ and
$\Phi^{\bm{\varepsilon}}_{u}(\bm{k})$ with
$\Phi^{\hat{\bm{\varepsilon}}}_{u}(\bm{k})$ to obtain
\begin{equation}
\label{FormP2}
\begin{array}{llll}
P_2^{\hat{\bm{\varepsilon}}}: &  \underset{\bm \varepsilon, \lambda_h^u, \lambda_h^{\ell}}{\rm max}  & \ \ \ \ \ \ \displaystyle  \frac{\lambda_h^u - \lambda_h^{\ell}}{\lambda_h^u + \lambda_h^{\ell}} \\ [1ex]   \\[1.5ex]
& \mbox{ s.t. }  & \Phi^{\hat{\bm{\varepsilon}\ast}}_{\ell}(\bm{k})[A_h(\bm{\varepsilon},\bm{k}) - \lambda_h^{\ell}  M_h(\bm{\varepsilon}) ] \Phi^{\hat{\bm{\varepsilon}}}_{\ell}(\bm{k}) \preceq 0, \quad & \forall \bm{k} \in \mathcal{S}_{n_k} , \\[1ex]
& & \Phi^{\hat{\bm{\varepsilon}\ast}}_{u}(\bm{k}) [A_h(\bm{\varepsilon},\bm{k})- \lambda_h^{u} M_h(\bm{\varepsilon})] \Phi^{\hat{\bm{\varepsilon}}}_{u}(\bm{k}) \succeq 0, \quad & \forall \bm{k} \in \mathcal{S}_{n_k} , \\[1ex]
& & \displaystyle \varepsilon_{\min} \le \varepsilon_i \le \varepsilon_{\max},   \quad & i = 1,\ldots, n_\varepsilon, \\[1ex]
& & \lambda_h^u \ , \ \lambda_h^{\ell}  > 0 . &
\end{array}
\end{equation}
Note in $P_2^{\hat{\bm{\varepsilon}}}$ that the subspaces
$\rg(\Phi^{\hat{\bm{\varepsilon}}}_{\ell}(\bm{k}))$ and
$\rg(\Phi^{\hat{\bm{\varepsilon}}}_{u}(\bm{k}))$ are approximations of the
subspaces $\rg(\Phi^{{\bm{\varepsilon}}}_{\ell}(\bm{k}))$ and
$\rg(\Phi^{{\bm{\varepsilon}}}_{u}(\bm{k}))$ and are no longer functions of
the decision variable vector $\bm \varepsilon$.

Note also that the semidefinite inclusions in
$P_2^{\hat{\bm{\varepsilon}}}$ are large-scale, i.e., the rank of
either the first or second inclusion is at least ${\cal N}/2$, for
each $\bm k \in {\cal S}_{n_k}$, and $\cal N$ will typically be quite
large.  In order to reduce the size of the inclusions, we reduce the
dimensions of the subspaces by considering only the ``important''
eigenvectors among $u_h^1(\bm{\varepsilon},\bm{k}) \, \ldots \,
u_h^m(\bm{\varepsilon},\bm{k}) , u_h^{m+1}(\bm{\varepsilon},\bm{k}) \,
\ldots \, u_h^{\mathcal{N}}(\bm{\varepsilon},\bm{k})$, namely those
$a_{\bm{k}}$ eigenvectors whose eigenvalues lie below but nearest to
$\lambda^m_h(\bm \varepsilon, \bm k)$ and those $b_{\bm{k}}$
eigenvectors whose eigenvalues lie above but nearest to
$\lambda^{m+1}_h(\bm \varepsilon, \bm k)$, for small values of
$a_{\bm{k}}$, $b_{\bm{k}}$, typically chosen in the range between $2$
and $5$, for each $\bm k \in {\cal S}_{n_k}$.  This yields reduced
matrices
\begin{equation*}
\Phi^{\hat{\bm{\varepsilon}}}_{a_{\bm{k}}+b_{\bm{k}}}(\bm{k}) :=
[\Phi^{\hat{\bm{\varepsilon}}}_{a_{\bm{k}}}(\bm{k}) \ | \ \Phi^{\hat{\bm{\varepsilon}}}_{b_{\bm{k}}}(\bm{k})] =
[u_h^{m-a_{\bm{k}}+1}(\hat{\bm{\varepsilon}},\bm{k}) \, \ldots  \, u_h^m(\hat{\bm{\varepsilon}},\bm{k}) \ | \ u_h^{m+1}(\hat{\bm{\varepsilon}},\bm{k}) \, \ldots \, u_h^{m+b_{\bm{k}}}(\hat{\bm{\varepsilon}},\bm{k})] .
\end{equation*}
Substituting $\Phi^{\hat{\bm{\varepsilon}}}_{a_{\bm{k}}}(\bm{k})$ in place of $\Phi^{\hat{\bm{\varepsilon}\ast}}_{\ell}(\bm{k})$ and $\Phi^{\hat{\bm{\varepsilon}}}_{b_{\bm{k}}}(\bm{k})$ in place of $\Phi^{\hat{\bm{\varepsilon}\ast}}_{u}(\bm{k})$ in  the formulation
$P_2^{\hat{\bm{\varepsilon}}}$ yields the following reduced
optimization formulation:
\begin{equation}
\label{FormP3}
\begin{array}{llll}
P_3^{\hat{\bm{\varepsilon}}}: &  \underset{\bm \varepsilon, \lambda_h^u, \lambda_h^{\ell}}{\rm max}  & \ \ \ \ \ \ \displaystyle  \frac{\lambda_h^u - \lambda_h^{\ell}}{\lambda_h^u + \lambda_h^{\ell}} \\ [1ex]   \\[1.5ex]
& \mbox{ s.t. }  & \Phi^{\hat{\bm{\varepsilon}\ast}}_{a_{\bm{k}}}(\bm{k})[A_h(\bm{\varepsilon},\bm{k}) - \lambda_h^{\ell}  M_h(\bm{\varepsilon}) ] \Phi^{\hat{\bm{\varepsilon}}}_{a_{\bm{k}}}(\bm{k}) \preceq 0, \quad & \forall \bm{k} \in \mathcal{S}_{n_k} , \\[1ex]
& & \Phi^{\hat{\bm{\varepsilon}\ast}}_{b_{\bm{k}}}(\bm{k}) [A_h(\bm{\varepsilon},\bm{k})- \lambda_h^{u} M_h(\bm{\varepsilon})] \Phi^{\hat{\bm{\varepsilon}}}_{b_{\bm{k}}}(\bm{k}) \succeq 0, \quad & \forall \bm{k} \in \mathcal{S}_{n_k} , \\[1ex]
& & \displaystyle \varepsilon_{\min} \le \varepsilon_i \le \varepsilon_{\max},   \quad & i = 1,\ldots, n_\varepsilon, \\[1ex]
& & \lambda_h^u \ , \ \lambda_h^{\ell}  > 0 . &
\end{array}
\end{equation}

In this way the formulation $P_3^{\hat{\bm{\varepsilon}}}$ seeks to
model only the anticipated ``active'' eigenvalue constraints, in exact
extension of active-set methods in nonlinear optimization.  The
integers $a_{\bm{k}}, b_{\bm{k}}$ are determined indirectly through
user-defined parameters $r_l>0$, and $r_u>0$, where we retain only
those eigenvectors whose eigenvalues are within $100r_l\%$ beneath
$\lambda^m_h(\hat{\bm \varepsilon}, \bm k)$ or whose eigenvalues are
within $100r_u\%$ above $\lambda^{m+1}_h(\hat{\bm \varepsilon}, \bm
k)$.  This translates to choosing $a_{\bm{k}}, b_{\bm{k}} \in
\mathbb{N}_{+}$ as the smallest integers that satisfy
\begin{eqnarray*}
\frac{\lambda_h^{m}(\hat{\bm{\varepsilon}},\bm k) - \lambda_h^{m-a_{\bm{k}}+1}(\hat{\bm{\varepsilon}},\bm k)}{\lambda_h^{m}(\hat{\bm{\varepsilon}},\bm k)} &\leq& r_l \ \le \frac{\lambda_h^{m}(\hat{\bm{\varepsilon}},\bm k) - \lambda_h^{m-a_{\bm{k}}}(\hat{\bm{\varepsilon}},\bm k)}{\lambda_h^{m}(\hat{\bm{\varepsilon}},\bm k)}, \\
\frac{\lambda_h^{m+b_{\bm{k}}}(\hat{\bm{\varepsilon}},\bm k) - \lambda_h^{m+1}(\hat{\bm{\varepsilon}},\bm k)}{\lambda_h^{m+1}(\hat{\bm{\varepsilon}},\bm k)} &\leq& r_u \ \le \frac{\lambda_h^{m+b_{\bm{k}}+1}(\hat{\bm{\varepsilon}},\bm k) - \lambda_h^{m+1}(\hat{\bm{\varepsilon}},\bm k)}{\lambda_h^{m+1}(\hat{\bm{\varepsilon}},\bm k)} \ .
\end{eqnarray*}
The dimensions of the resulting subspaces
$\rg(\Phi^{\hat{\bm{y}}}_{a_{\bm{k}}}(\bm{k}))$ and
$\rg(\Phi^{\hat{\bm{y}}}_{b_{\bm{k}}}(\bm{k}))$ are typically very small
($a_{\bm{k}},b_{\bm{k}} \sim 2,\ldots,5$).  Furthermore, the subspaces
are well-spanned by including all relevant eigenvectors corresponding
to those eigenvalues with multiplicity at or near the current min/max
values.

We observe that $P_3^{\hat{\bm{\varepsilon}}}$ has significantly smaller
semidefinite inclusions than if the {\em full} subspaces were used.
Also, the subspaces are kept {\em fixed} at $\hat{\bm{\varepsilon}}$
in order to reduce the nonlinearity of the underlying problem.
Furthermore, we show below that for the TE and TM polarizations that
$P_3^{\hat{\bm{\varepsilon}}}$ can be easily re-formulated as a linear
fractional semidefinite program, and hence is solvable using modern
interior-point methods.

\subsection{Fractional SDP Formulations for TE and TM Polarizations}

We now show that by a simple change of variables for each of the TE
and TM polarizations, problem $P_3^{\hat{\bm{\varepsilon}}}$ can be
converted to a linear fractional semidefinite program and hence can be
further converted to a linear semidefinite program.

\paragraph{TE polarization.} We introduce the following new decision variable
notation for convenience:
\begin{equation*}
\bm{y} := (y_1,y_2,\ldots,y_{n_y}) :=
(1/\varepsilon_1,\ldots,1/\varepsilon_{n_{\varepsilon}}, \lambda^{\ell}_h, \lambda^{u}_h) \ ,
\end{equation*}
and set $y_{\min} = 1/\varepsilon_{\max}$ and $y_{\max} =
1/\varepsilon_{\min}$.  We also amend our notation to write various
functional dependencies on $\bm y$ instead of $\bm \varepsilon$ such
as $\Phi^{\hat{\bm{y}}}_{\ell}(\bm{k})$, etc.  Utilizing
(\ref{Affine}), we re-write $P_3^{\hat{\bm{\varepsilon}}}$ for the TE
polarization as
\begin{equation}
\label{FormTE}
\begin{array}{llll}
P_{\rm TE}^{\hat{\bm{y}}}: &  \underset{\bm{y}}{\rm max}  & \ \ \ \ \ \ \displaystyle  \frac{y_{n_y} - y_{n_y-1}}{y_{n_y} + y_{n_y-1}} \\ [1ex]   \\[1.5ex]
& \mbox{ s.t. } & \Phi^{\hat{\bm{y}}\ast}_{a_{\bm{k}}}(\bm{k}) \left[\sum_{i=1}^{n_y-2} y_{i} A_{h,i}^{\rm TE}(\bm{k}) - y_{n_y-1}  M_h^{\rm TE} \right] \Phi^{\hat{\bm{y}}}_{a_{\bm{k}}}(\bm{k}) \preceq 0, \quad & \forall \bm{k} \in \mathcal{S}_{n_k} , \\[1ex]
& & \Phi^{\hat{\bm{y}}\ast}_{b_{\bm{k}}}(\bm{k}) \left[\sum_{i=1}^{n_y-2} y_{i} A_{h,i}^{\rm TE}(\bm{k}) - y_{n_y} M_h^{\rm TE}\right] \Phi^{\hat{\bm{y}}}_{b_{\bm{k}}}(\bm{k}) \succeq 0, \quad & \forall \bm{k} \in \mathcal{S}_{n_k} , \\[1ex]
& & \displaystyle y_{\min} \le y_i \le y_{\max},   \quad & i = 1,\ldots, n_y -2, \\[1ex]
& & y_{n_y-1} \ , \   y_{n_y}  > 0 . &
\end{array}
\end{equation}
We note that the objective function is a linear fractional expression
and the constraint functions are linear functions of the variables
$\bm{y}$.  Therefore $P_{\rm TE}^{\hat{\bm{y}}}$ is a linear
fractional SDP.  Using a standard homogenization \cite{charnes1962plf,craven1973dfl}, a
linear fractional SDP can be converted to a linear
SDP.\footnote{Indeed, for notational simplicity consider a linear
fractional optimization problem of the form $\max_x \frac{c^Tx}{d^Tx}$
subject to $b-Ax \in K_1$, $x \in K_2$, where $d^Tx >0$ for all
feasible $x$ and $K_1$, $K_2$ are convex cones.  Then this problem is
equivalent to the problem $\max_{w,\theta} c^Tw$ subject to $b\theta
-Aw \in K_1$, $w \in K_2$, $d^Tw=1$, $\theta \ge 0$, under the
elementary transformations $x \leftarrow (w/\theta)$ and $(w,\theta)
\leftarrow (x/d^Tx, 1/d^Tx)$, see \cite{charnes1962plf,craven1973dfl}.}

\paragraph{TM polarization.}
We introduce slightly different decision variable notation for convenience:
\begin{equation*}
\bm{z} := (z_1,z_2,\ldots,z_{n_z}) :=
(\varepsilon_1,\ldots,\varepsilon_{n_{\varepsilon}}, 1/\lambda^{\ell}_h, 1/\lambda^{u}_h),
\end{equation*}
and set $z_{\min} = \varepsilon_{\min}$ and $z_{\max} = \varepsilon_{\max}$.  Similar to the TE case, we amend our notation to write various functional dependencies on $\bm z$ instead of $\bm \varepsilon$ such as $\Phi^{\hat{\bm{z}}}_{\ell}(\bm{k})$, etc.  Noting that
$$\frac{\lambda_h^u - \lambda_h^{\ell}}{\lambda_h^u + \lambda_h^{\ell}} = \frac{z_{n_z-1} - z_{n_z}}{z_{n_z-1} + z_{n_z}} \ , $$ utilizing (\ref{Affine}), and multiplying the semidefinite inclusions of (\ref{FormP3}) by $z_{n_z-1}$ and $z_{n_z}$, respectively, we re-write $P_3^{\hat{\bm{\varepsilon}}}$ for the TM polarization as
\begin{equation}
\label{FormTM}
\begin{array}{llll}
P_{\rm TM}^{\hat{\bm{z}}}: &  \underset{\bm{z}}{\rm max}  & \ \ \ \ \ \ \displaystyle  \frac{z_{n_z-1} - z_{n_z}}{z_{n_z-1} + z_{n_z}} \\ [1ex]   \\[1.5ex]
& \mbox{ s.t. } & \Phi^{\hat{\bm{z}}\ast}_{a_{\bm{k}}}(\bm{k}) \left[z_{n_z-1}A_{h}^{\rm TM}(\bm{k}) - \sum_{i=1}^{n_z-2} z_{i}   M_{h,i}^{\rm TM} \right] \Phi^{\hat{\bm{z}}}_{a_{\bm{k}}}(\bm{k}) \preceq 0, \quad & \forall \bm{k} \in \mathcal{S}_{n_k} , \\[1ex]
& & \Phi^{\hat{\bm{z}}\ast}_{b_{\bm{k}}}(\bm{k}) \left[z_{n_z}A_{h}^{\rm TM}(\bm{k}) - \sum_{i=1}^{n_z-2} z_{i}   M_{h,i}^{\rm TM} \right] \Phi^{\hat{\bm{z}}}_{b_{\bm{k}}}(\bm{k}) \succeq 0, \quad & \forall \bm{k} \in \mathcal{S}_{n_k} , \\[1ex]
& & \displaystyle z_{\min} \le z_i \le z_{\max},   \quad & i = 1,\ldots, n_z -2, \\[1ex]
& & z_{n_z-1} \ , \   z_{n_z}  > 0 . &
\end{array}
\end{equation}
Here again the objective function is a linear fractional form and the
constraint functions are linear functions of the variables $\bm{z}$.
Therefore $P_{\rm TM}^{\hat{\bm{z}}}$ is a linear fractional SDP with
format similar to that of $P_{\rm TE}^{\hat{\bm{y}}}$.

Since both $P_{\rm TE}^{\hat{\bm{y}}}$ and $P_{\rm TM}^{\hat{\bm{z}}}$
are linear fractional semidefinite programs, they can be solved very
efficiently by using modern interior point methods. Here we use the
SDPT3 software~\cite{tutuncu2003ssq} for this task.

\subsection{Main Algorithm}

We summarize our numerical approach for solving the band gap
optimization problem of the TE polarization in the following
table. Essentially the same algorithm (with the modifications
described in the previous section) is used to solve the band gap
optimization problem of the TM polarization.
\begin{table}[hbt]
\begin{center}
\begin{tabular}{|l|}
\hline
\\[-2ex]
\hspace{4.7cm} {\bf Implementation Steps} \\[1ex]
\hline
\\[-2ex]
{\bf Step 1.}  Start with an initial guess $\bm{y}^0$ and an error tolerance $\epsilon_{\rm tol}$, and set $\hat{\bm{y}} := \bm{y}^0$.
\\[1ex]
{\bf Step 2.}  For each wave vector $\bm{k} \in \mathcal{S}_{n_k}$, do:
\\[1ex]
\hspace{1.45cm} Determine the subspace dimensions $a_{\bm{k}}$ and $b_{\bm{k}}$.
\\[1ex]
\hspace{1.45cm} Compute the matrices $\Phi^{\hat{\bm{y}}}_{a_{\bm{k}}}(\bm{k})$ and $\Phi^{\hat{\bm{y}}}_{b_{\bm{k}}}(\bm{k})$.
\\[1ex]
{\bf Step 3.} Form the semidefinite program $P_{\rm TE}^{\hat{\bm{y}}}$.
\\[1ex]
{\bf Step 4.}  Solve $P_{\rm TE}^{\hat{\bm{y}}}$ for an optimal solution ${\bm y}^\ast$.
\\[1ex]
{\bf Step 5.}  If $\|{\bm y}^\ast - \hat{\bm y}\| \le \epsilon_{\rm tol}$, stop and return the optimal solution ${\bm y}^\ast$.
\\[1ex]
\hspace{1.45cm} Else update $\hat{\bm y} \leftarrow {\bm y}^\ast$ and go to {\bf Step 2.}
\\[1ex]
\hline
\end{tabular}\caption{Main algorithm for solving the band gap optimization problem.}
\label{Algorithm01}
\end{center}
\end{table}

\section{Results and Discussions}
\label{sec_OR}
\subsection{Model Setup}
We consider a two-dimensional photonic crystal confined in the
computational domain of a unit cell of the square lattice, and with
square domain $\Omega \equiv [-1,1] \times [-1,1]$. The domain
$\Omega$ is decomposed into a uniform quadrilateral (in particular, we
use square elements for the square lattice) grid of dimensions $64
\times 64$, which yields a mesh size of $h=1/32$ and $4096$ linear
square elements.

\begin{figure}[h]
\centering
\includegraphics[scale = 0.6]{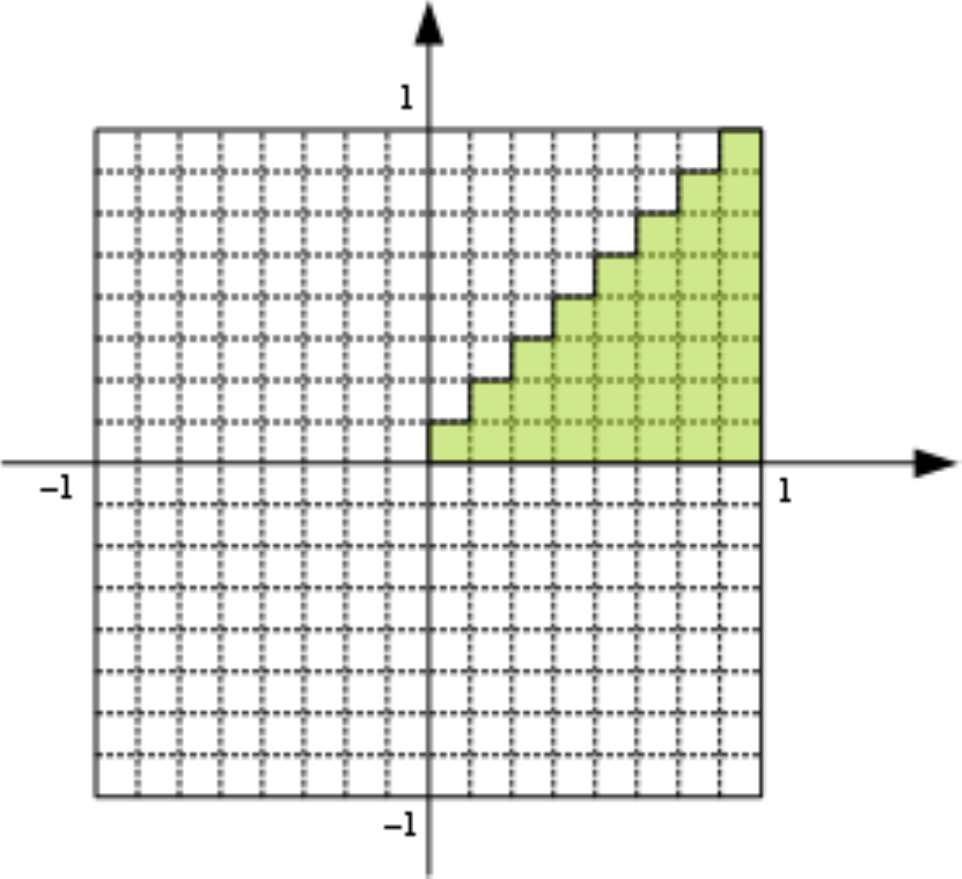}
\caption{An illustration of a coarse mesh ($16\times 16$) and dielectric function for the
square lattice. The shaded cells indicate the decision variables relating the dielectric material ($\varepsilon_i, i=1,2,\ldots,n_{\varepsilon}$). Note that the actual computational mesh ($64\times 64$) is finer than this one.}
\label{figOR_SquLat}
\end{figure}

The dielectric function $\bm{\varepsilon}$ is composed of two
materials with dielectric constants $\varepsilon_{\min} = 1$ (air) and
$\varepsilon_{\max} = 11.4$ (GaAs). As mentioned earlier in
Section~\ref{sec:theoptproblem}, the symmetry of the lattice can be
exploited to further reduce the dielectric function to be defined in
only $1/8$ of the computation domain. The number of decision variables
relating the dielectric material ($\varepsilon_i$, $i= 1, 2,\ldots,
n_{\varepsilon}$) is thus reduced to {$n_{\varepsilon} = (1+32)\times
32/2 = 528$}.  Figure~\ref{figOR_SquLat} shows an illustration of a
coarse mesh ($16\times 16$) and dielectric function for the square lattice to aid
visualization; note that the actual computational mesh ($64\times 64$) is finer
than this one. The shaded cells represent those modeled by
$\varepsilon$, and the rest are obtained through
symmetry. Furthermore, in this case, the irreducible Brillouin zone
$\mathcal{B}$ is the triangle shown in Figure~\ref{SquLat_BriZone},
with $n_k = 12$ $\bm{k}$-points taken along the boundary of this
region ($\partial \mathcal{B}$). Band diagrams plotted in the figures
below show the eigenvalues moving along the boundary of $\mathcal{B}$,
from $\Gamma$ to $\text{X}$ to $\text{M}$ and back to $\Gamma$.

\subsection{Choices of Parameters}
\subsubsection{Initial configuration}

\begin{figure}[hbt]
\centering
\subfigure[Initial crystal configuration $\#1$]{
\includegraphics[scale=0.35]{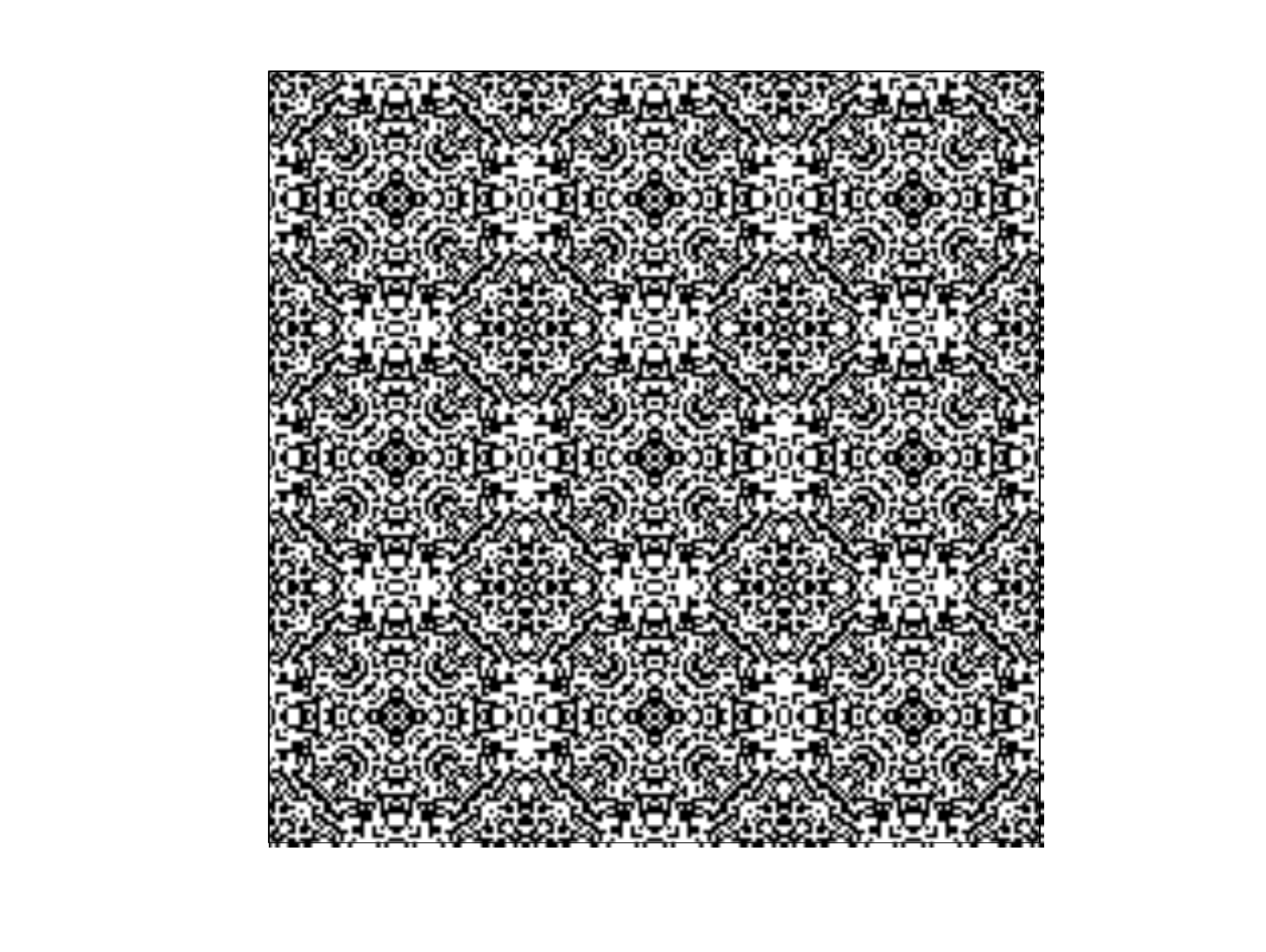}}
\subfigure[Optimized crystal structure $\#1$]{
\includegraphics[scale=0.35]{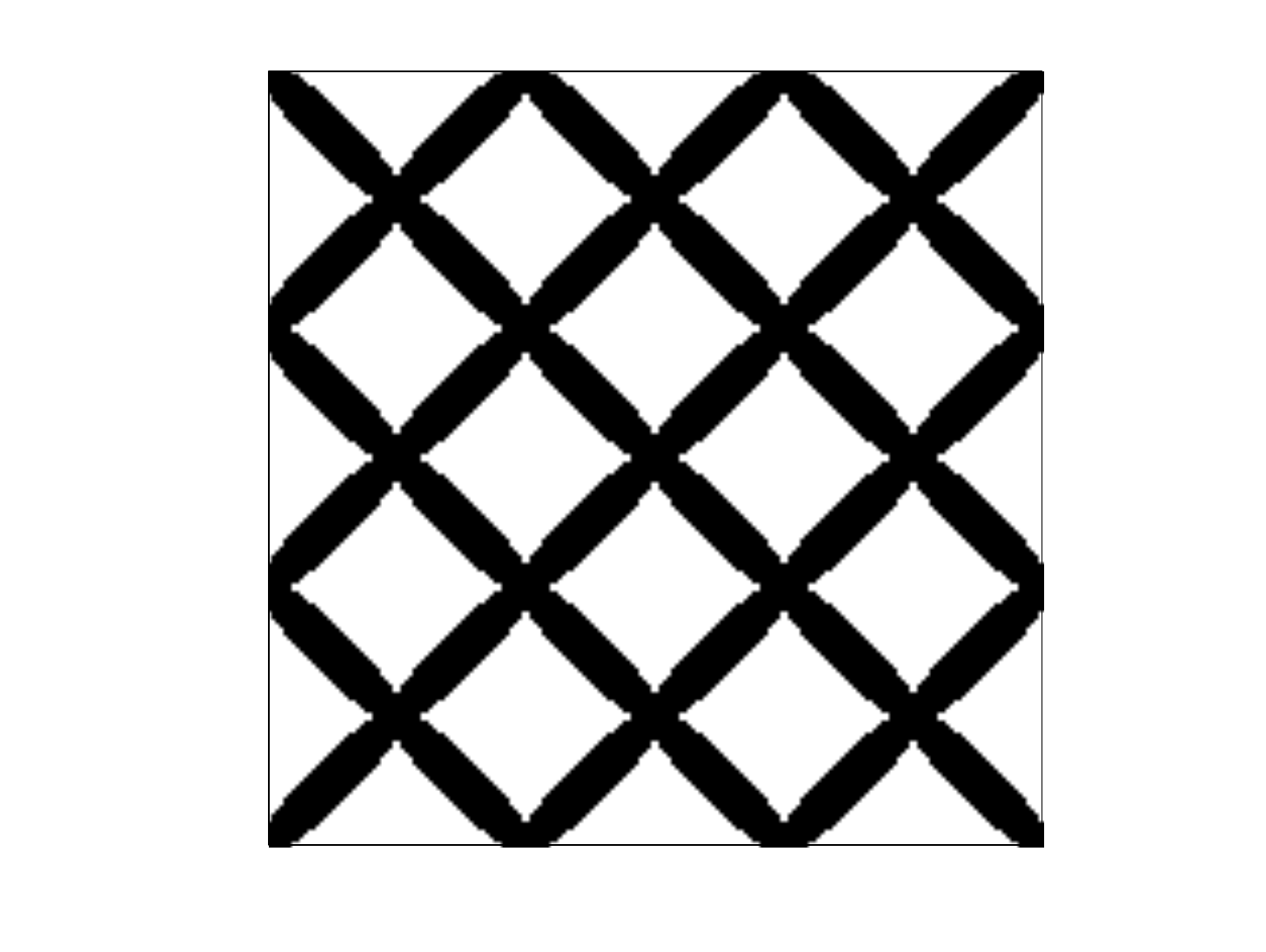}}
\subfigure[Optimized band structure $\#1$]{
\includegraphics[scale=0.35]{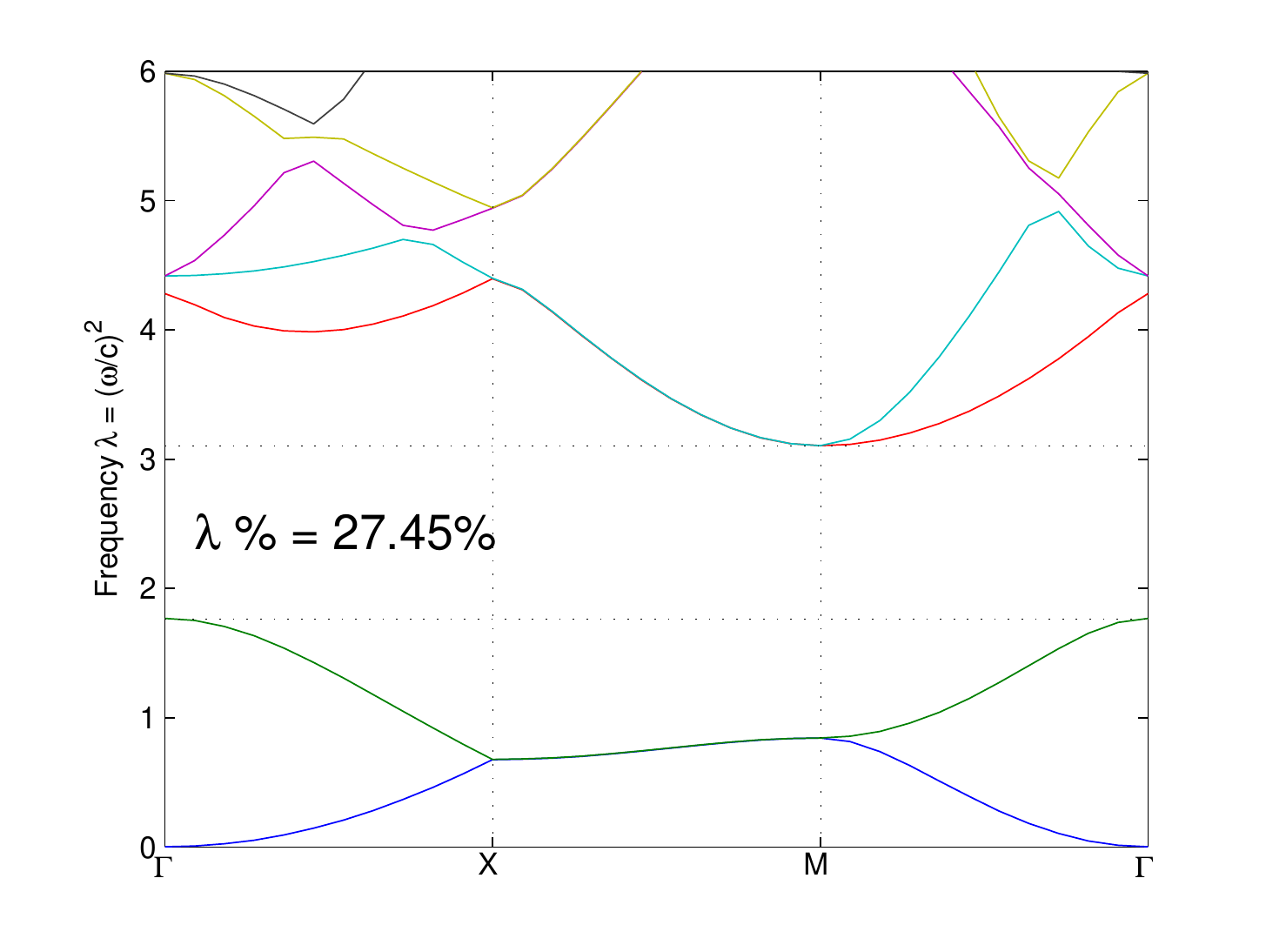}}
\subfigure[Initial crystal configuration $\#2$]{
\includegraphics[scale=0.35]{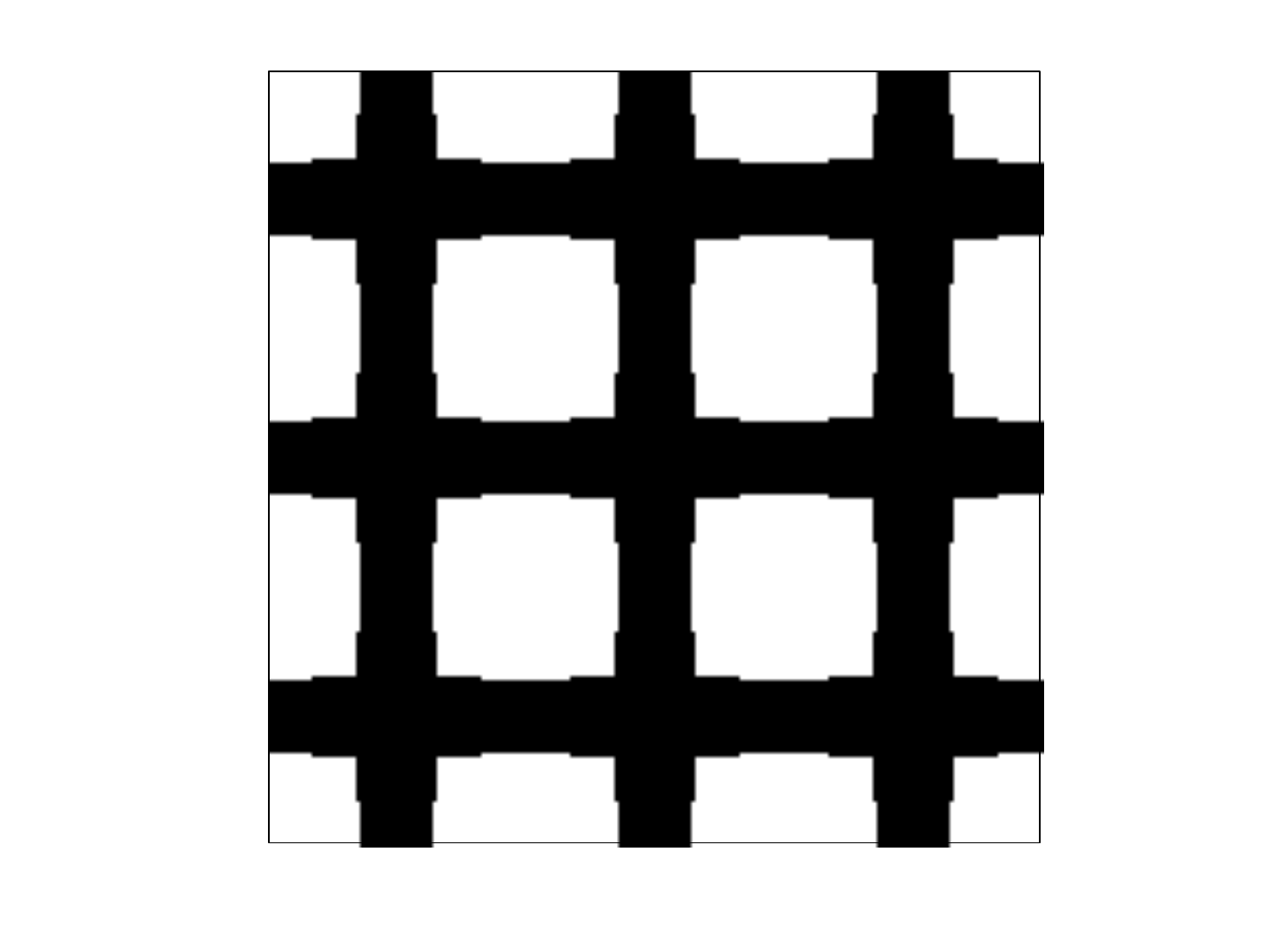}}
\subfigure[Optimized crystal structure $\#2$]{
\includegraphics[scale=0.35]{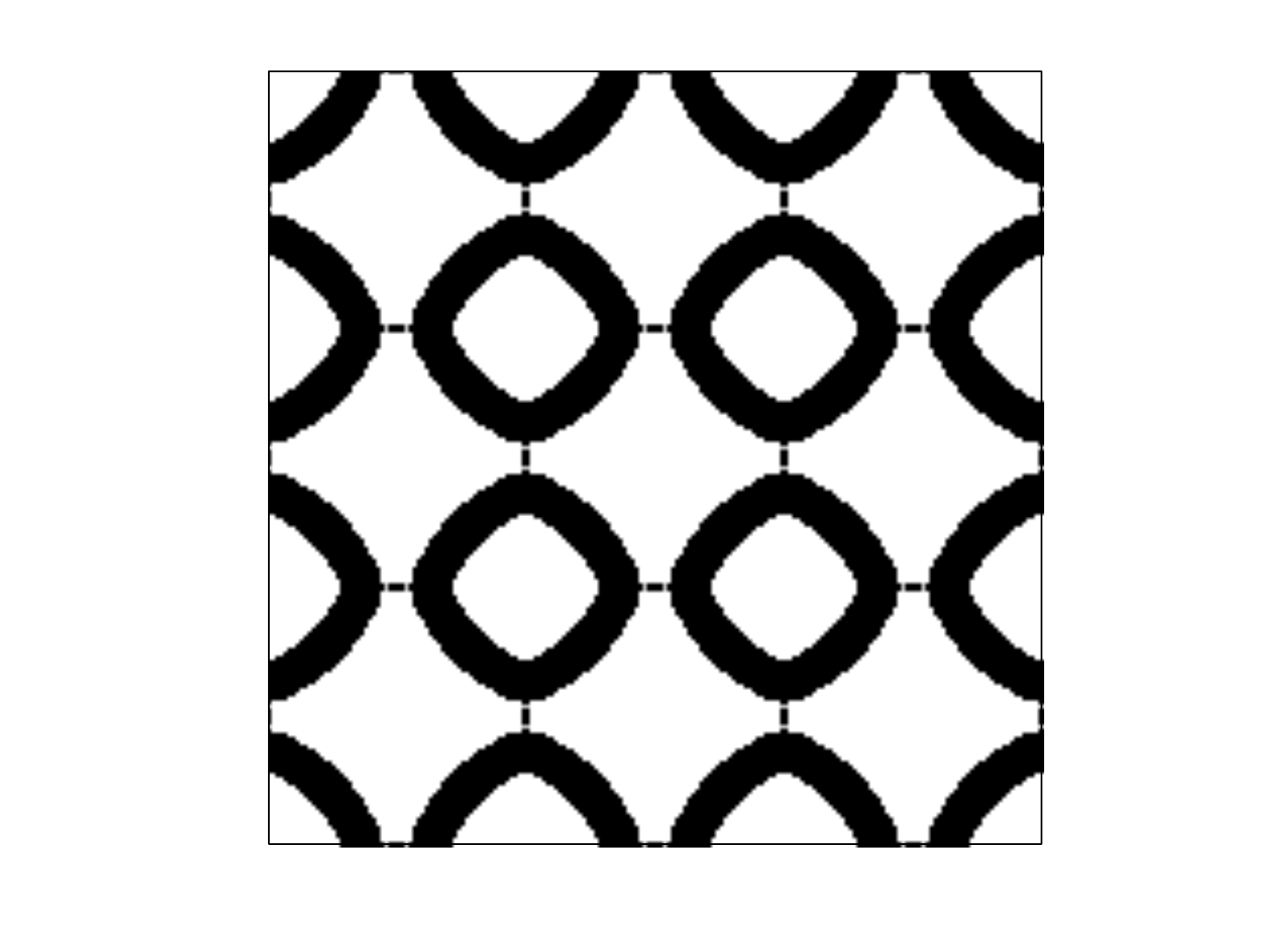}}
\subfigure[Optimized band structure $\#2$]{
\includegraphics[scale=0.35]{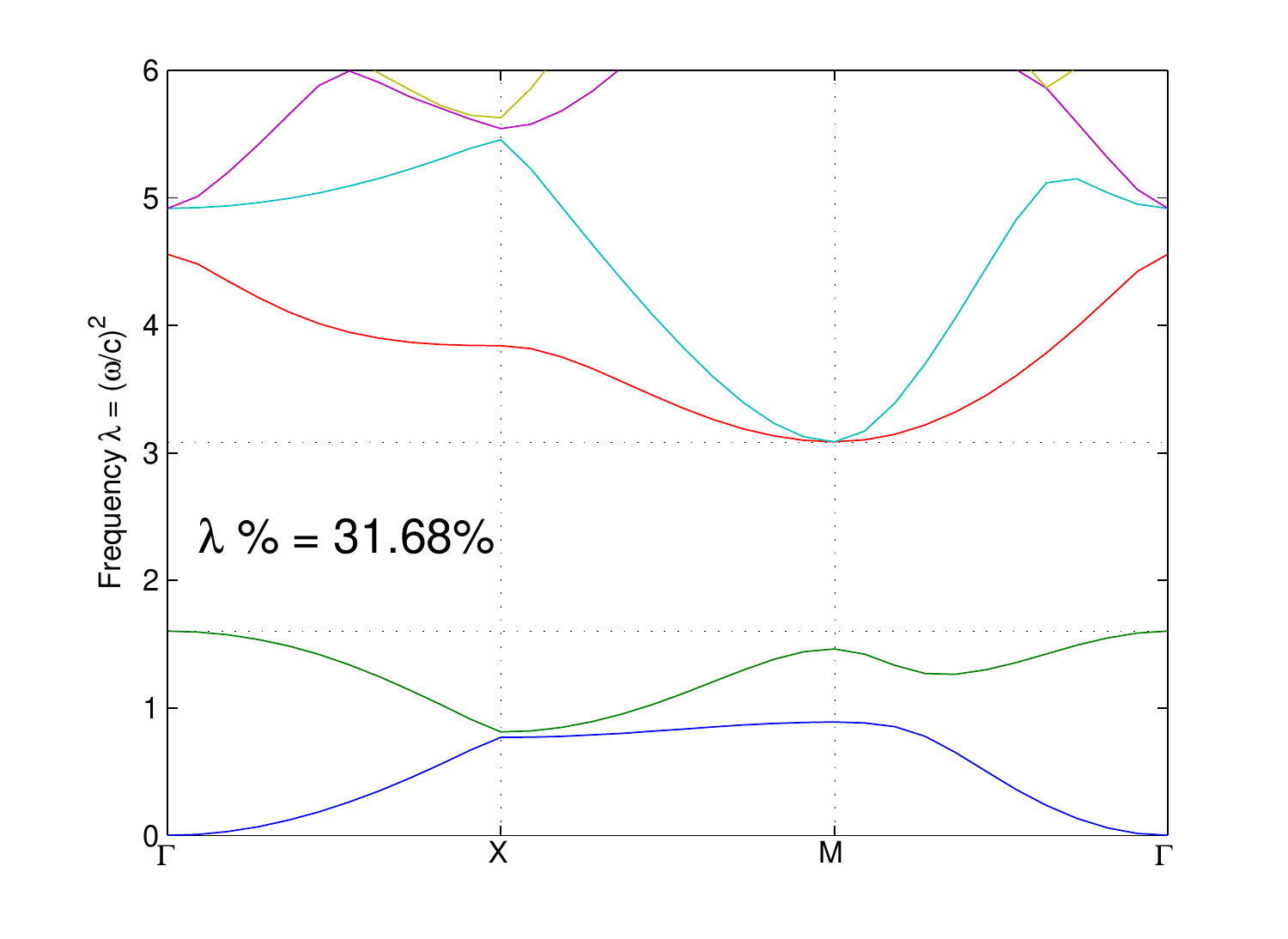}}
\caption{Two locally optimal band gaps between $\lambda^2_{\text{TE}}$ and $\lambda^3_{\text{TE}}$ in the square lattice}
\label{figOR_TElocal}
\end{figure}

\begin{figure}[hbt]
\centering
\subfigure[Initial crystal configuration $\#1$]{
\includegraphics[scale=0.35]{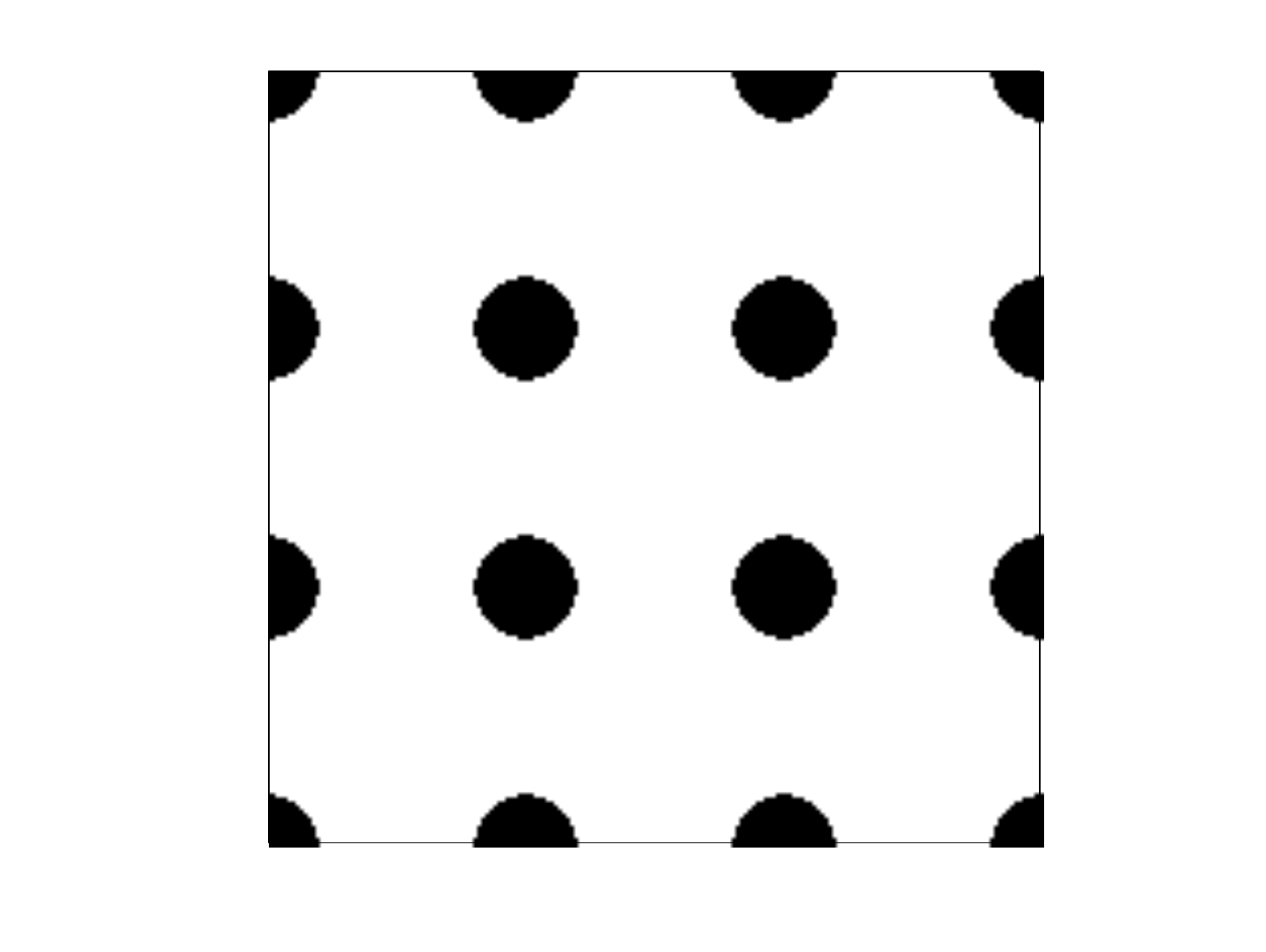}}
\subfigure[Optimized crystal structure $\#1$]{
\includegraphics[scale=0.35]{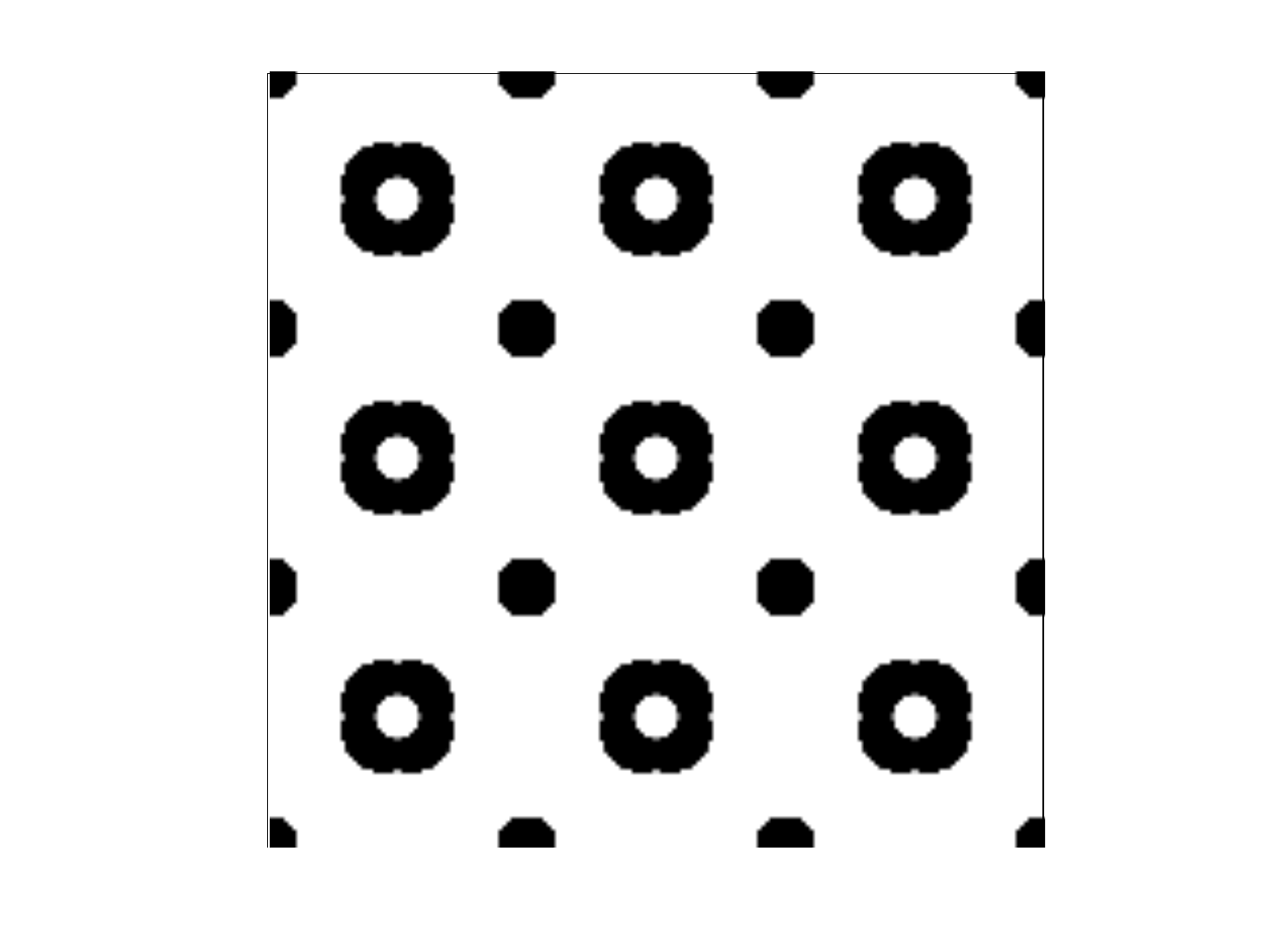}}
\subfigure[Optimized band structure $\#1$]{
\includegraphics[scale=0.35]{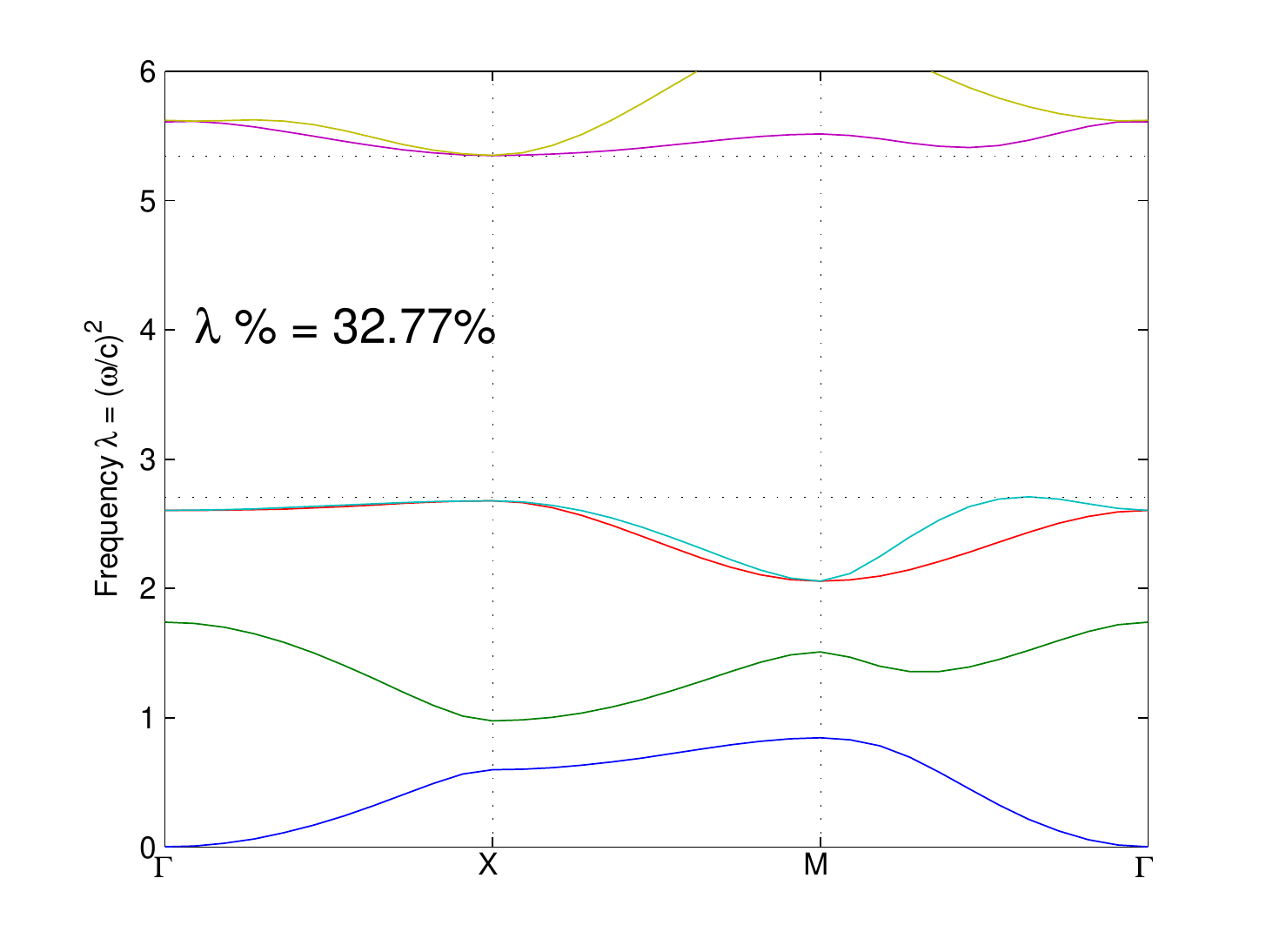}}
\subfigure[Initial crystal configuration $\#2$]{
\includegraphics[scale=0.35]{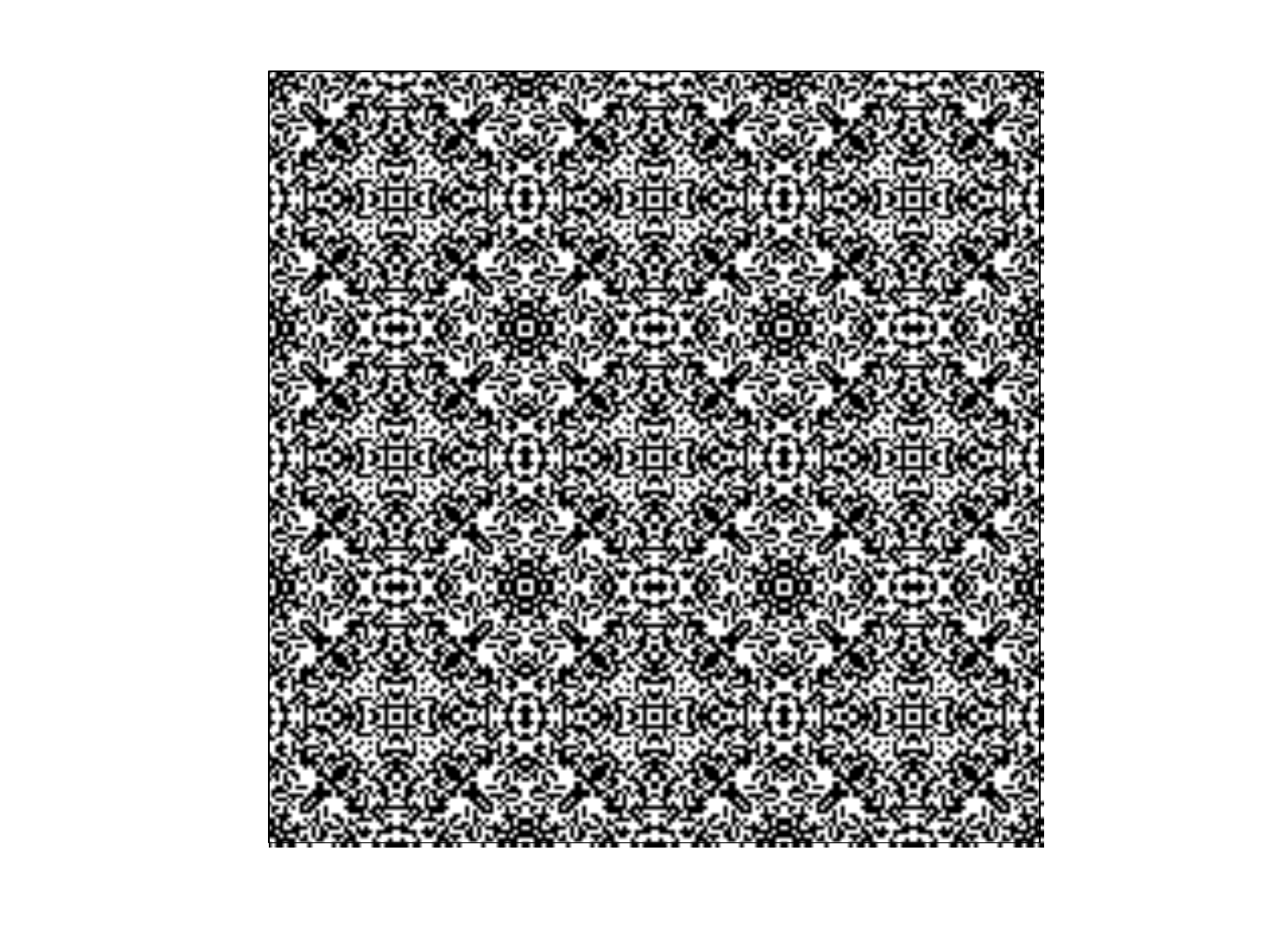}}
\subfigure[Optimized crystal structure $\#2$]{
\includegraphics[scale=0.35]{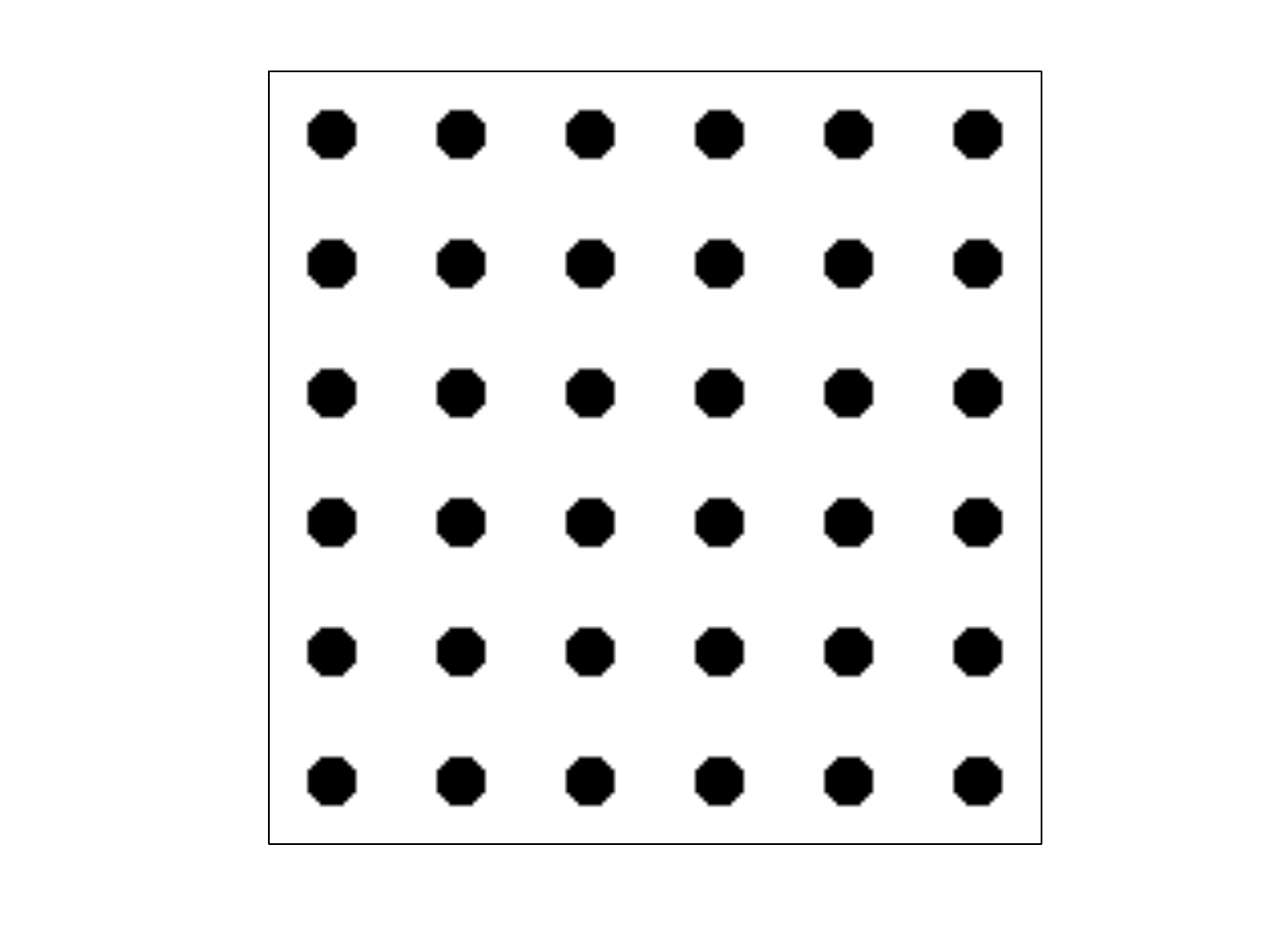}}
\subfigure[Optimized band structure $\#2$]{
\includegraphics[scale=0.35]{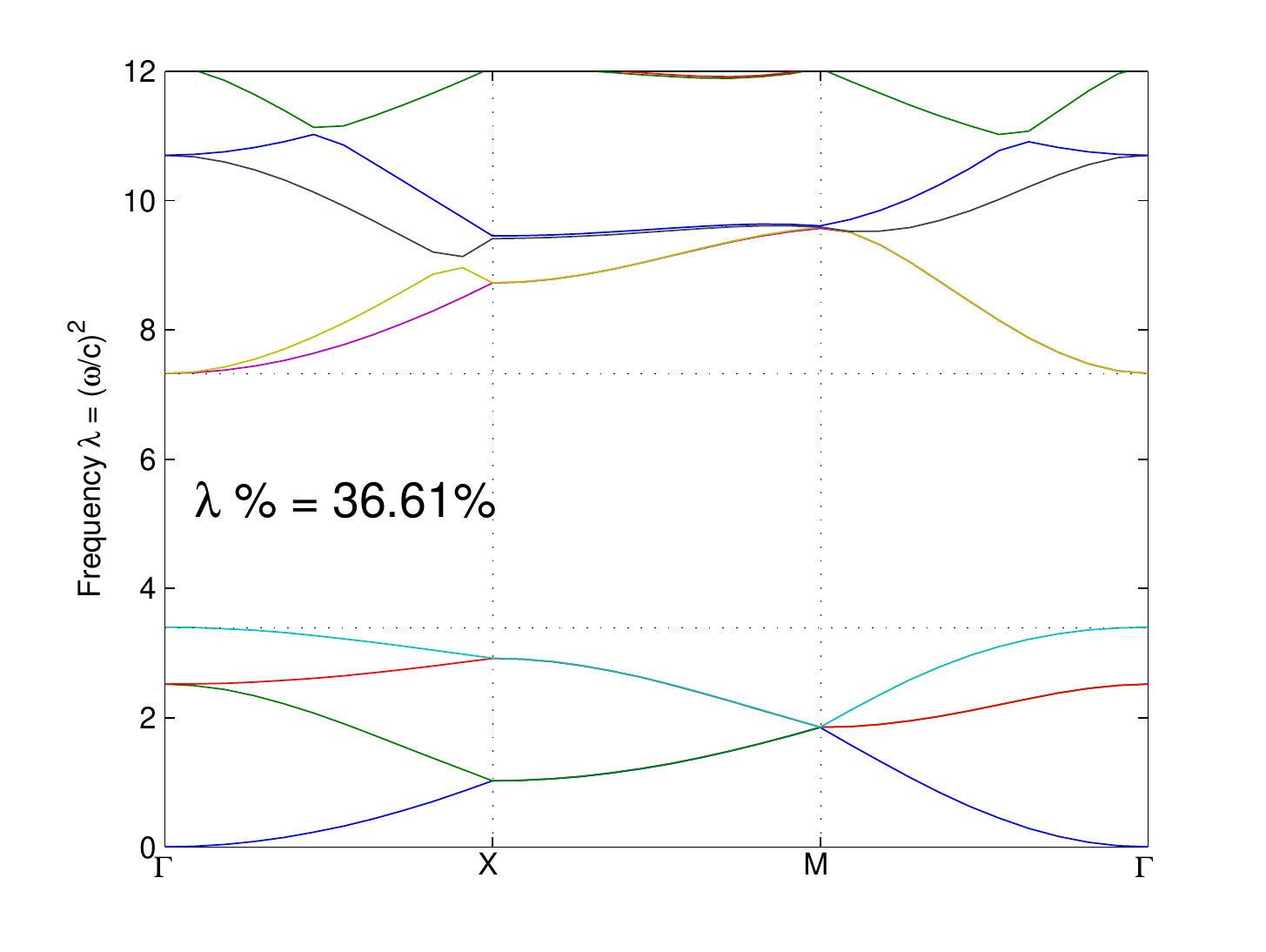}}
\caption{Two locally optimal band gaps between $\lambda^4_{\text{TM}}$ and $\lambda^5_{\text{TM}}$ in the square lattice}
\label{figOR_TMlocal}
\end{figure}

Because the underlying optimization problem may have many local optima, the performance of our
method can be sensitive to the choice of the initial values of the decision variables $\bm y^0$, which in turn depend on the initial
configuration $\bm \varepsilon^0$. Indeed, different initial configurations do lead to different local optima as shown in Figure~\ref{figOR_TElocal} for the second TE band gap and in Figure~\ref{figOR_TMlocal} for the fourth TM band gap. Therefore, the choice of the initial configuration is important. We examine here two different types of initial configurations: photonic crystals exhibiting band gaps at the low frequency spectrum and random distribution.

The well-known photonic crystals (e.g., dielectric rods in air -- Figure \ref{figOR_TMlocal}(a), air holes in dielectric material, orthogonal dielectric veins -- Figure \ref{figOR_TElocal}(d)) exhibit band-gap structures at the low frequency spectrum. Such a distribution seems to be a sensible choice for the initial configuration as it resembles various known optimal structures
\cite{burger87ipt}. When these well-known photonic crystals are used as the initial configuration, our method easily produces the band-gap structures at the low frequency mode (typically, the first three TE and TM modes).  On the other hand, maximizing the band gap at the high frequency mode (typically, above the first three TE and TM modes) tends to produce more complicated structures which are very different from the known photonic crystals mentioned above. As a result, when these photonic crystals are used as the initial configurations for maximizing the band gap at the high frequency mode, the obtained results are less satisfactory.


Random initial configurations such as Figures \ref{figOR_TElocal}(a) and \ref{figOR_TMlocal}(d)) have very high spatial variation and may thus be suitable for maximizing the band gap at the high frequency mode. Indeed, we observe that random distributions often yield larger band gaps (better results) than the known photonic crystals for the high frequency modes. Of course, the random initialization does not eliminate the possibility of multiple local optima intrinsic to the physical problem. In view of this effect, we use multiple random distributions to initialize our method. In particular, we start our main algorithm with a number of \emph{uniformly} random distributions as initial configurations to obtain the optimal structures in our numerical results discussed below.

\subsubsection{Subspace dimensions}

The dimensions of the subspaces $\rg(\Phi_{a_{\bm k}}^{\hat{\bm y}}(\bm
k))$ and $\rg(\Phi_{b_{\bm k}}^{\hat{\bm y}}(\bm k))$ are determined
indirectly by the parameters $r_l$ and $r_u$. A good choice of $r_l$
(and $r_u$) is one that returns $a_{\bm k} \ll \mathcal{N}$ (and
$b_{\bm k} \ll \mathcal{N}$), and at the same time includes the
``important'' eigenvectors to enhance convergence to an optimum. In our numerical experiments, we
choose $r_u = r_l = 0.1$ which in turn leads to the resulting subspace dimensions $a_{\bm k}$ and $b_{\bm k}$ in the range of
$[2,5]$. Moreover, we find that choosing larger values of $r_u$ and $r_l$ (e.g., $r_{l} = r_u = 0.2$), which in turn increases $a_{\bm k}$ and $b_{\bm k}$ and hence increases computational cost, does not yield fewer iterations than choosing $r_u=r_l=0.1$.

\subsection{Computational Cost}
With all the programs implemented in MATLAB and the computation
performed on a Linux PC with Dual Core AMD Opteron 270, $1.99$GHz, a
successful run of the algorithm can typically be done in $2$--$30$
minutes including $5$--$30$ outer iterations, i.e., passes of Steps 2-5
of the main algorithm in Table~\ref{Algorithm01}. An example of the
computational cost and outer iterations for different band gap
optimization is shown in Table~\ref{tabCompCost} as a general
illustration of our computational experience.

We point out that these numbers merely represent one set of
possibilities; variations in the numerical results are likely to occur
with different random initial configurations. Nevertheless, the
computation cost does serve as an indication of the general level of
difficulty of finding a solution in each problem. In general, lower
eigenvalue band gap optimization problems are easier to solve (at
least to local optima). Moreover, the table illustrates that TM
problems usually solve faster and require fewer outer iterations. This
latter observation is consistent with the result reported in
\cite{kao2005mbg}, and is possibly due to the high non-convexity of the
original TE optimization problem.

\begin{table}[hbt]
\begin{center}
{\footnotesize
\begin{tabular}{l||c|c|c|c|c|c|c|c|c|c}
\hline
\hline
&$\Delta \lambda^{TE}_{1,2}$ & $\Delta \lambda^{TE}_{2,3}$ & $\Delta \lambda^{TE}_{3,4}$ & $\Delta \lambda^{TE}_{4,5}$ & $\Delta \lambda^{TE}_{5,6}$ & $\Delta \lambda^{TE}_{6,7}$ & $\Delta \lambda^{TE}_{7,8}$ & $\Delta \lambda^{TE}_{8,9}$ & $\Delta \lambda^{TE}_{9,10}$ & $\Delta \lambda^{TE}_{10,11}$\\
\hline
Execution time (min)& $5.7$ & $2.5$ & $8.9$ & $20.4$ & $17.9$ & $20.5$ & $19.4$ & $27.3$ & $26.4$ & $25.8$\\
\hline
Outer Iterations & $11$ & $8$ & $29$ & $26$ & $18$ & $25$ & $15$ & $27$ & $19$ & $23$ \\
\hline
\hline
& $\Delta \lambda^{TM}_{1,2}$ & $\Delta \lambda^{TM}_{2,3}$ & $\Delta \lambda^{TM}_{3,4}$ & $\Delta \lambda^{TM}_{4,5}$ & $\Delta \lambda^{TM}_{5,6}$ & $\Delta \lambda^{TM}_{6,7}$ & $\Delta \lambda^{TM}_{7,8}$ & $\Delta \lambda^{TM}_{8,9}$ & $\Delta \lambda^{TM}_{9,10}$ & $\Delta \lambda^{TM}_{10,11}$\\
\hline
Execution time (min)& $1.8$ & $5.6$ & $3.5$ & $5.4$ & $11.7$ & $9.5$ & $10.8$ & $3.9$ & $11.2$ & $9.5$\\
\hline
Outer Iterations & $4$ & $9$ & $5$ & $7$ & $16$ & $9$ & $9$ & $9$ & $12$ & $10$ \\
\hline
\end{tabular}}
\caption{Example of computation time and the number of outer iterations of a successful run for optimizing various band gaps, for both TE and TM polarization.  Here $\Delta \lambda^{TE}_{i,i+1}$ denotes the gap-midgap ratio between the $i^{\rm th}$ and $(i+1)^{\rm st}$ eigenvalue for the TE polarization.}
\label{tabCompCost}
\end{center}
\end{table}

Before ending this section, we point out some possible ways to improve
the computational cost of our procedure. For the eigenvalue
calculation, it is probably helpful to apply a more efficient
eigensolver (we used MATLAB's {\tt eigs} function in the current
implementation). Another promising approach is to explore mesh
adaptivity and incorporate non-uniform meshing for the representation
of the dielectric function, as well as the eigenvalue
calculation. As further discussed in Section~\ref{conclusion}, it
should be possible to significantly reduce the number of decision
variables and computation cost with this approach.

\subsection{Optimal Structures}
\begin{figure}
\centering
\subfigure[]{
\includegraphics[scale=0.5]{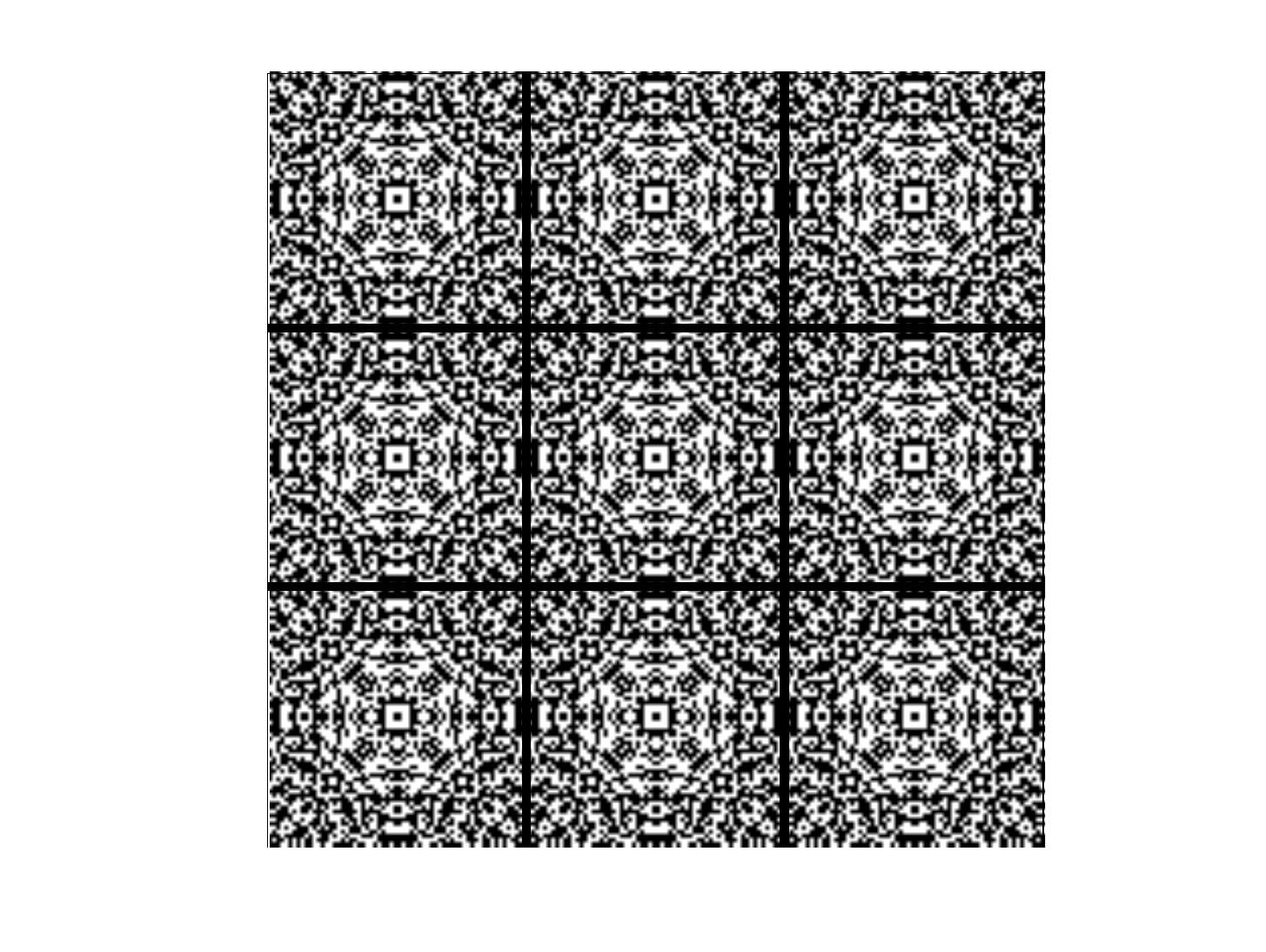}}
\subfigure[]{
\includegraphics[scale=0.5]{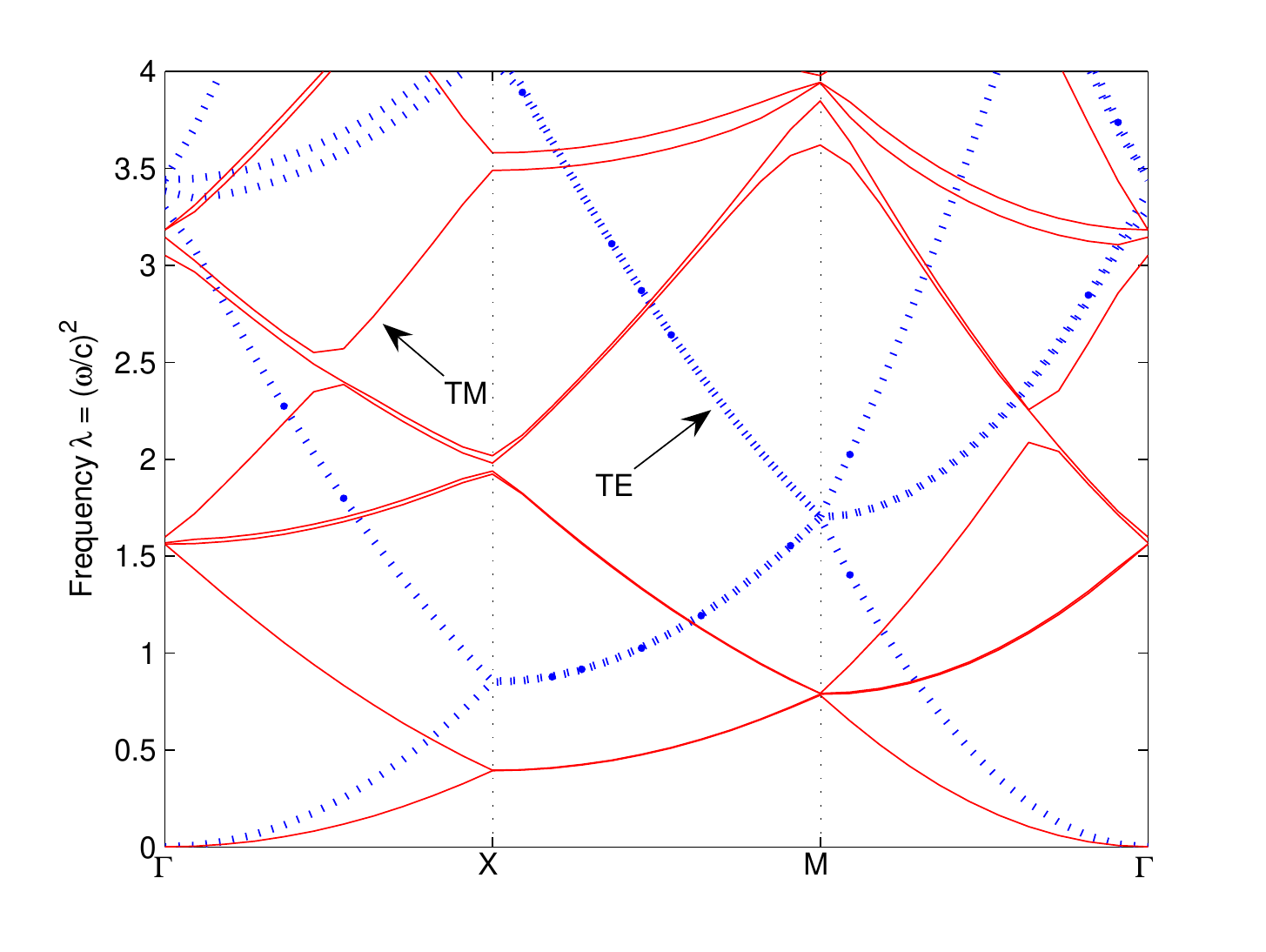}}
\caption{(a) Random starting structure with translation, rotation, and reflection symmetry, $3\times 3$ unit cells in square lattice. (b) Band structure before optimization. }
\label{figOR_beforeOpt}
\end{figure}

\begin{figure}
\centering
\subfigure[]{
\includegraphics[scale=0.35]{Figures/Sq_before_geometry}}
\subfigure[]{
\includegraphics[scale=0.35]{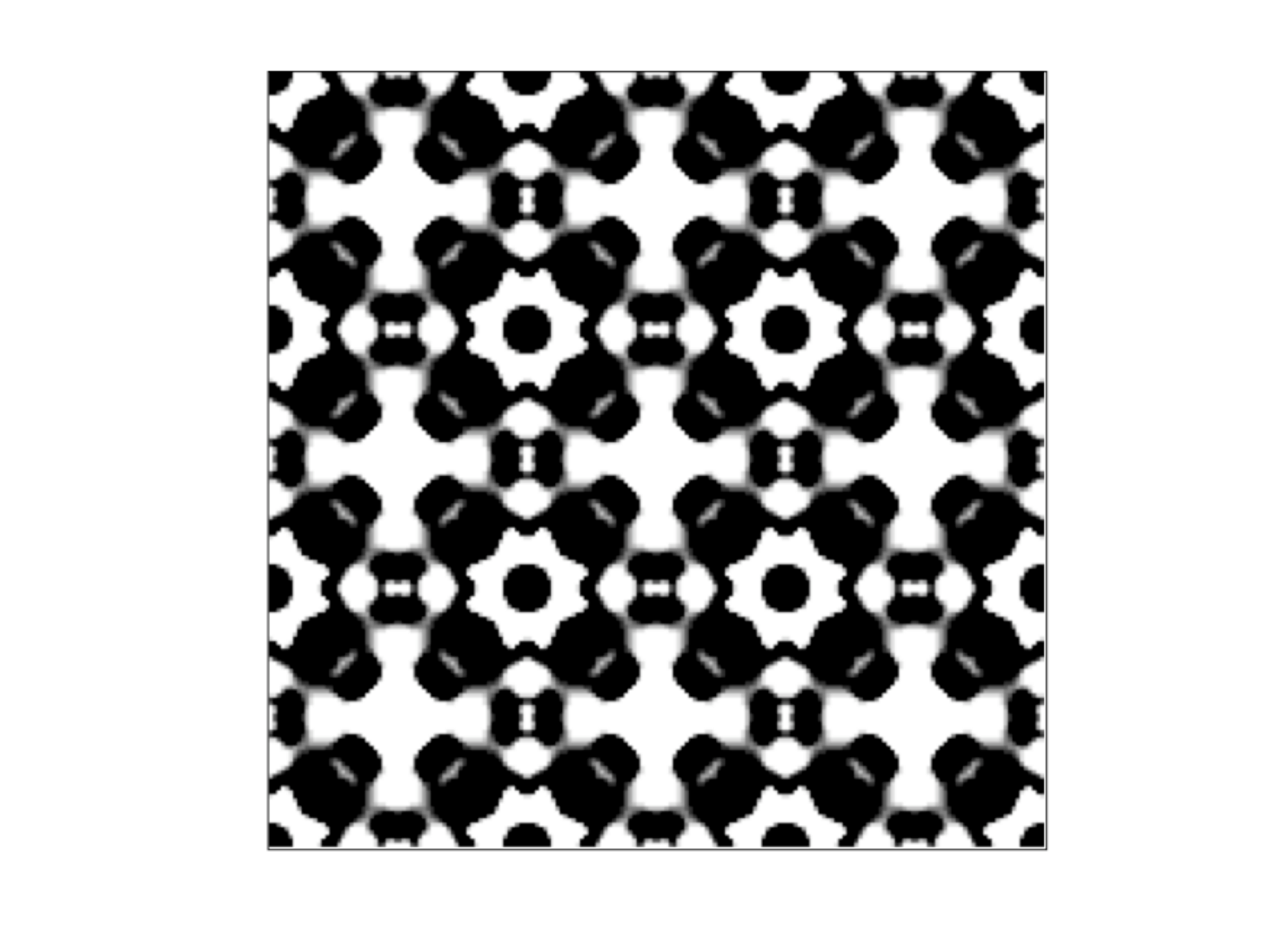}}
\subfigure[]{
\includegraphics[scale=0.35]{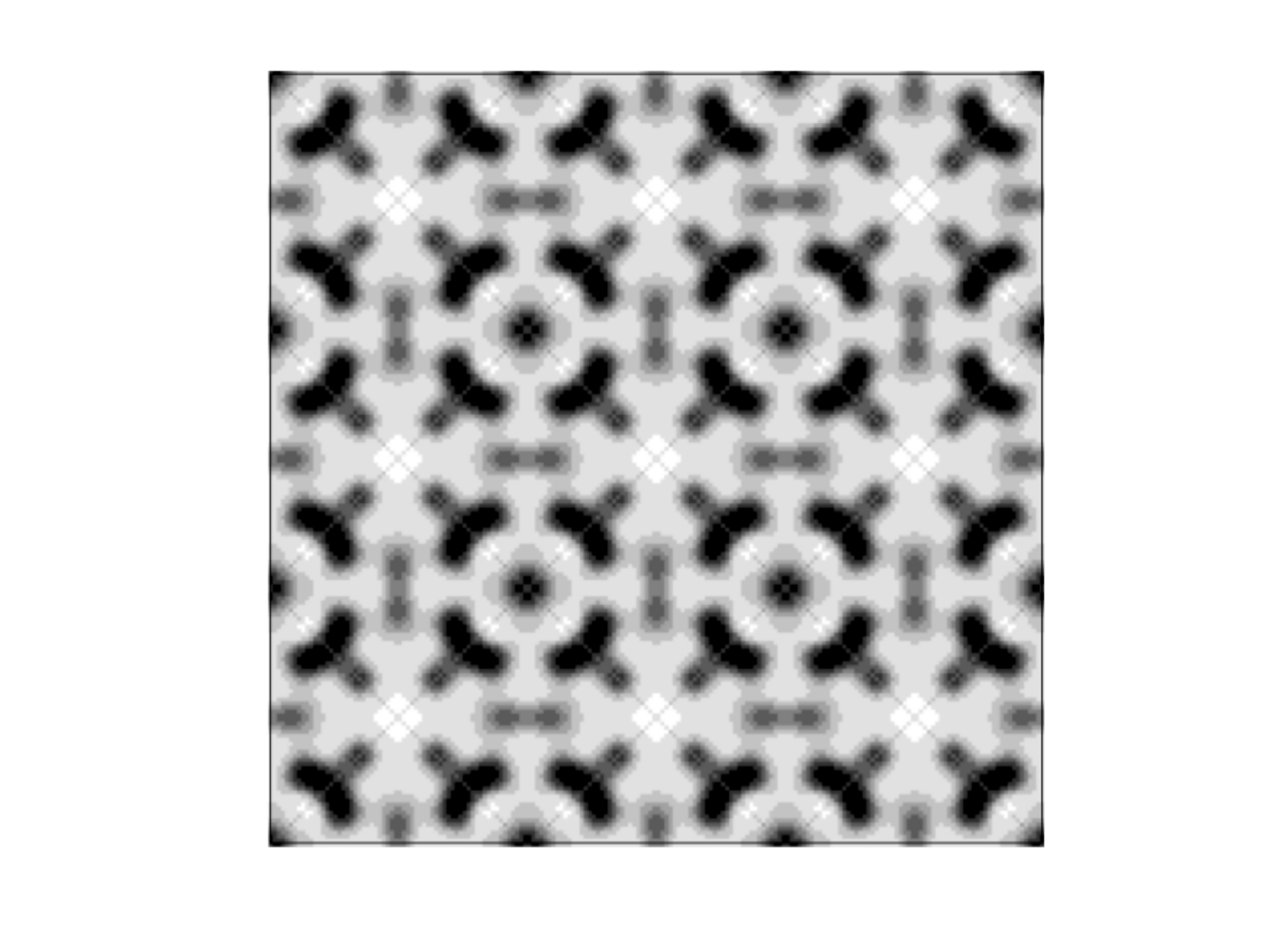}}
\subfigure[]{
\includegraphics[scale=0.35]{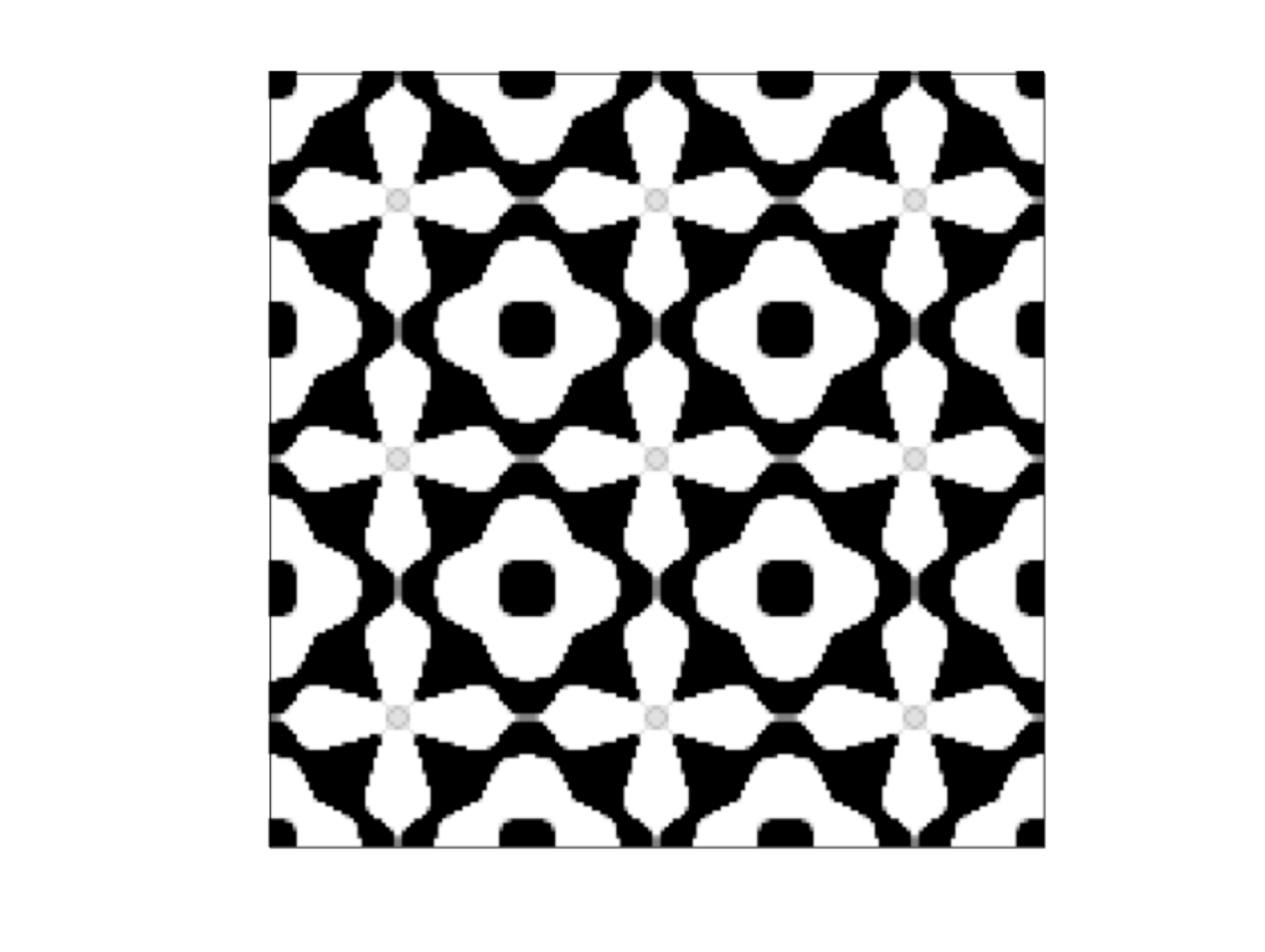}}
\subfigure[]{
\includegraphics[scale=0.35]{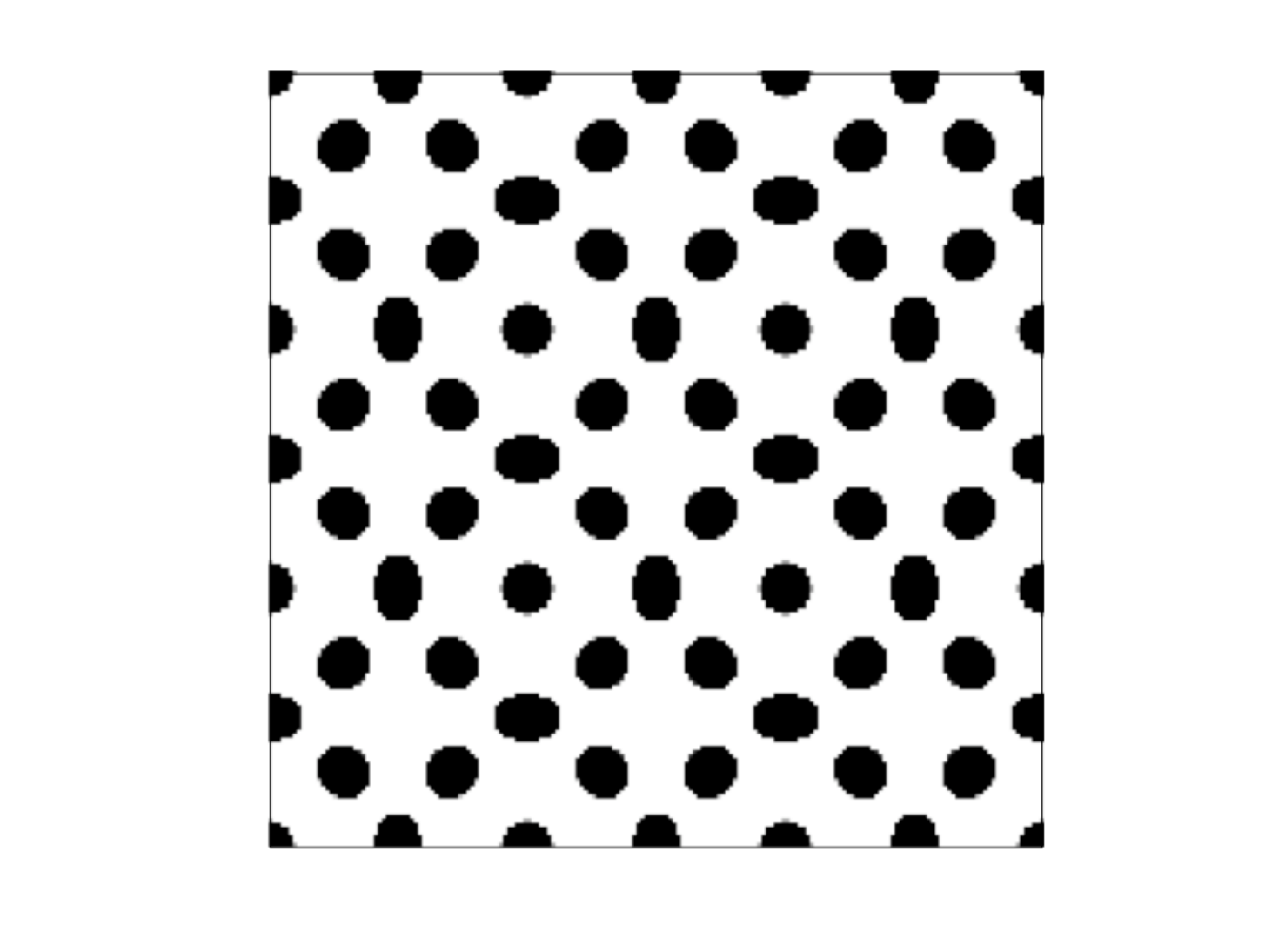}}
\subfigure[]{
\includegraphics[scale=0.35]{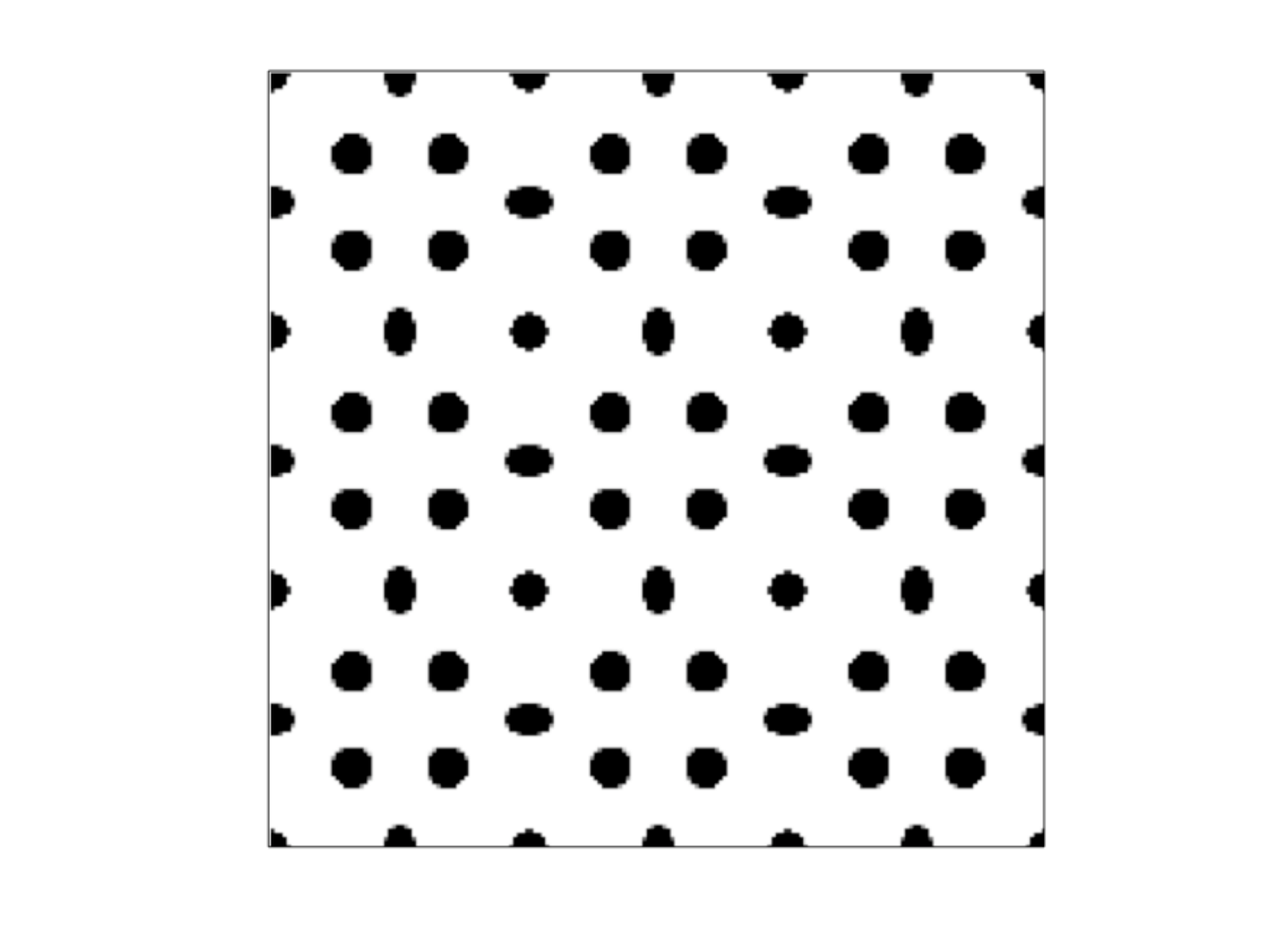}}
\caption{The evolution of the square lattice crystal structure for optimizing the gap-midgap ratio between $\lambda^7_{\text{TM}}$ and $\lambda^8_{\text{TM}}$.}
\label{figOR_3a}
\end{figure}

\begin{figure}
\centering
\subfigure[]{
\includegraphics[scale=0.3]{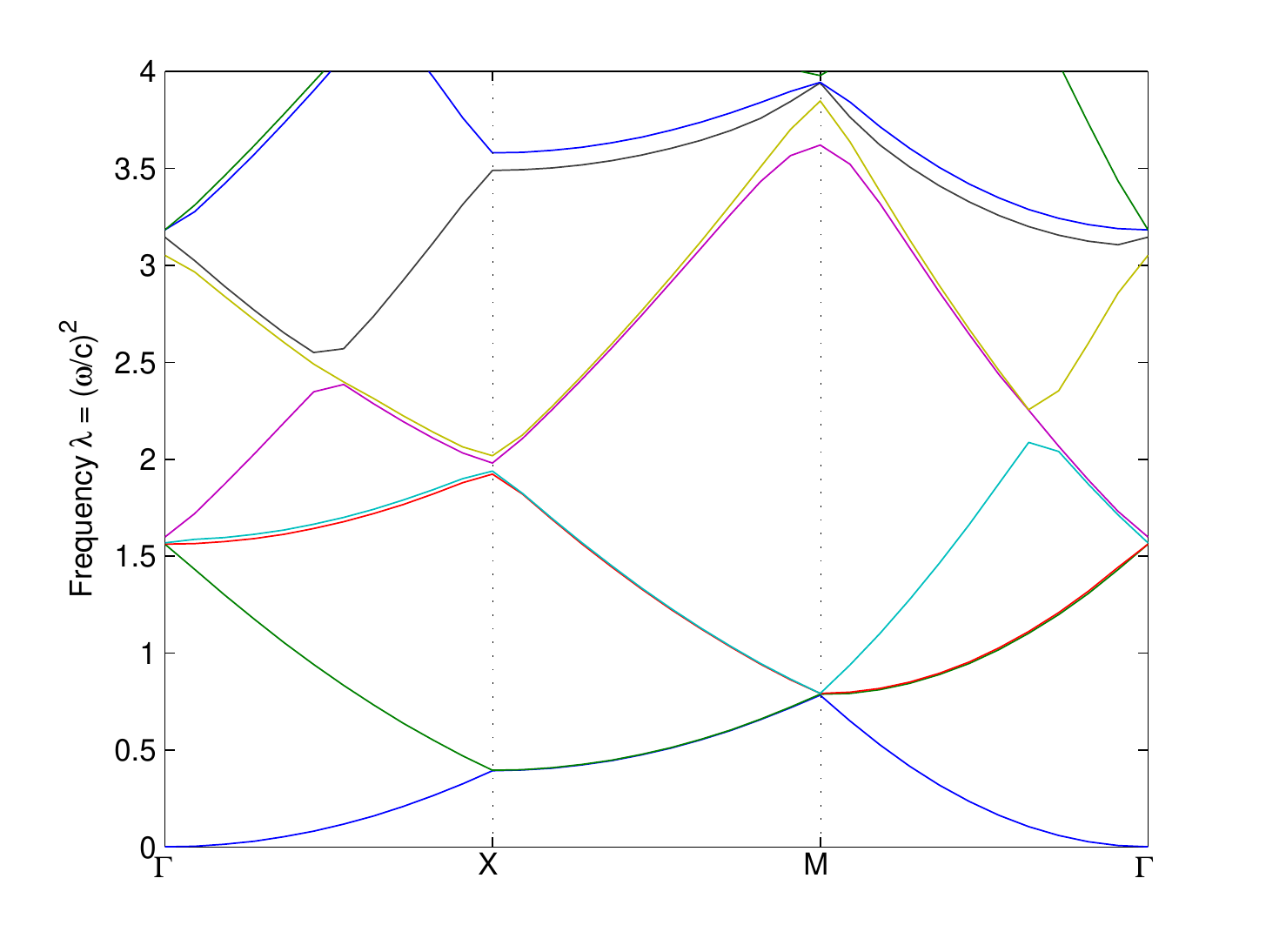}}
\subfigure[]{
\includegraphics[scale=0.3]{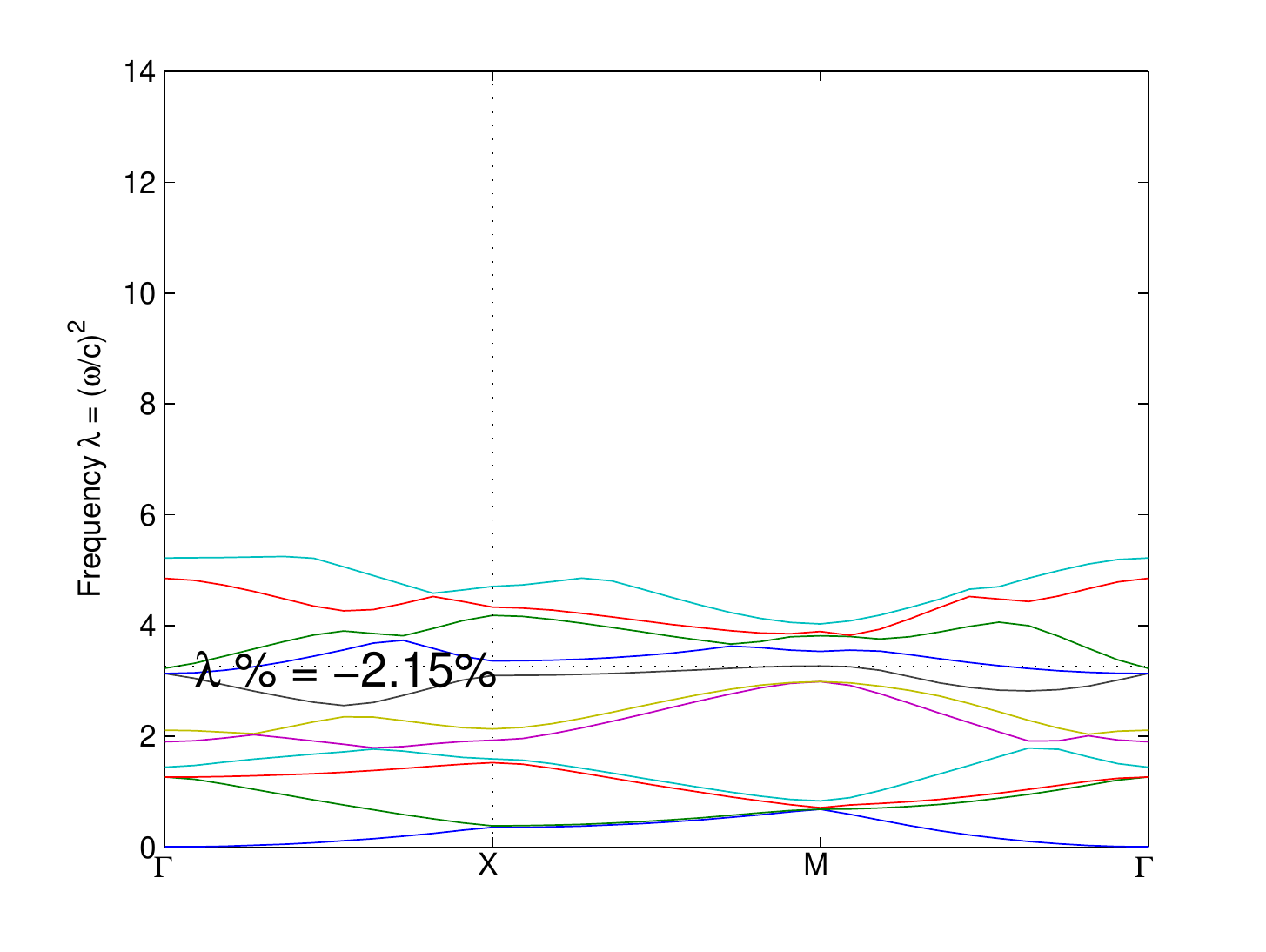}}
\subfigure[]{
\includegraphics[scale=0.3]{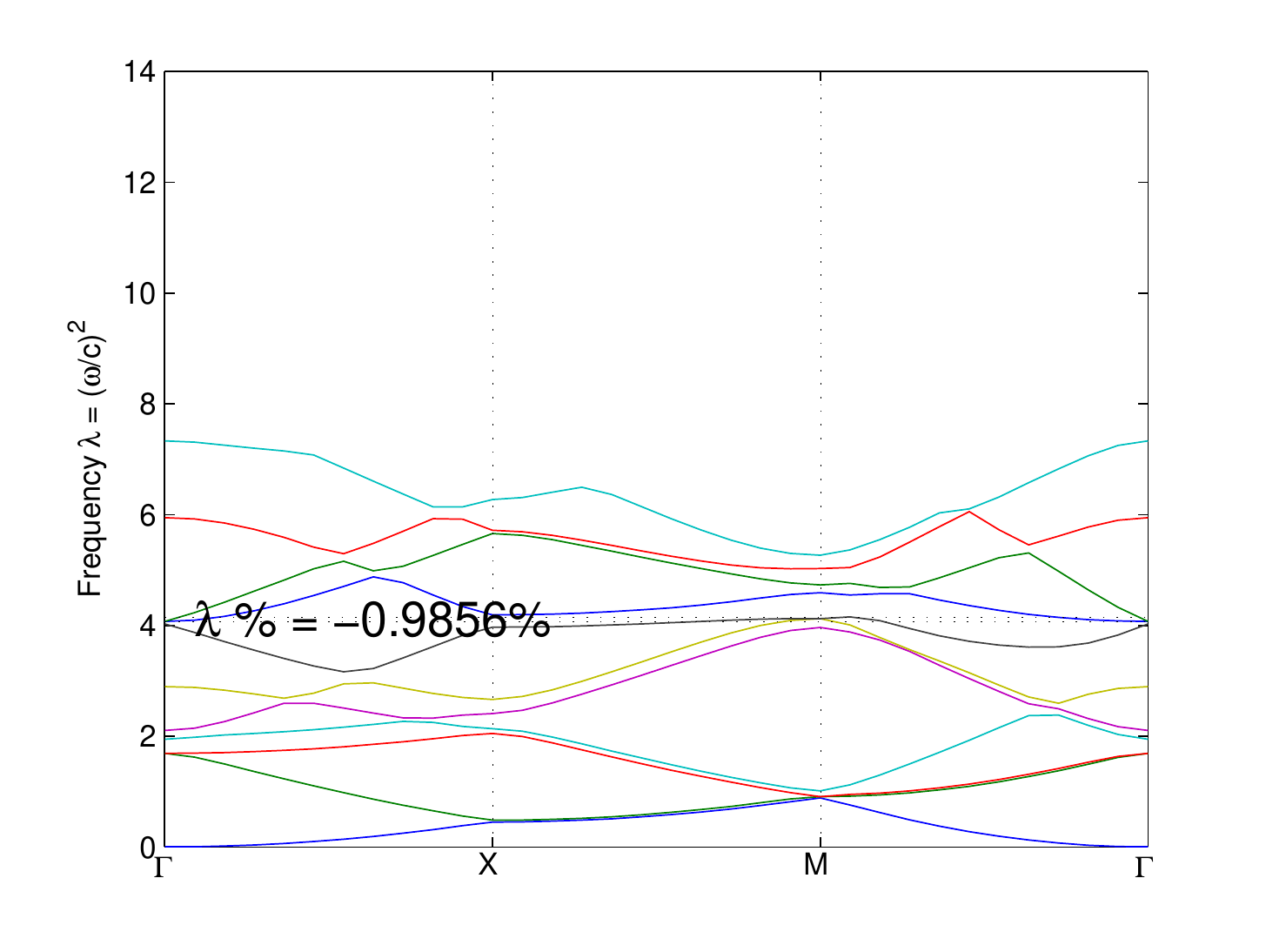}}
\subfigure[]{
\includegraphics[scale=0.3]{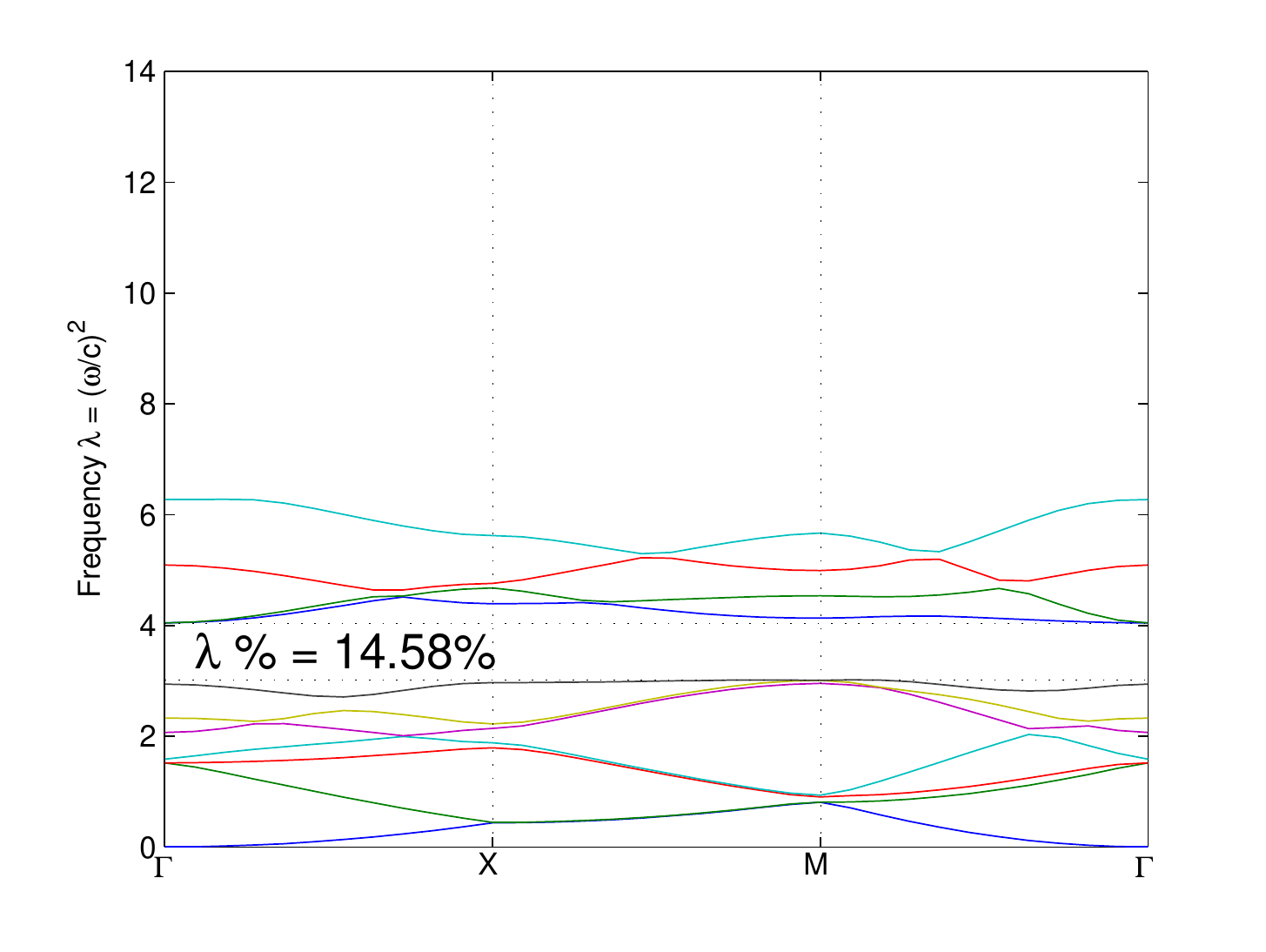}}
\subfigure[]{
\includegraphics[scale=0.3]{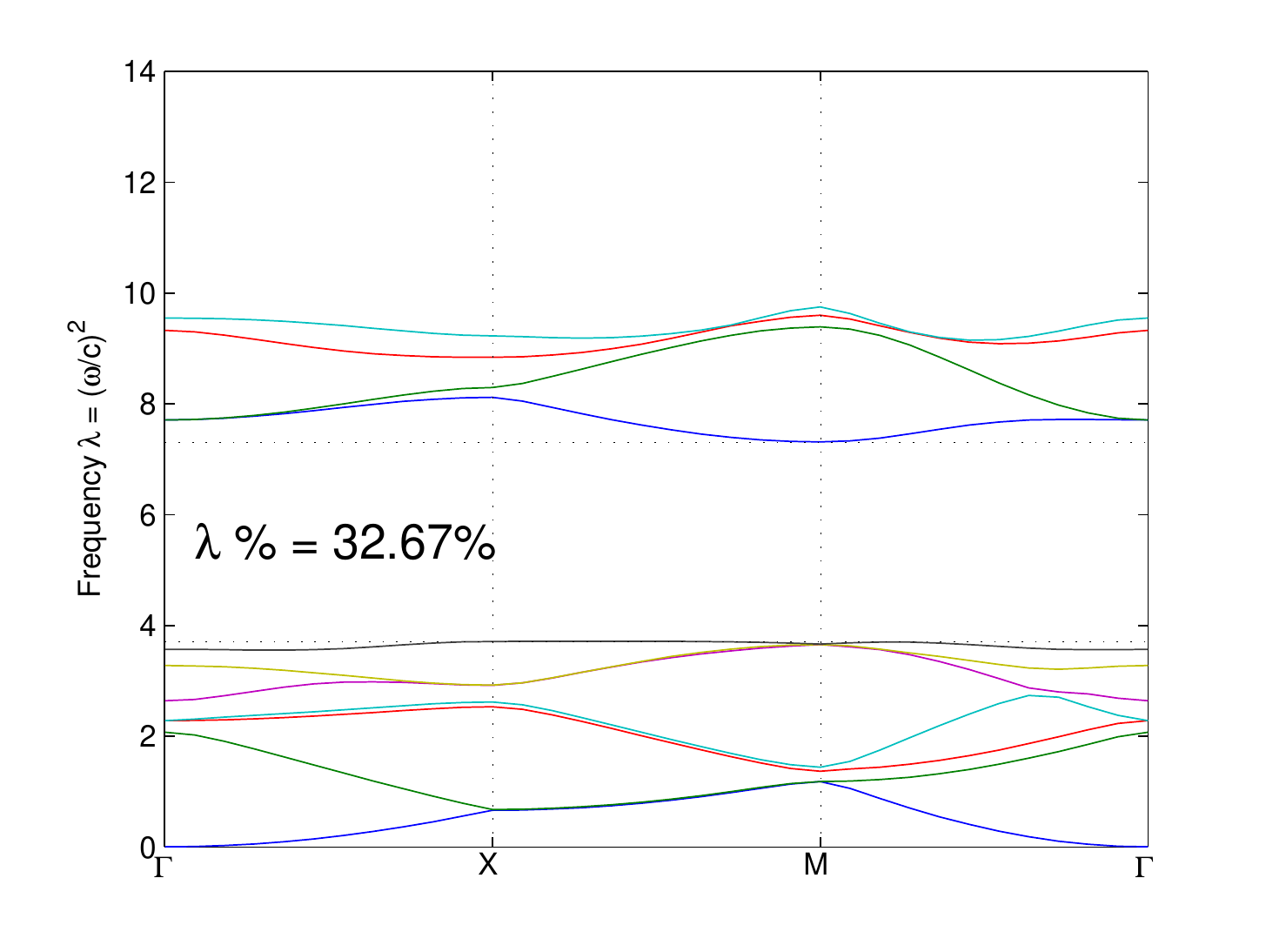}}
\subfigure[]{
\includegraphics[scale=0.3]{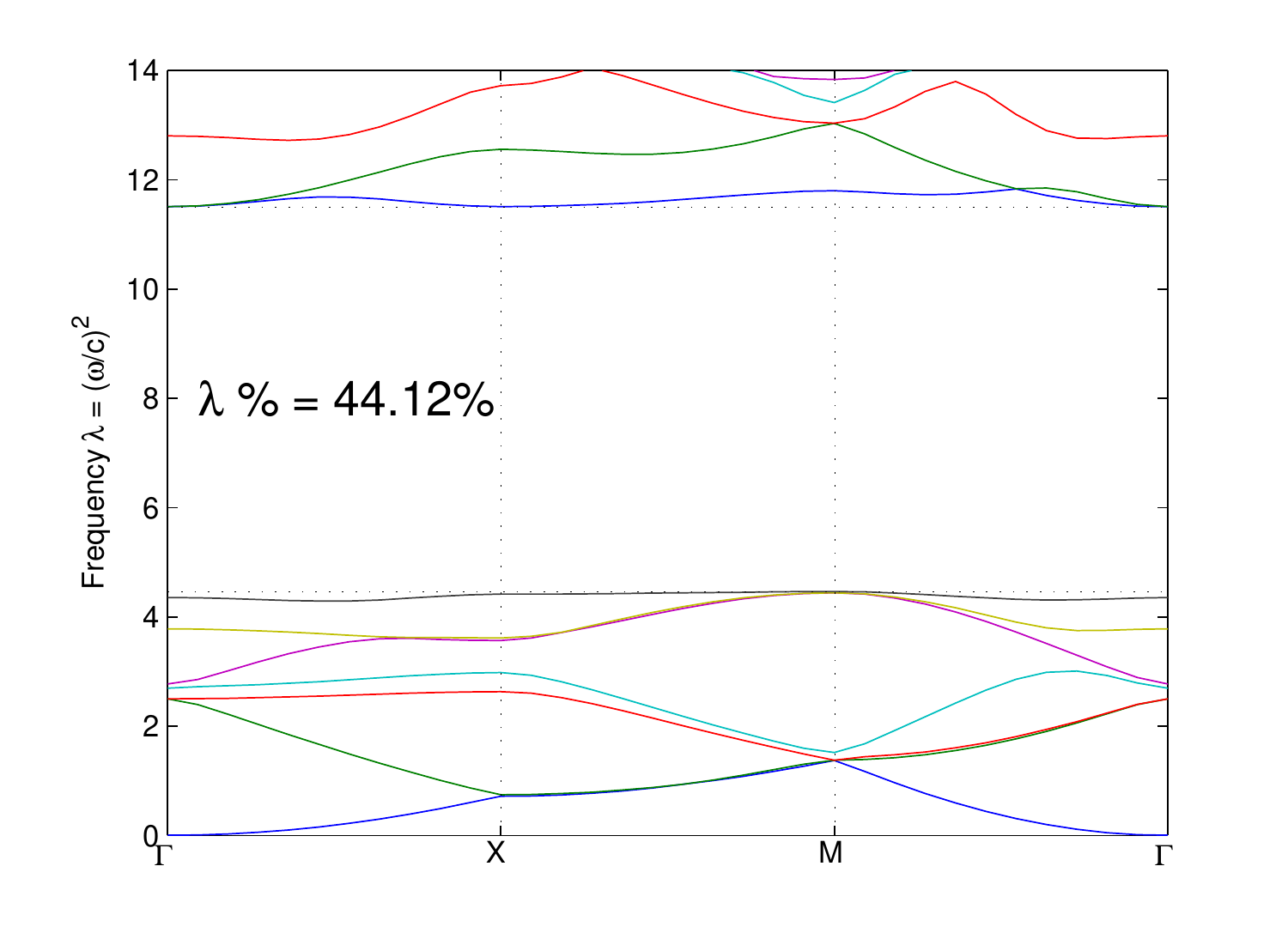}}
\caption{The corresponding band structure (of Figure~\ref{figOR_3a})) and the gap-midgap ratio between $\lambda^7_{\text{TM}}$ and $\lambda^8_{\text{TM}}$ in the square lattice.}
\label{figOR_3b}
\end{figure}

\begin{figure}
\centering
\subfigure[]{
\includegraphics[scale=0.35]{Figures/Sq_before_geometry}}
\subfigure[]{
\includegraphics[scale=0.35]{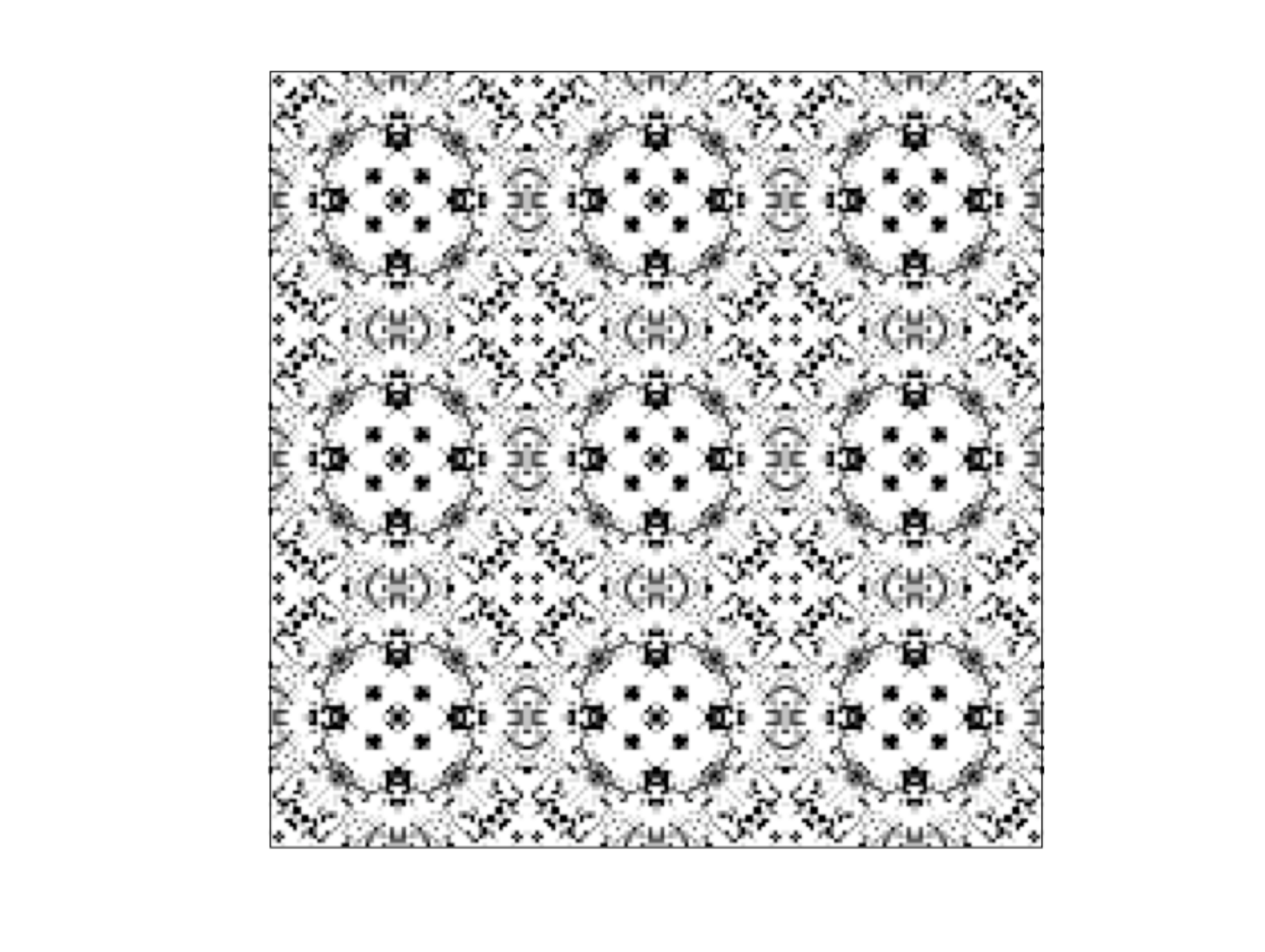}}
\subfigure[]{
\includegraphics[scale=0.35]{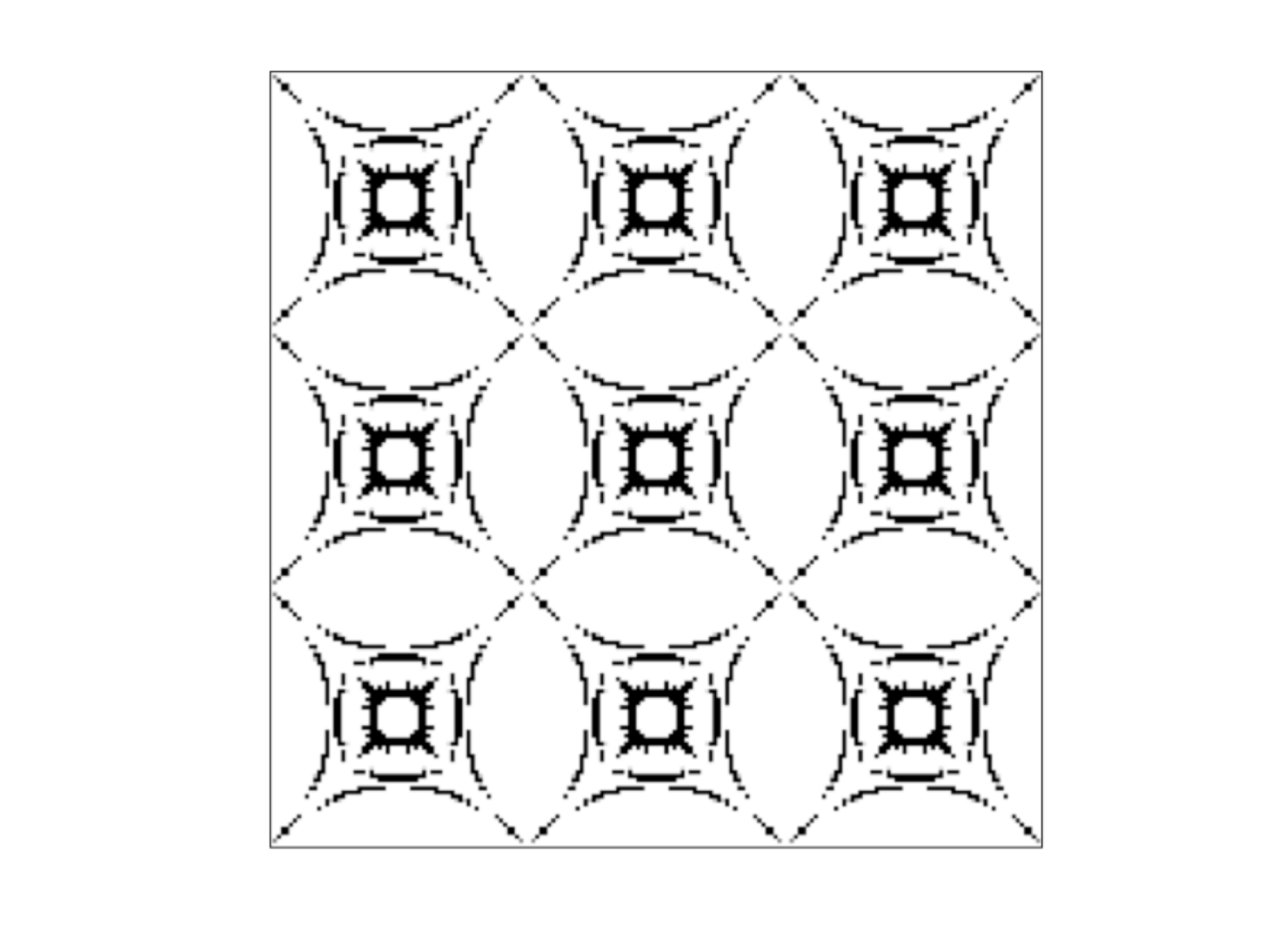}}
\subfigure[]{
\includegraphics[scale=0.35]{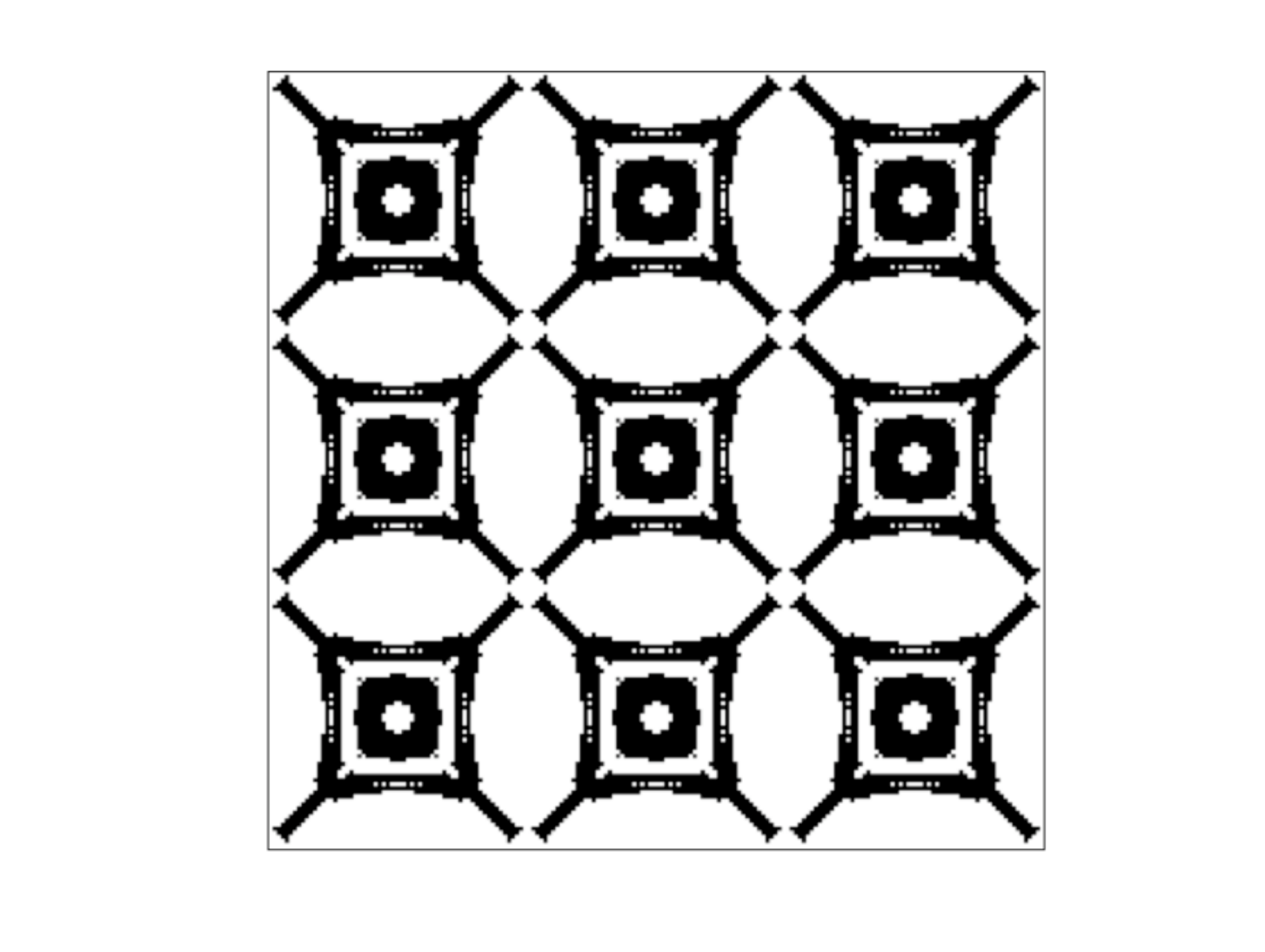}}
\subfigure[]{
\includegraphics[scale=0.35]{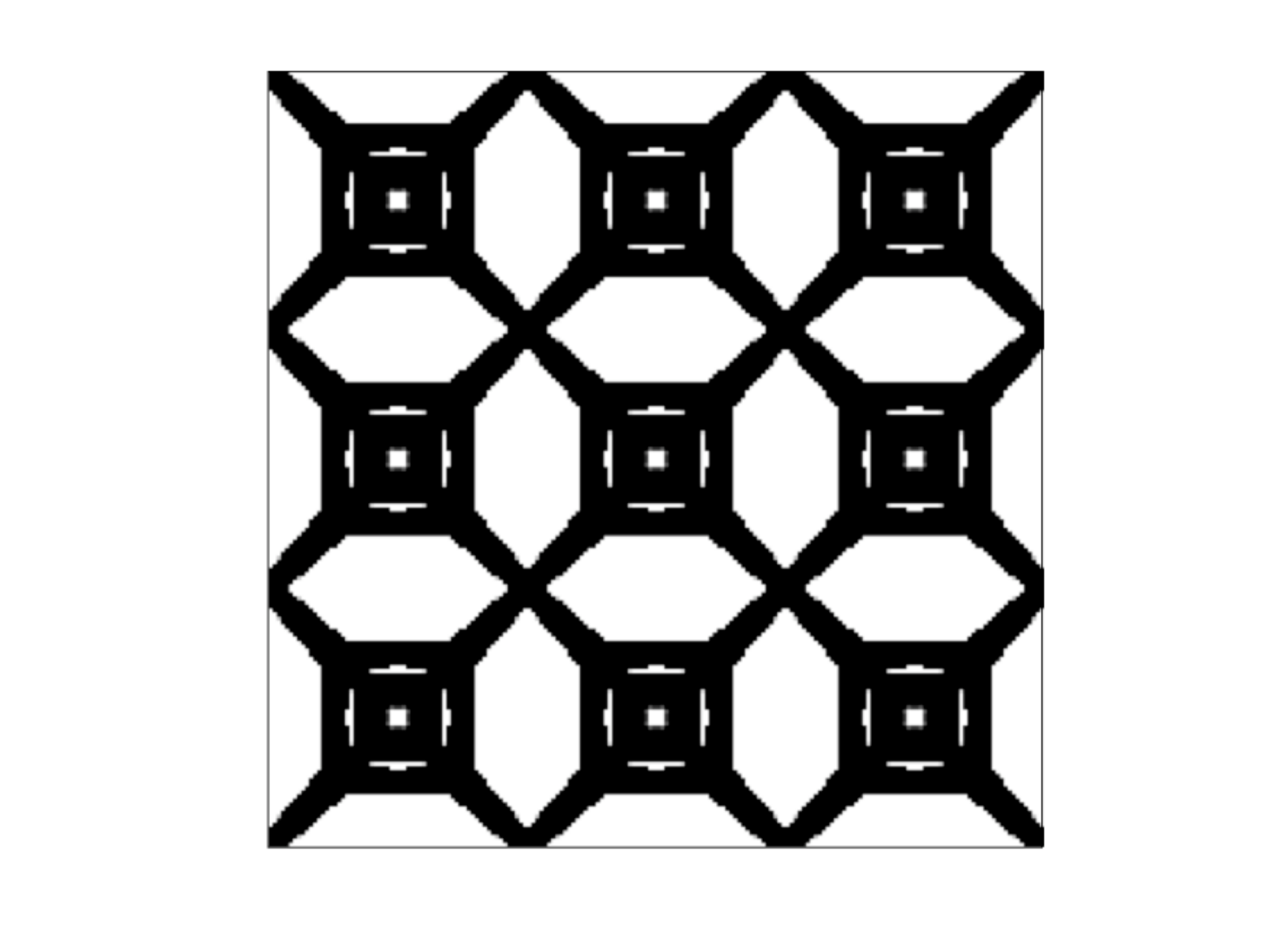}}
\subfigure[]{
\includegraphics[scale=0.35]{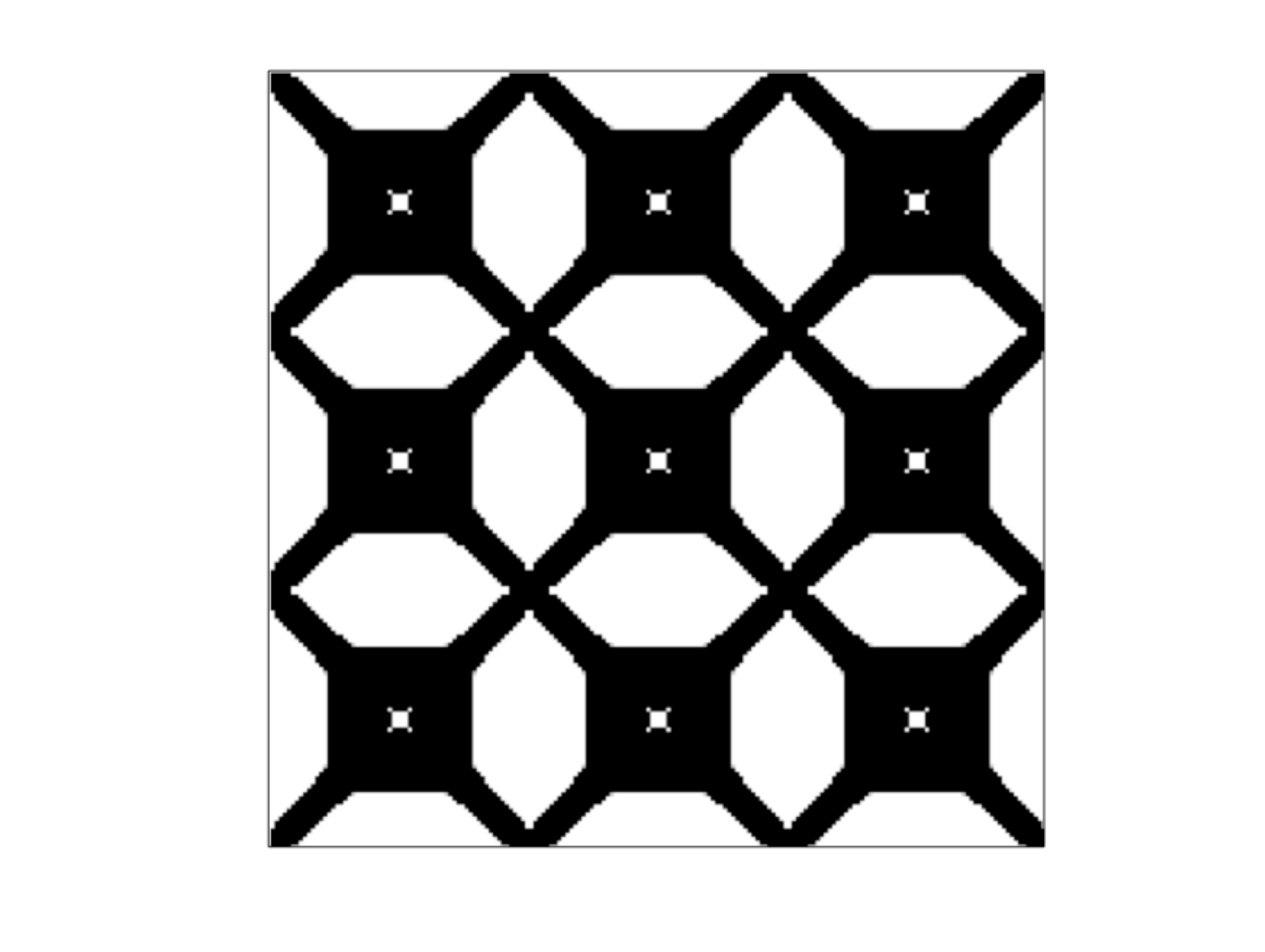}}
\caption{The evolution of the square lattice crystal structure for optimizing the gap-midgap ratio between $\lambda^3_{\text{TE}}$ and $\lambda^4_{\text{TE}}$.}
\label{figOR_3c}
\end{figure}

\begin{figure}
\centering
\subfigure[]{
\includegraphics[scale=0.3]{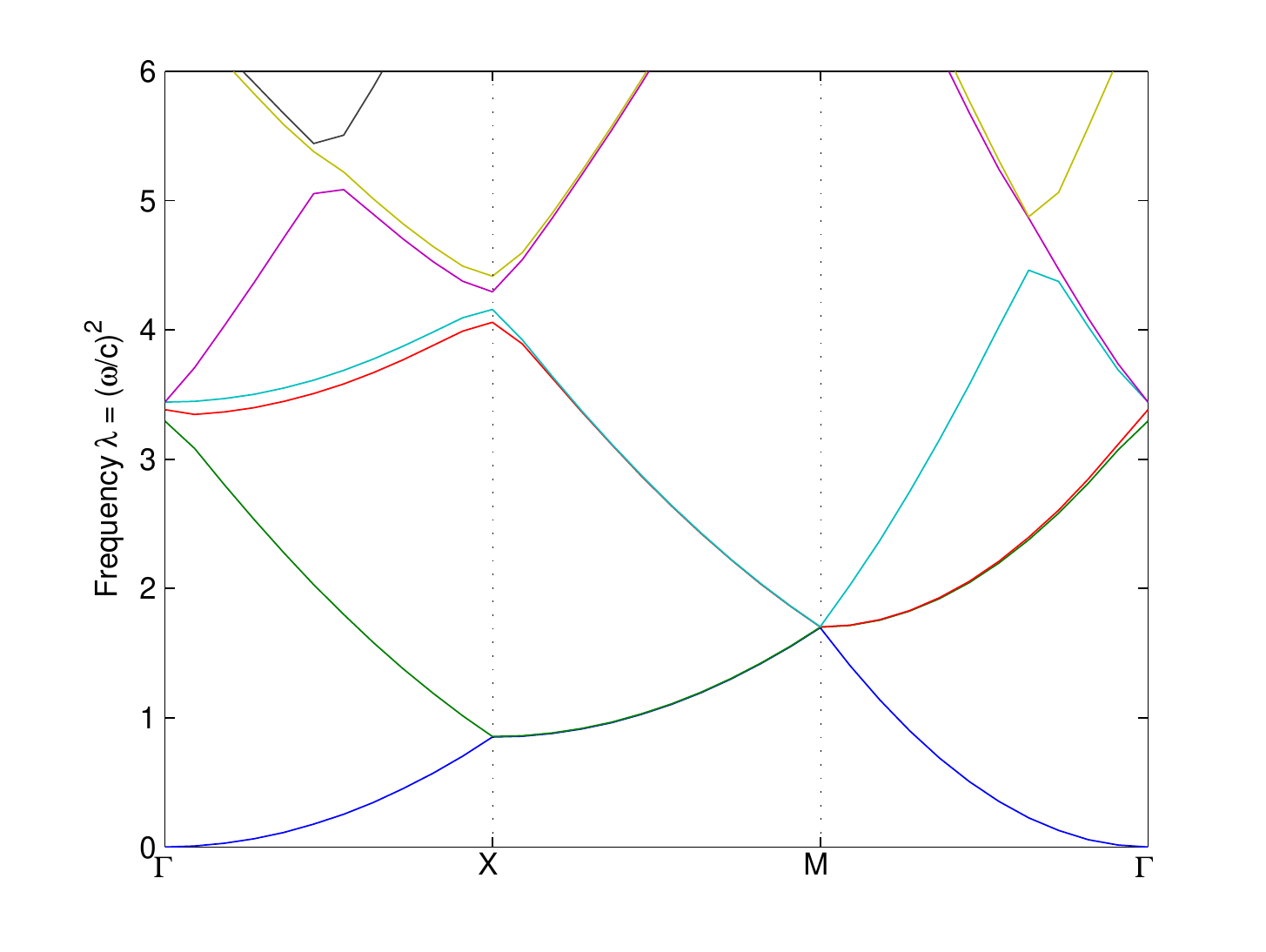}}
\subfigure[]{
\includegraphics[scale=0.3]{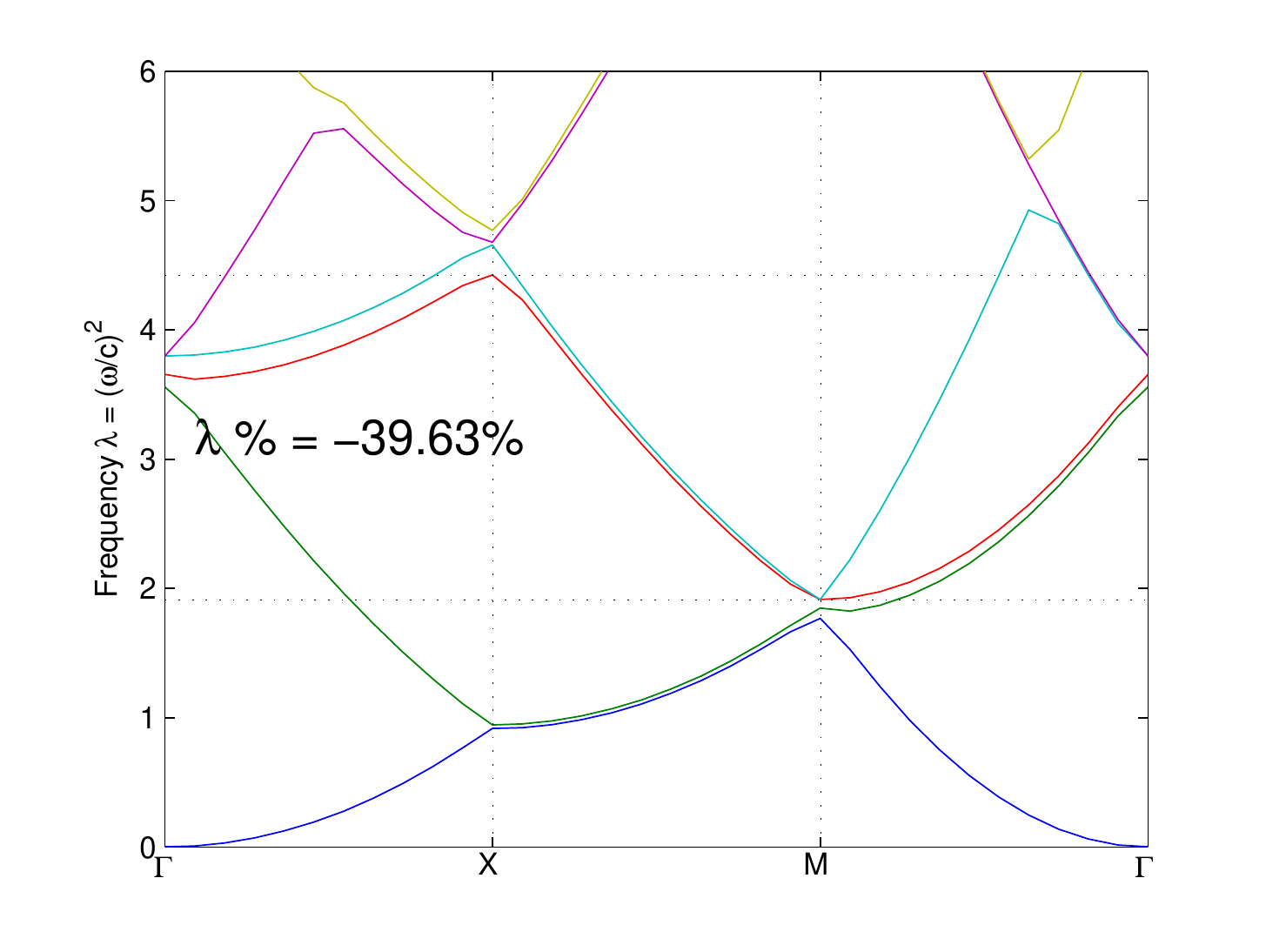}}
\subfigure[]{
\includegraphics[scale=0.3]{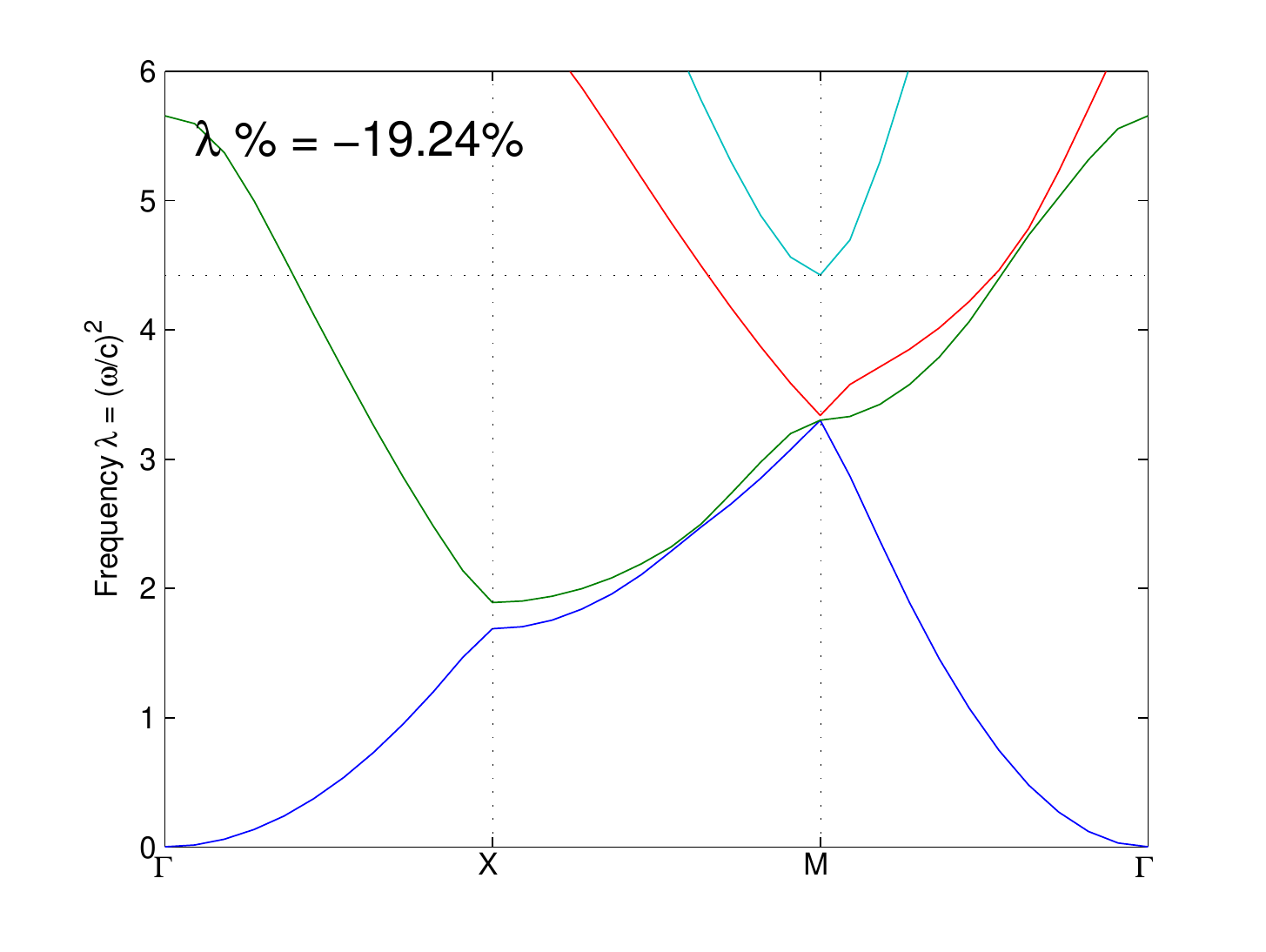}}
\subfigure[]{
\includegraphics[scale=0.3]{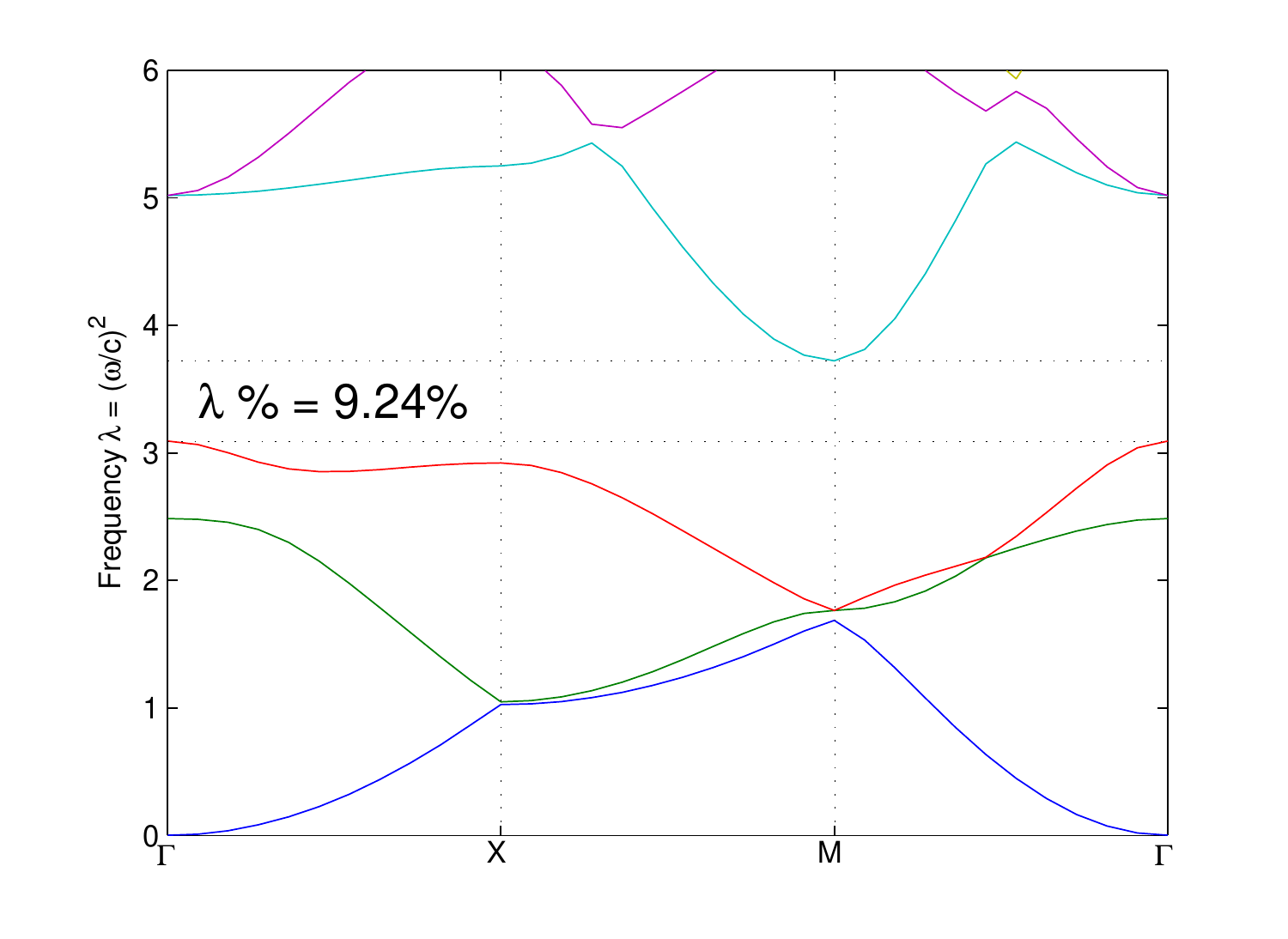}}
\subfigure[]{
\includegraphics[scale=0.3]{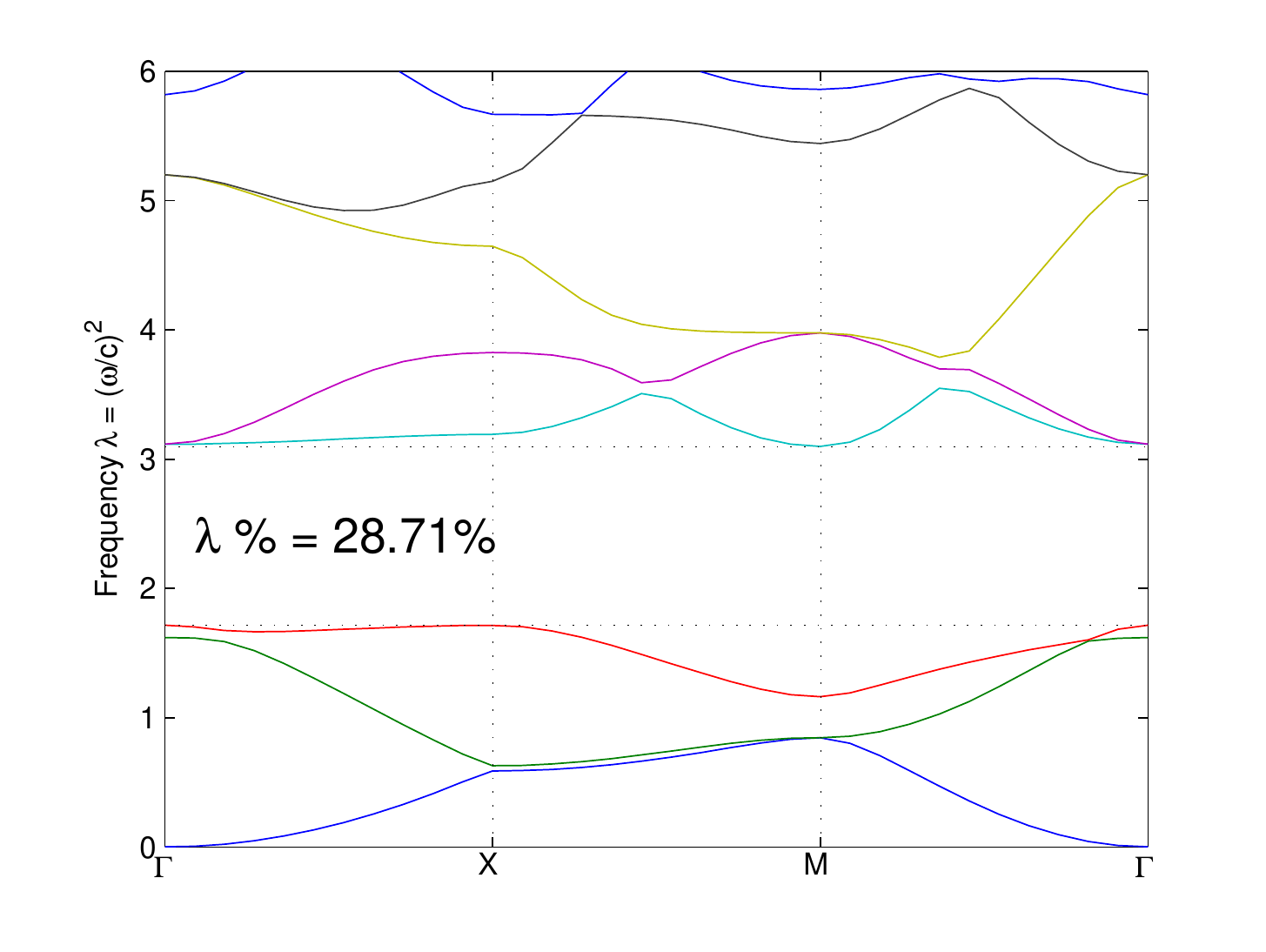}}
\subfigure[]{
\includegraphics[scale=0.3]{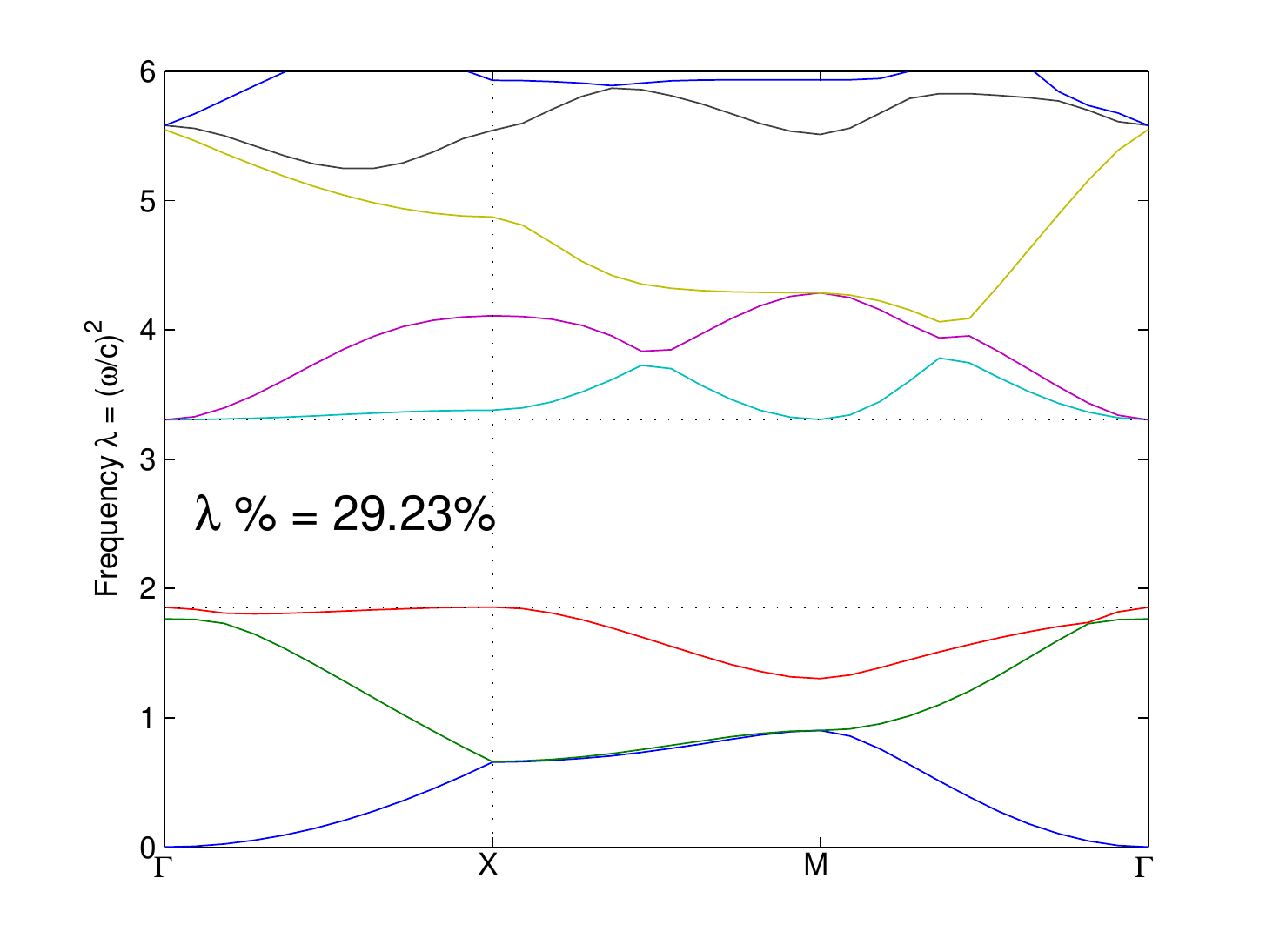}}
\caption{The corresponding band structure (of Figure~\ref{figOR_3c})) and the gap-midgap ratio between $\lambda^3_{\text{TE}}$ and $\lambda^4_{\text{TE}}$ in the square lattice.}
\label{figOR_3d}
\end{figure}

We start the optimization procedure with a random distribution of the
dielectric, such as the one shown in
Figure~\ref{figOR_beforeOpt}(a). The corresponding band structures of
the TE and TM fields are shown in Figure \ref{figOR_beforeOpt}(b).  In
Figure~\ref{figOR_3a}, we present an example of the evolution of the
crystal structure as the optimization process progresses. (The light
color indicates the low dielectric constant and the dark color denotes
the high dielectric constant.) As illustrated in
Figure~\ref{figOR_3b}, the gap-midgap ratio starts from a negative
value ($-8.93\%$) corresponding to the random configuration (Figure
\ref{figOR_3b}(a)) and increases up to $+43.90\%$ corresponding to the
optimal configuration (Figure \ref{figOR_3b}(f)) at which time the
optimization process terminates successfully. Another example of the
optimization evolution for TE polarization is shown in
Figure~\ref{figOR_3c} and Figure~\ref{figOR_3d}, in which the
gap-midgap ratio increases from $-39.21\%$ to $+29.23\%$.

\begin{figure}
\centering
\subfigure[Optimal crystal structure]{
\includegraphics[scale=0.45]{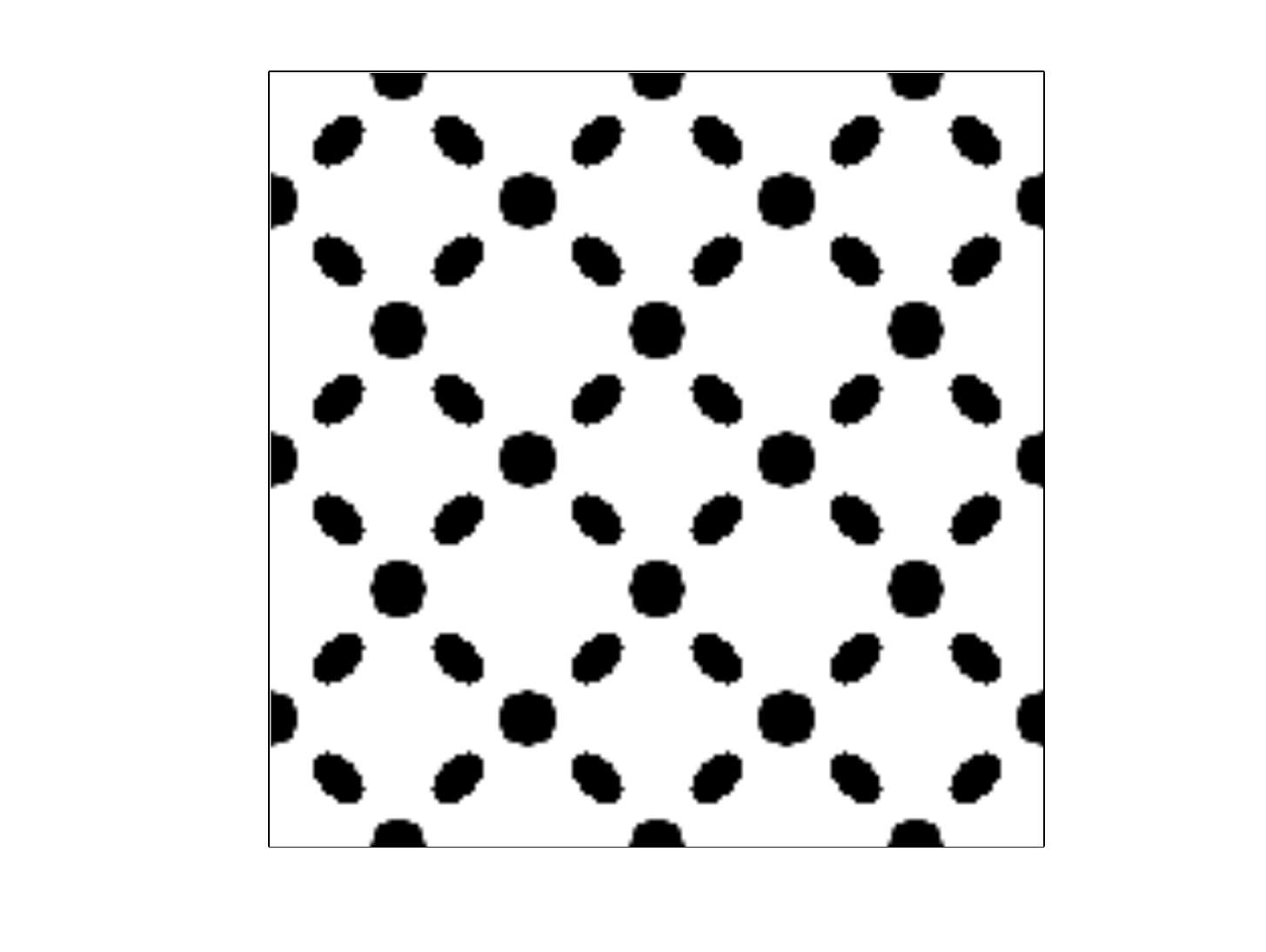}}
\subfigure[Optimal band structure]{
\includegraphics[scale=0.45]{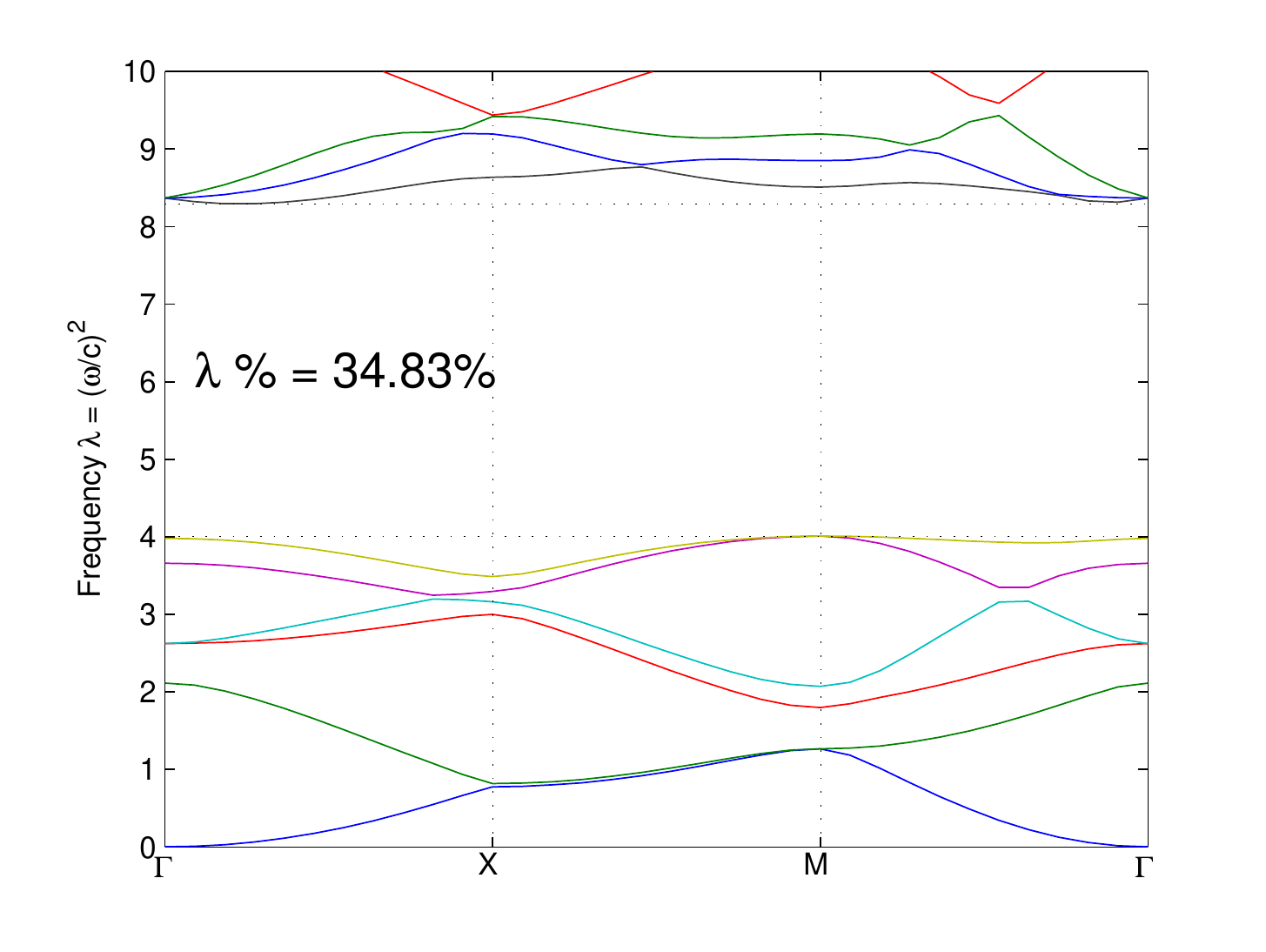}}
\caption{Optimization of band gap between $\lambda^6_{\text{TM}}$ and $\lambda^7_{\text{TM}}$ in the square lattice.}
\label{figOR_TM67}
\end{figure}

\begin{figure}
\centering
\subfigure[Optimal crystal structure]{
\includegraphics[scale=0.45]{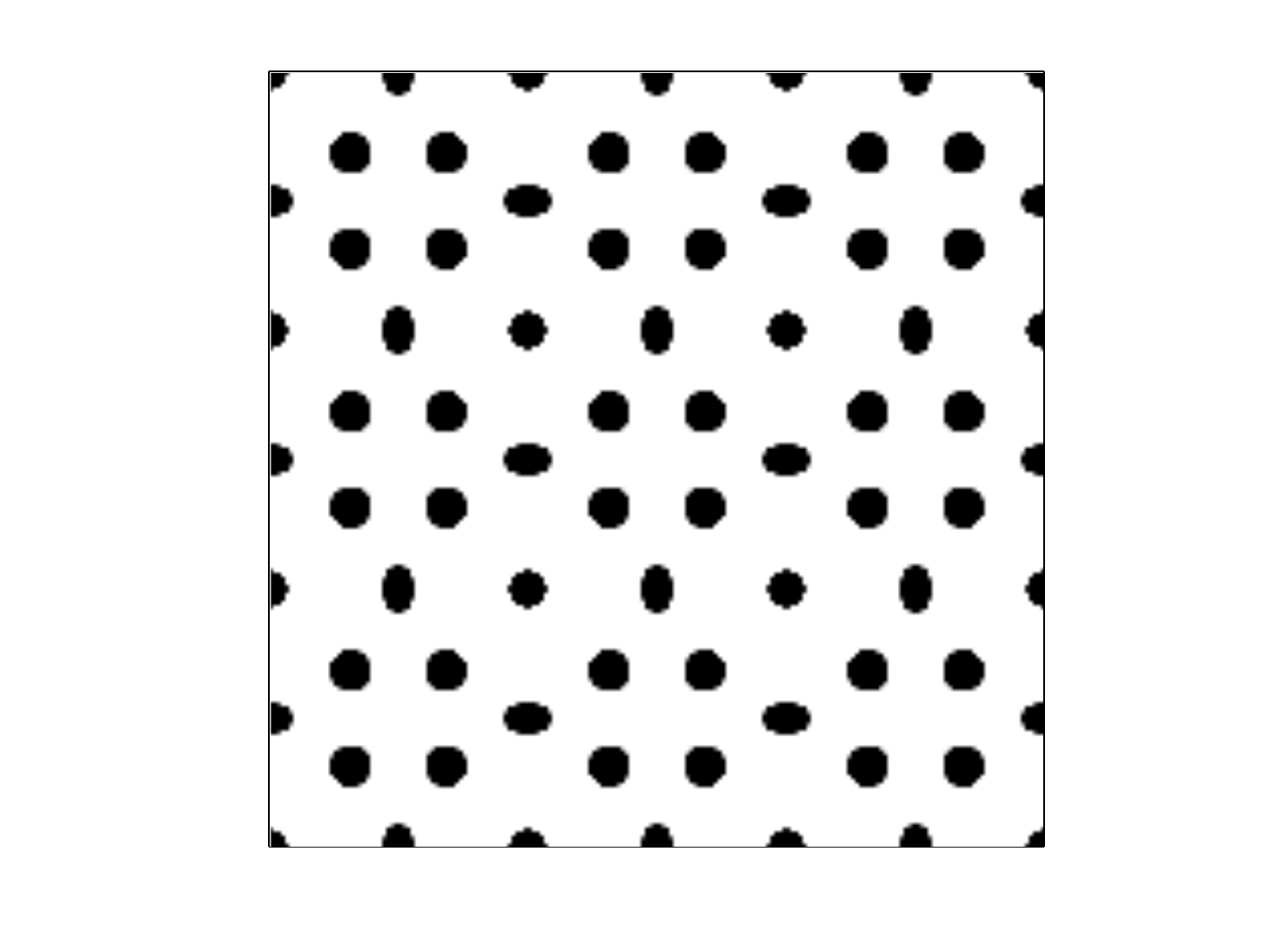}}
\subfigure[Optimal band structure]{
\includegraphics[scale=0.45]{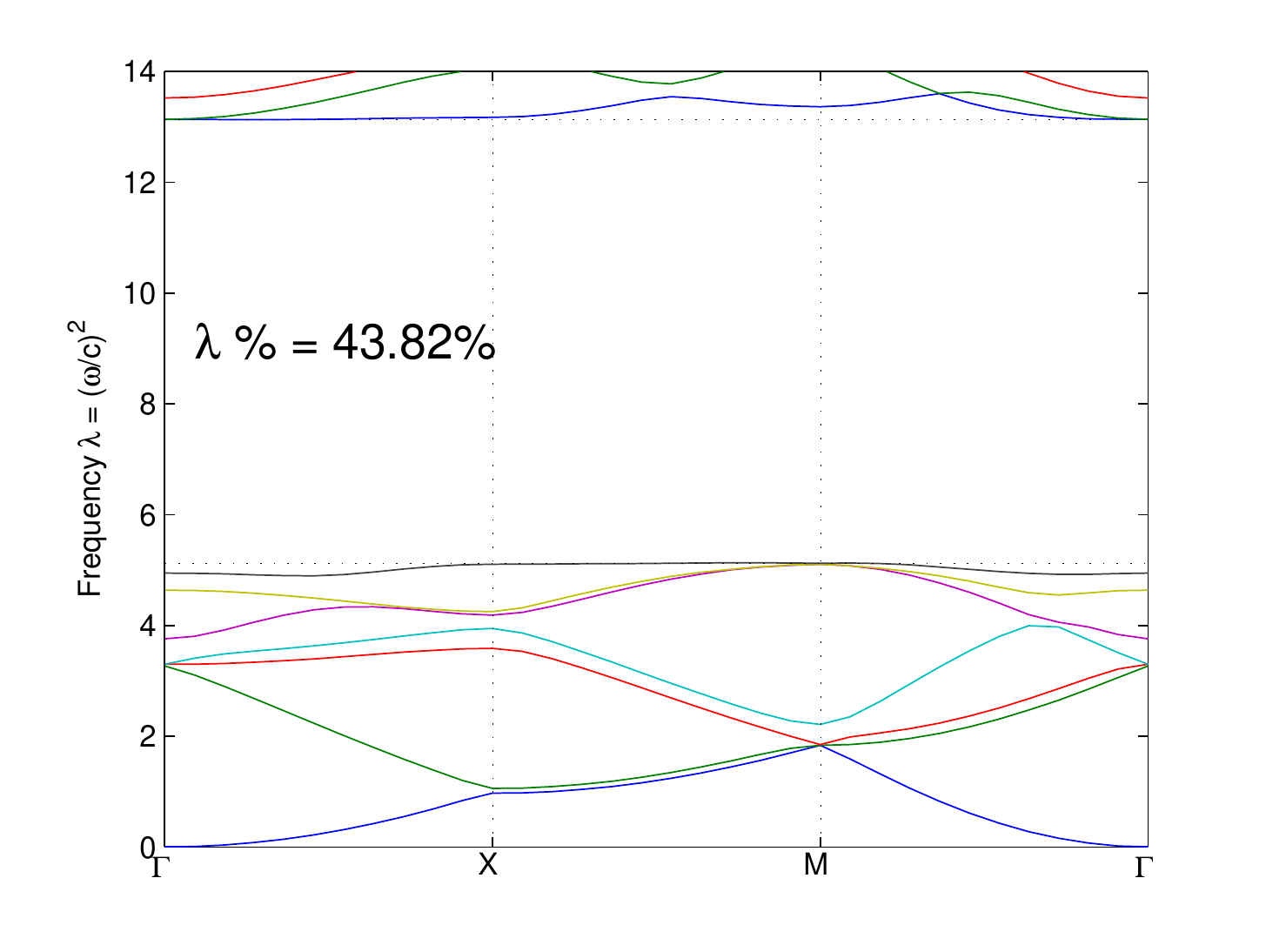}}
\caption{Optimization of band gap between $\lambda^7_{\text{TM}}$ and $\lambda^8_{\text{TM}}$ in the square lattice.}
\label{figOR_TM78}
\end{figure}

\begin{figure}
\centering
\subfigure[Optimal crystal structure]{
\includegraphics[scale=0.45]{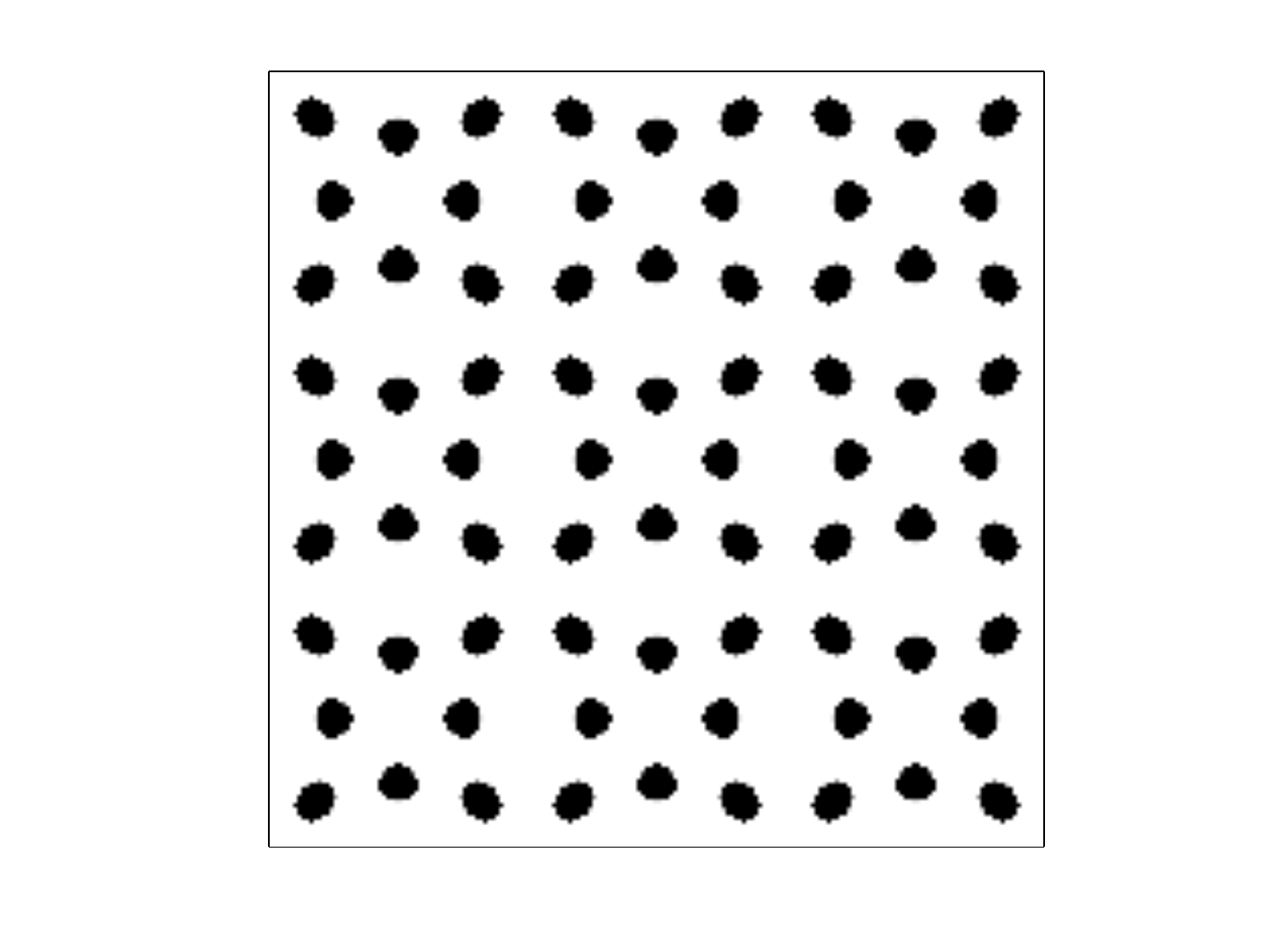}}
\subfigure[Optimal band structure]{
\includegraphics[scale=0.45]{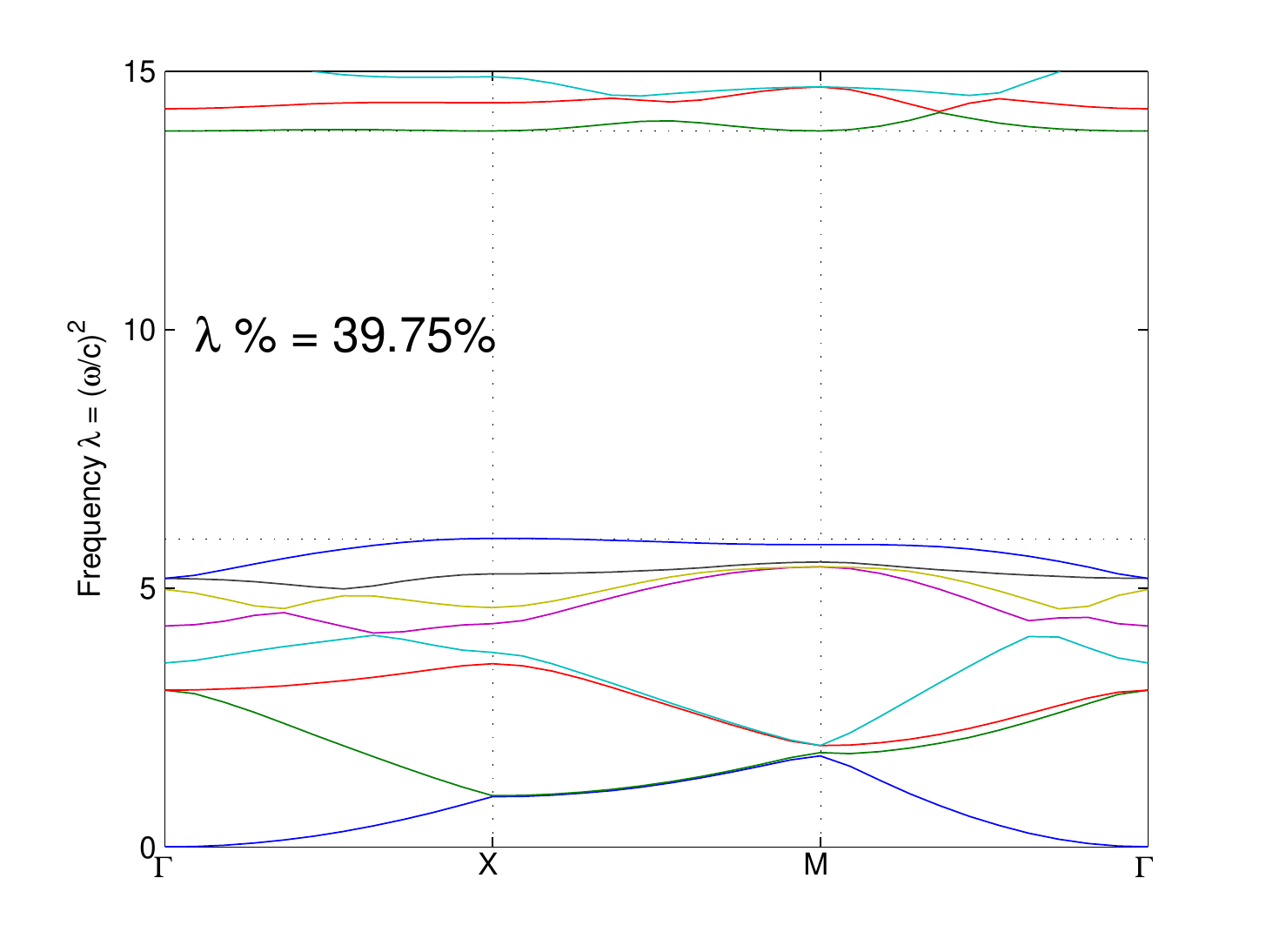}}
\caption{Optimization of band gap between $\lambda^8_{\text{TM}}$ and $\lambda^9_{\text{TM}}$ in the square lattice.}
\label{figOR_TM89}
\end{figure}

\begin{figure}
\centering
\subfigure[Optimal crystal structure]{
\includegraphics[scale=0.45]{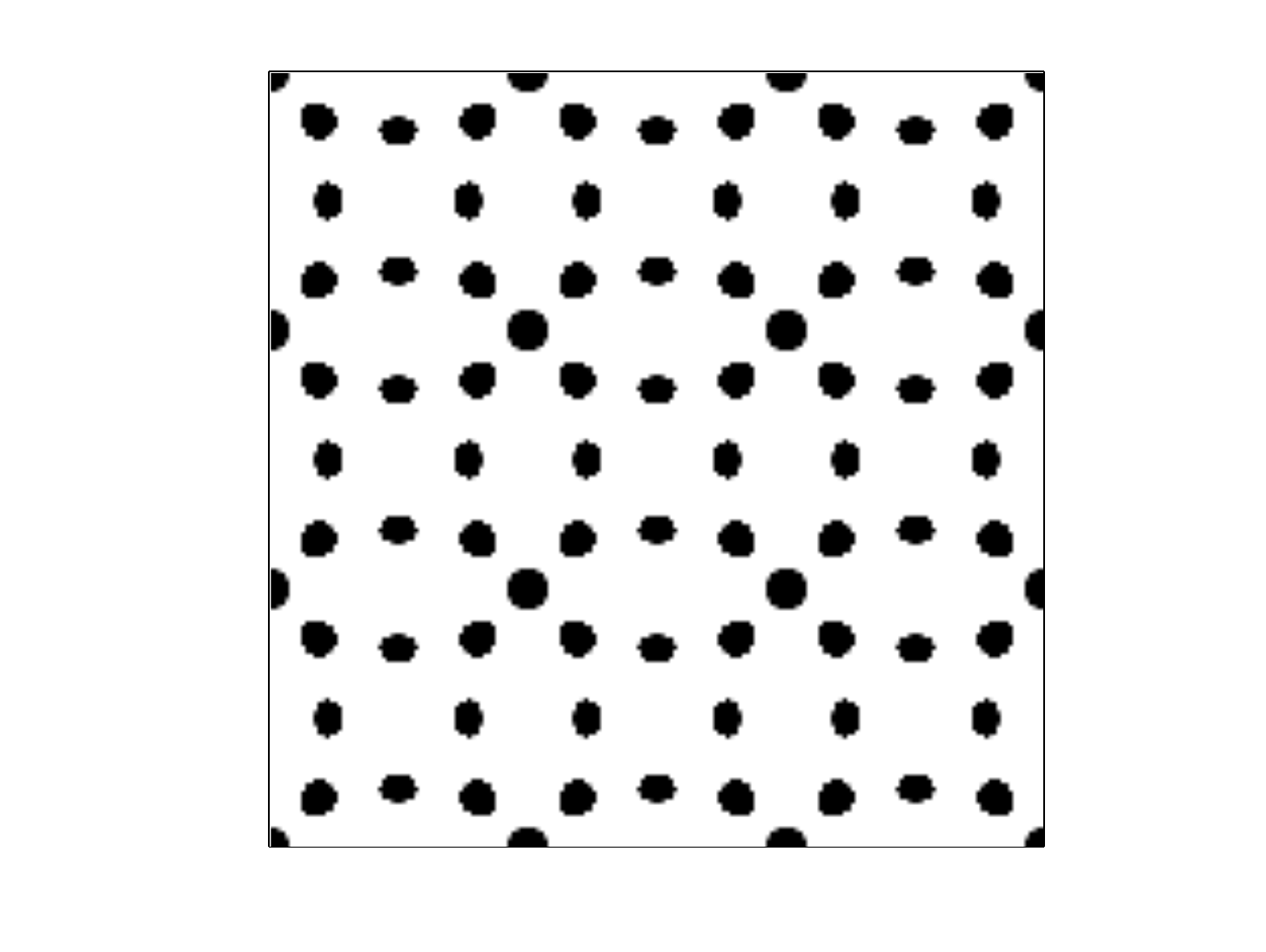}}
\subfigure[Optimal band structure]{
\includegraphics[scale=0.45]{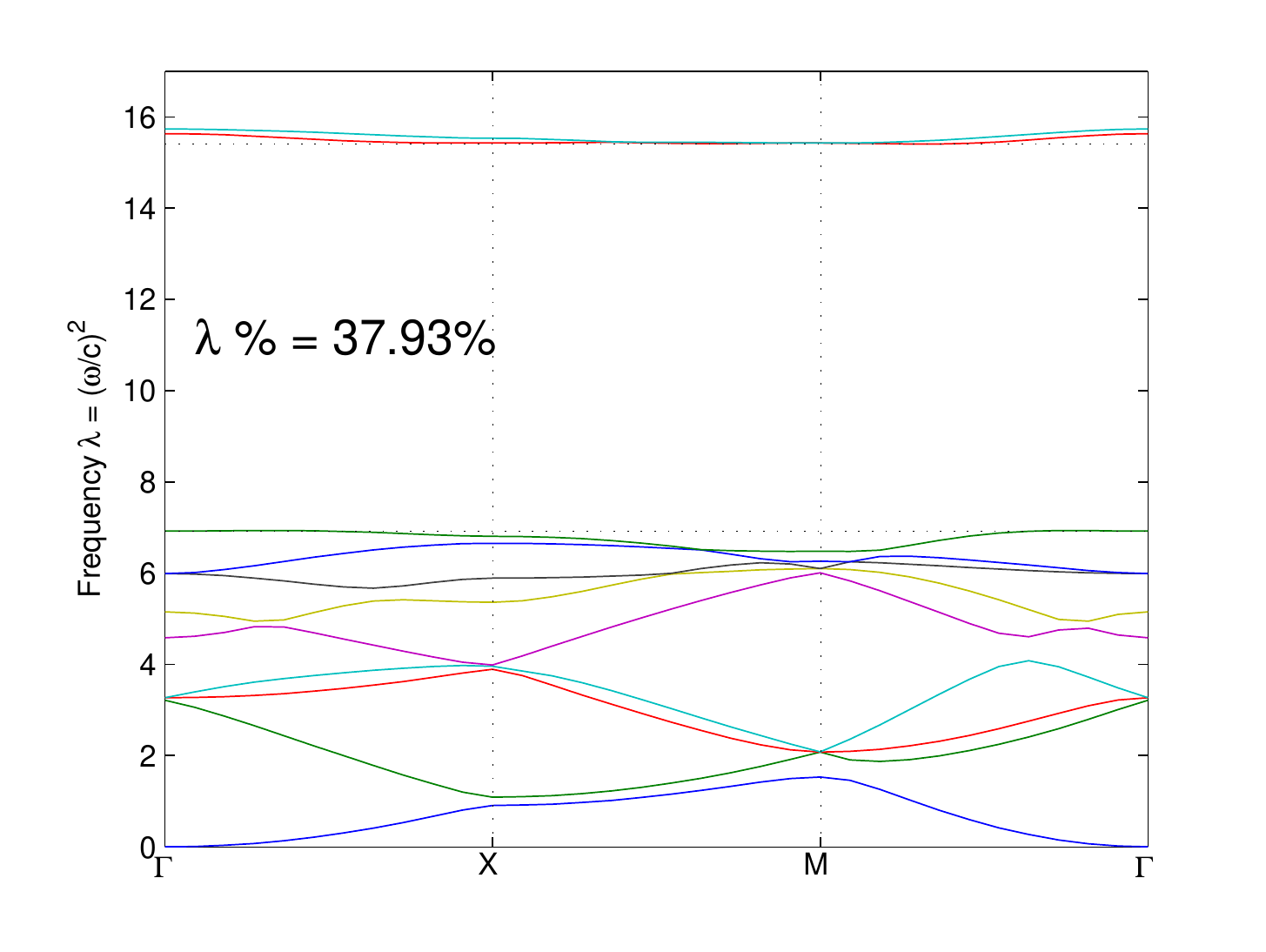}}
\caption{Optimization of band gap between $\lambda^9_{\text{TM}}$ and $\lambda^{10}_{\text{TM}}$ in the square lattice.}
\label{figOR_TM910}
\end{figure}

\begin{figure}
\centering
\subfigure[Optimal crystal structure]{
\includegraphics[scale=0.45]{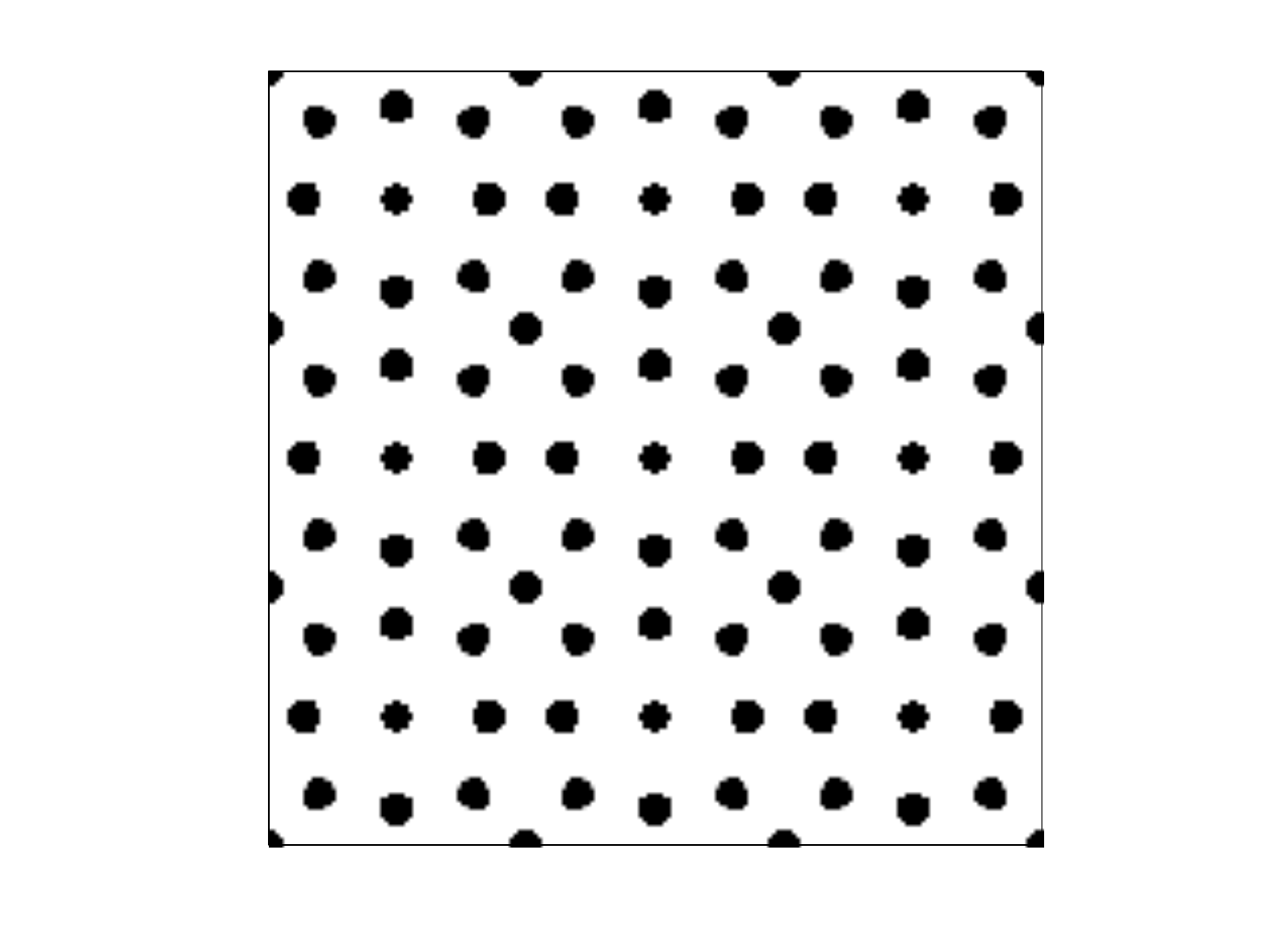}}
\subfigure[Optimal band structure]{
\includegraphics[scale=0.45]{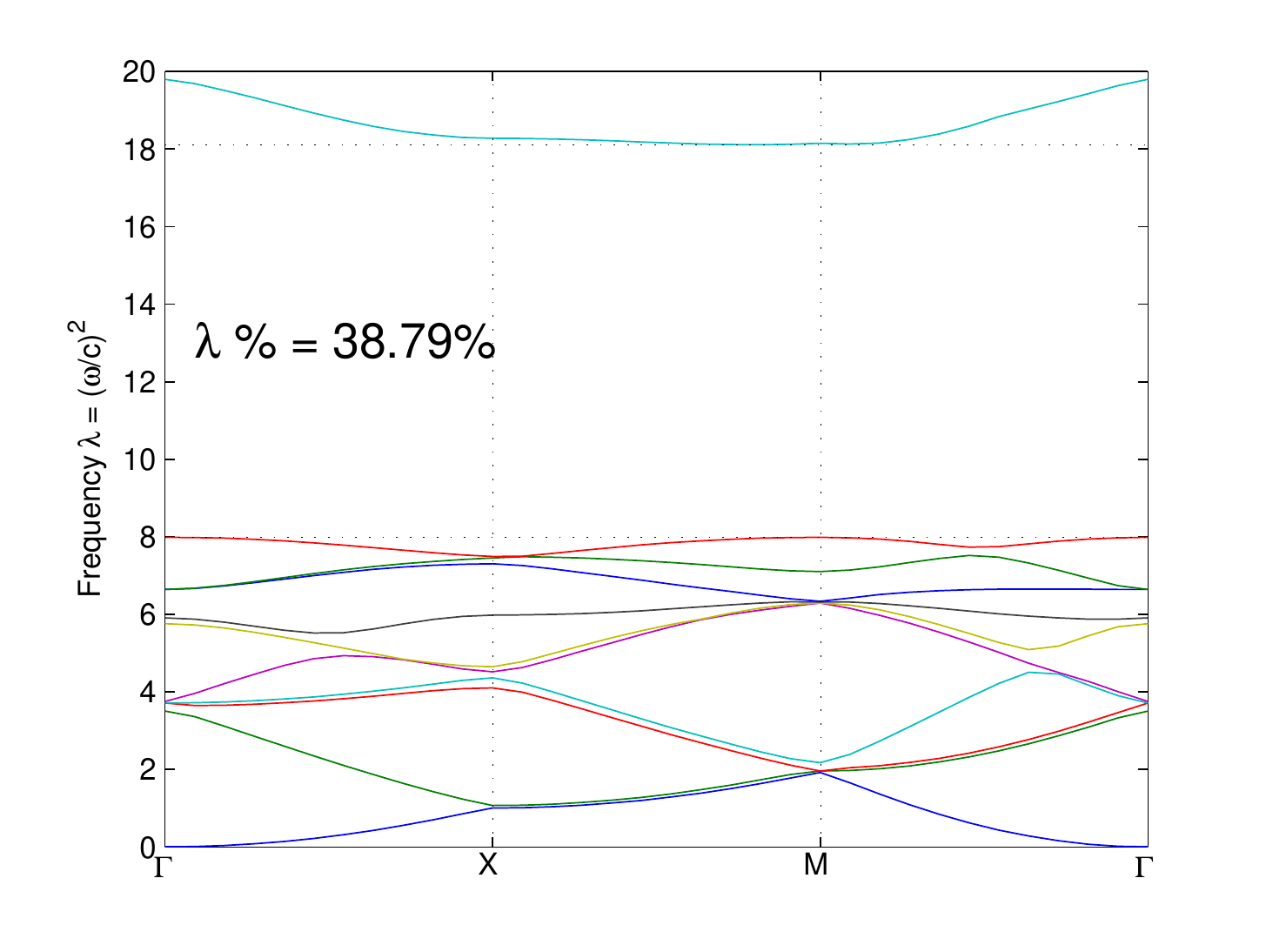}}
\caption{Optimization of band gap between $\lambda^{10}_{\text{TM}}$ and $\lambda^{11}_{\text{TM}}$ in the square lattice.}
\label{figOR_TM1011}
\end{figure}

\begin{figure}
\centering
\subfigure[Optimal crystal structure]{
\includegraphics[scale=0.45]{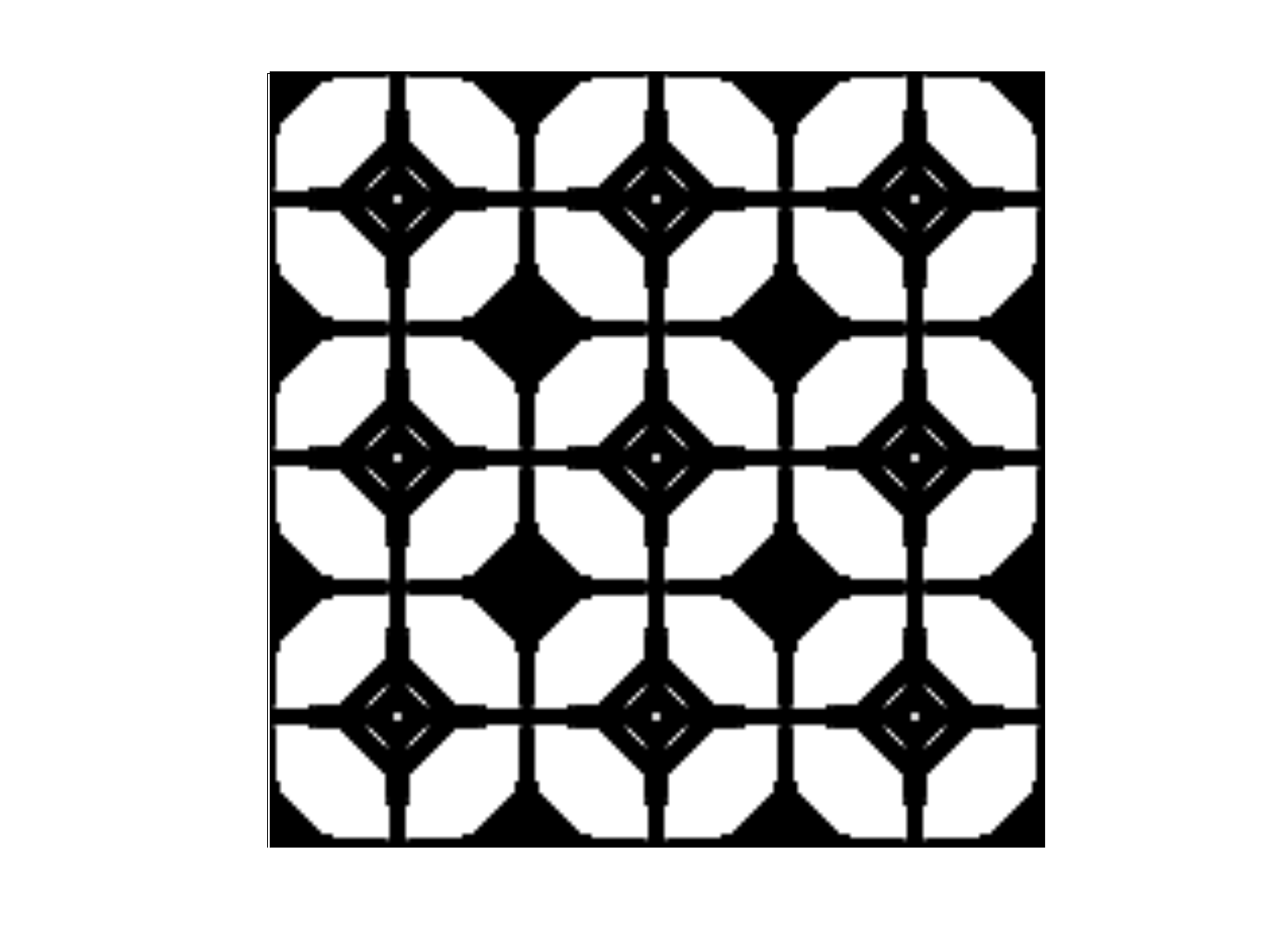}}
\subfigure[Optimal band structure]{
\includegraphics[scale=0.45]{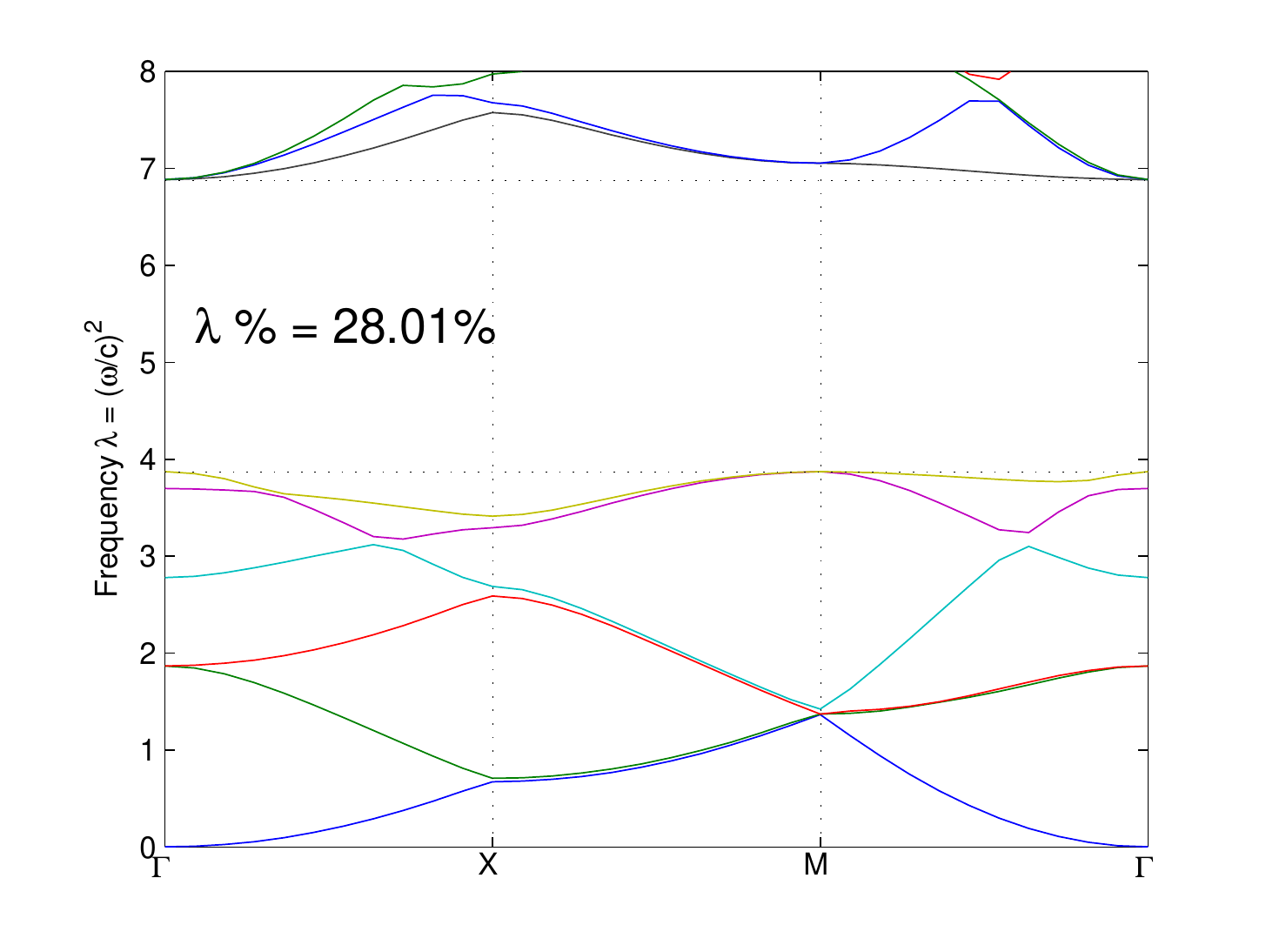}}
\caption{Optimization of band gap between $\lambda^6_{\text{TE}}$ and $\lambda^7_{\text{TE}}$ in the square lattice.}
\label{figOR_TE67}
\end{figure}

\begin{figure}
\centering
\subfigure[Optimal crystal structure]{
\includegraphics[scale=0.45]{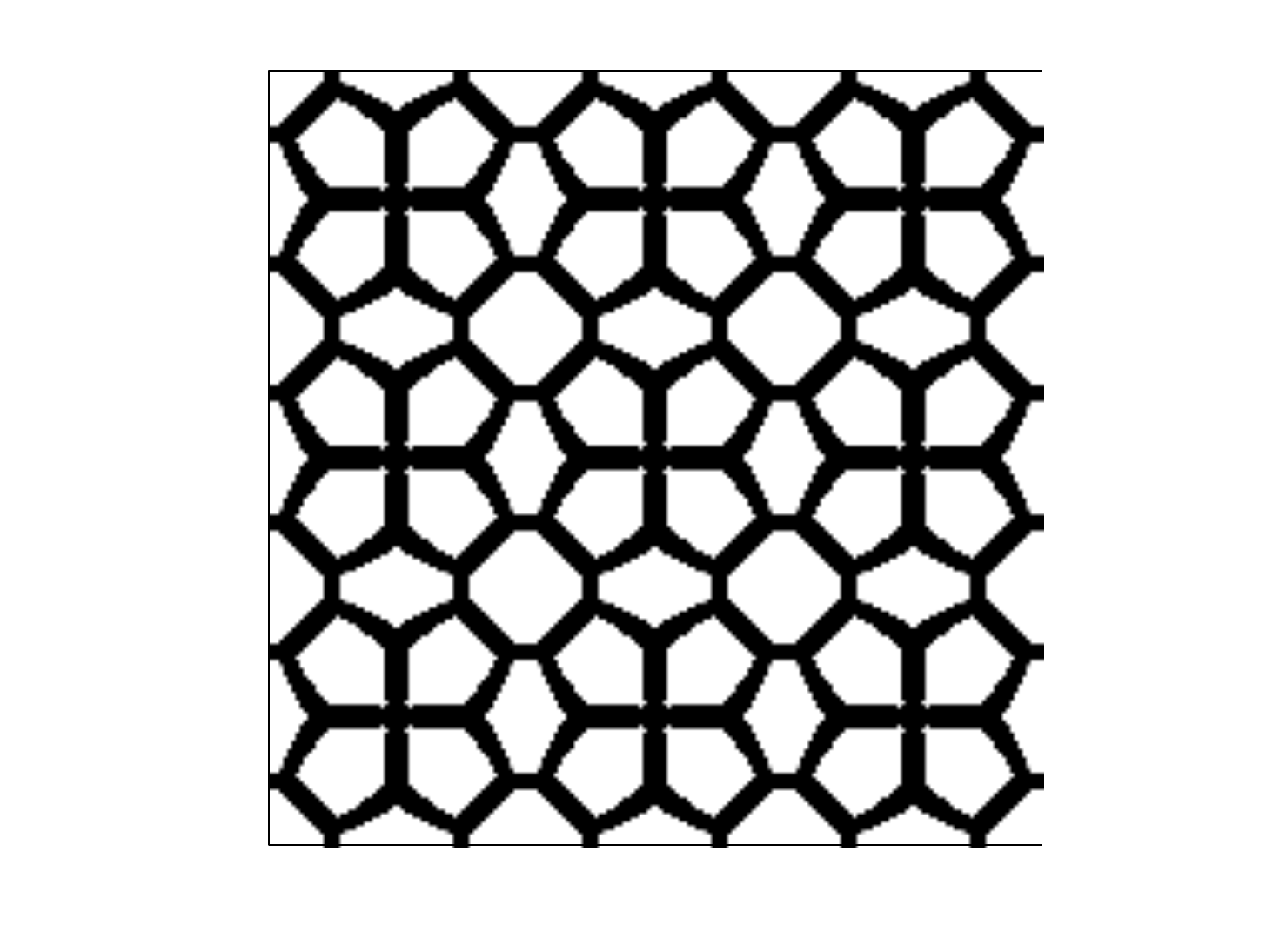}}
\subfigure[Optimal band structure]{
\includegraphics[scale=0.45]{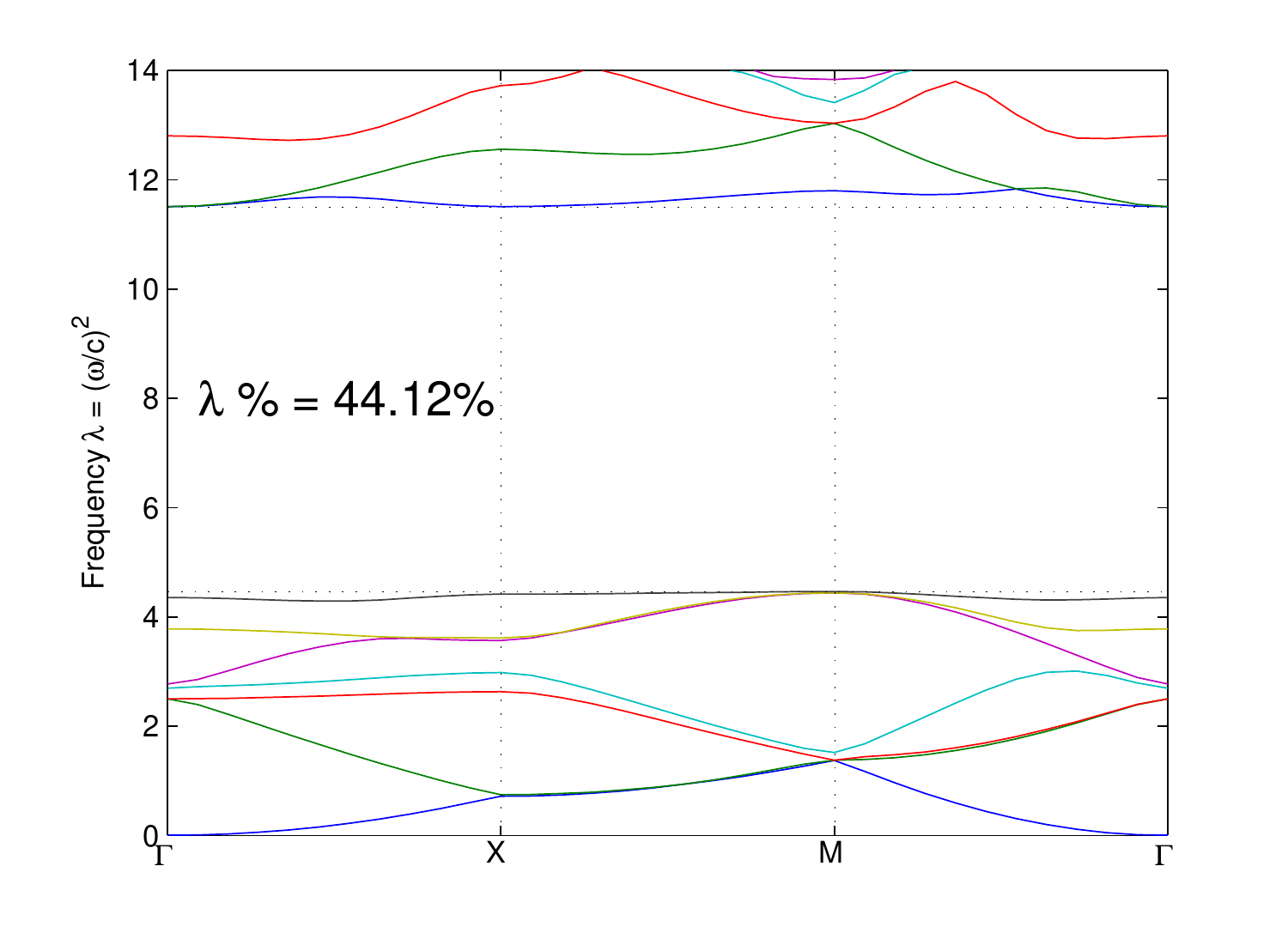}}
\caption{Optimization of band gap between $\lambda^7_{\text{TE}}$ and $\lambda^8_{\text{TE}}$ in the square lattice.}
\label{figOR_TE78}
\end{figure}

\begin{figure}
\centering
\subfigure[Optimal crystal structure]{
\includegraphics[scale=0.45]{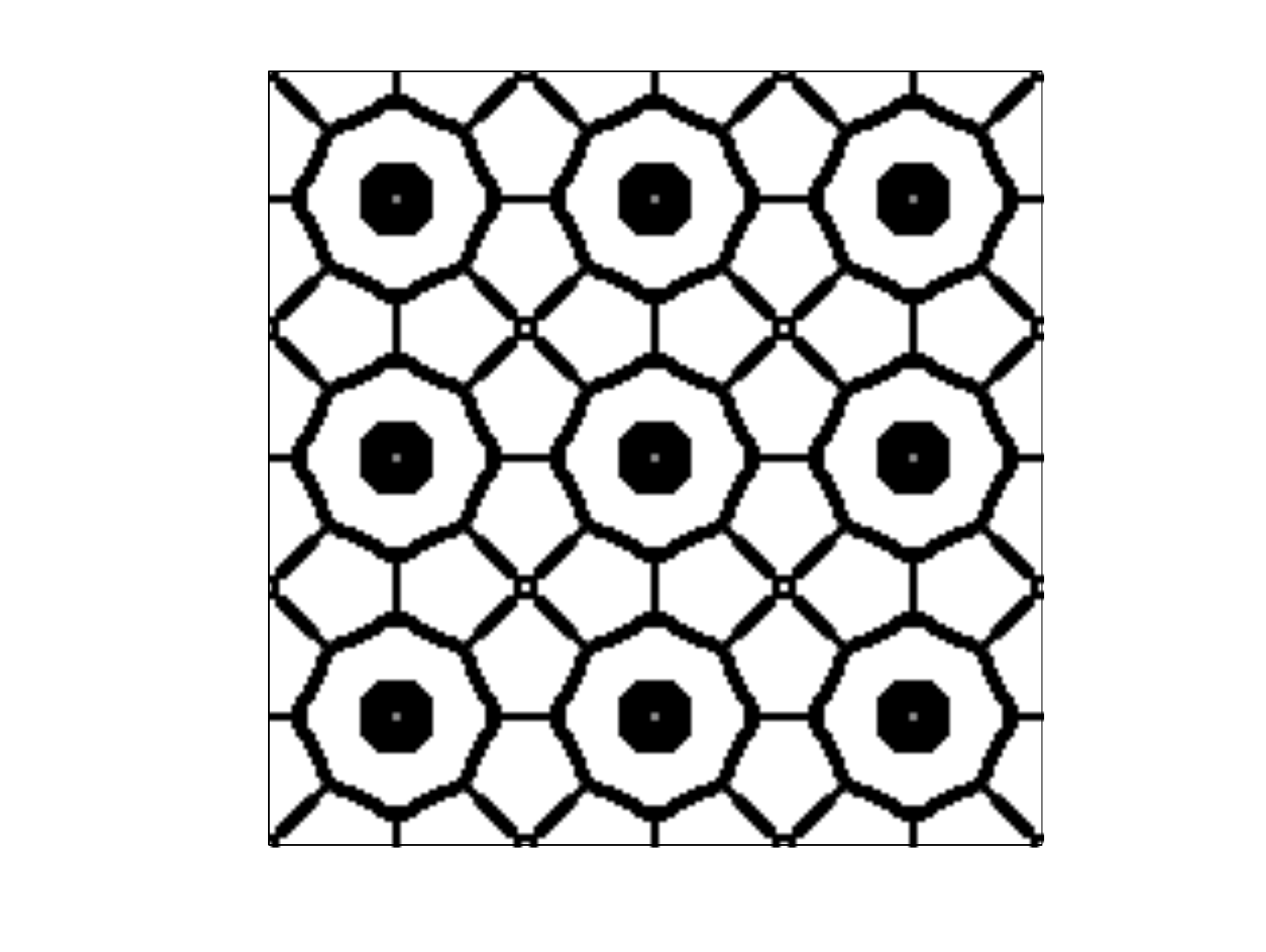}}
\subfigure[Optimal band structure]{
\includegraphics[scale=0.45]{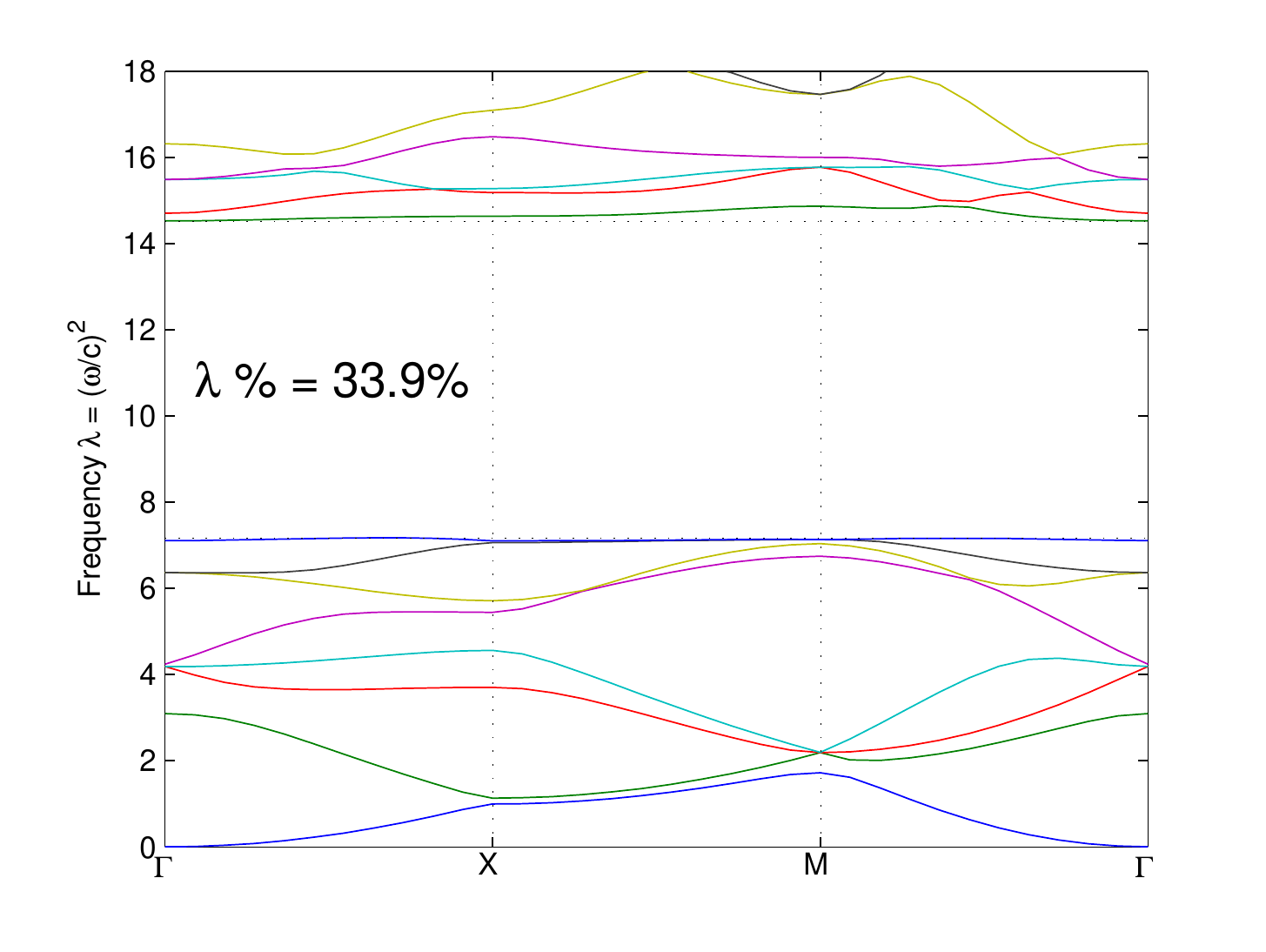}}
\caption{Optimization of band gap between $\lambda^8_{\text{TE}}$ and $\lambda^9_{\text{TE}}$ in the square lattice.}
\label{figOR_TE89}
\end{figure}

\begin{figure}
\centering
\subfigure[Optimal crystal structure]{
\includegraphics[scale=0.45]{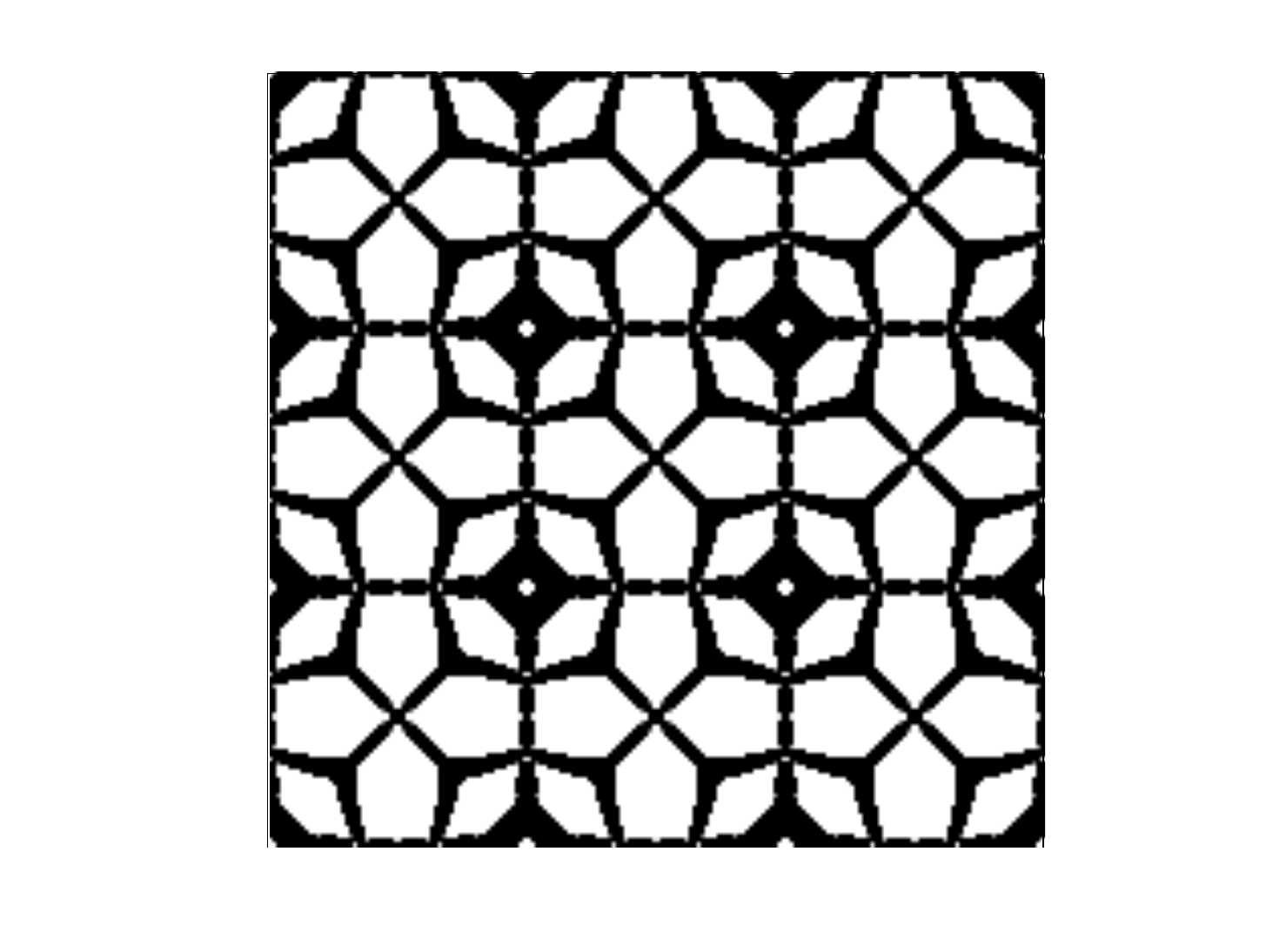}}
\subfigure[Optimal band structure]{
\includegraphics[scale=0.45]{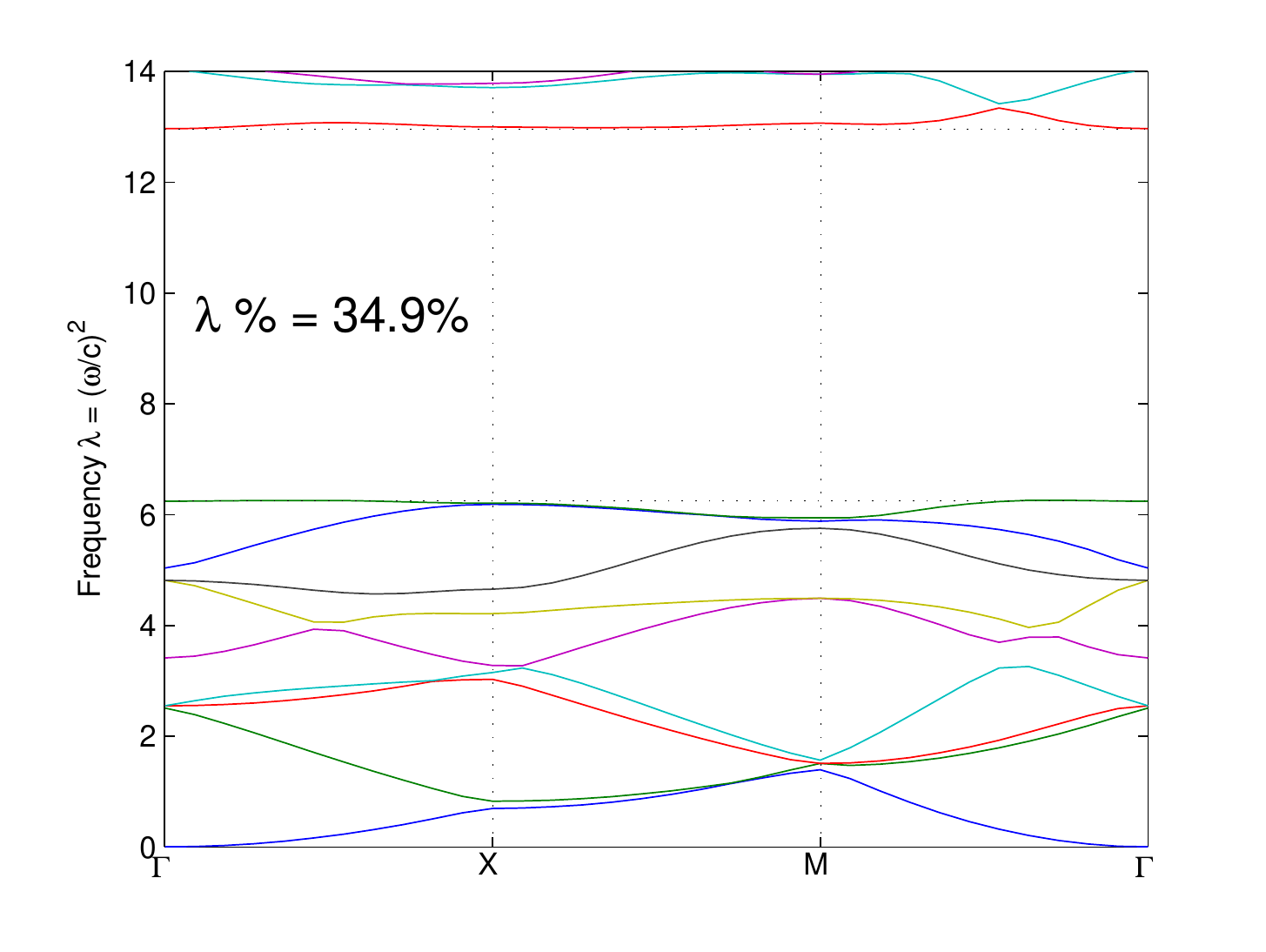}}
\caption{Optimization of band gap between $\lambda^9_{\text{TE}}$ and $\lambda^{10}_{\text{TE}}$ in the square lattice.}
\label{figOR_TE910}
\end{figure}

\begin{figure}
\centering
\subfigure[Optimal crystal structure]{
\includegraphics[scale=0.45]{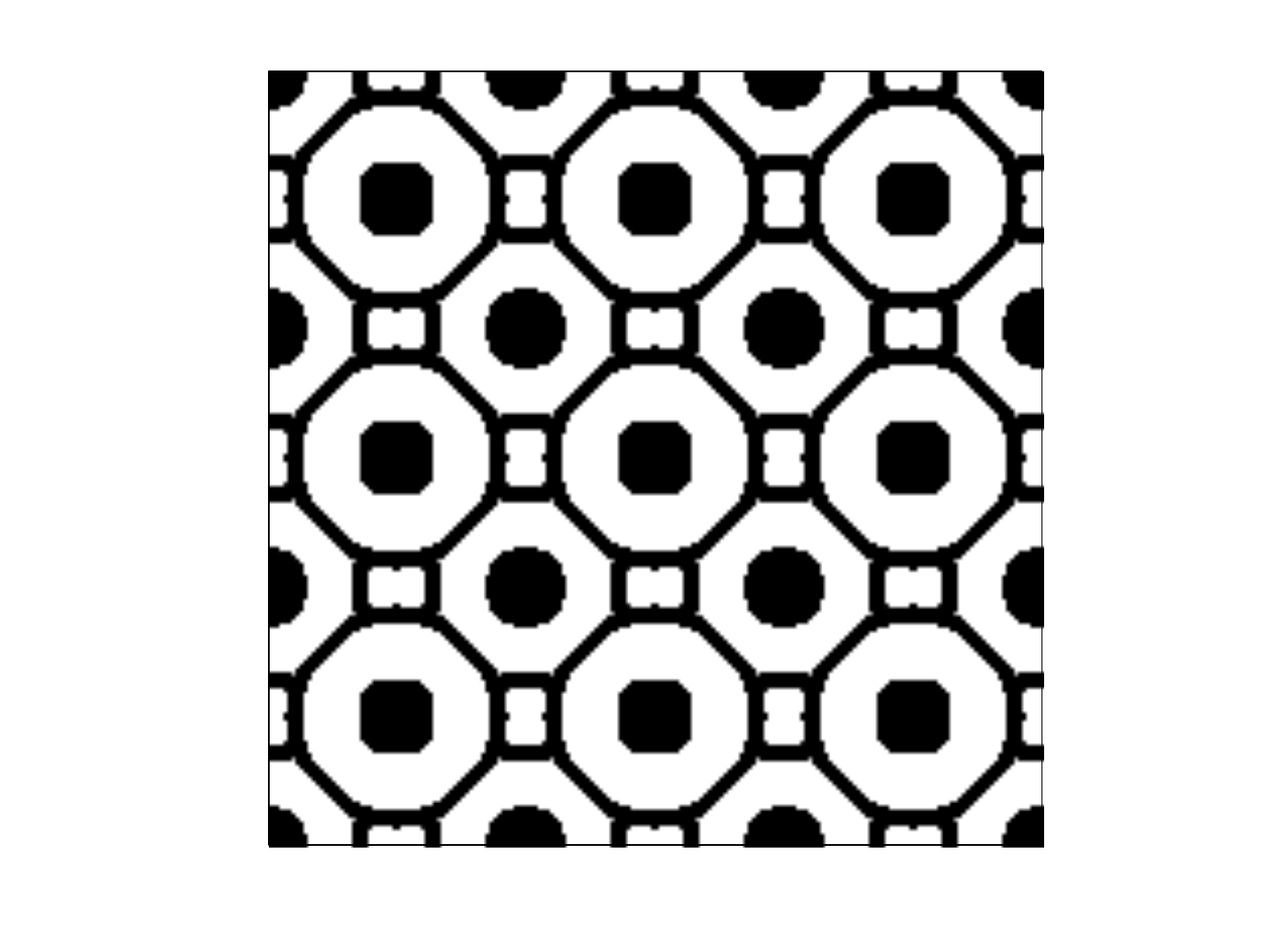}}
\subfigure[Optimal band structure]{
\includegraphics[scale=0.45]{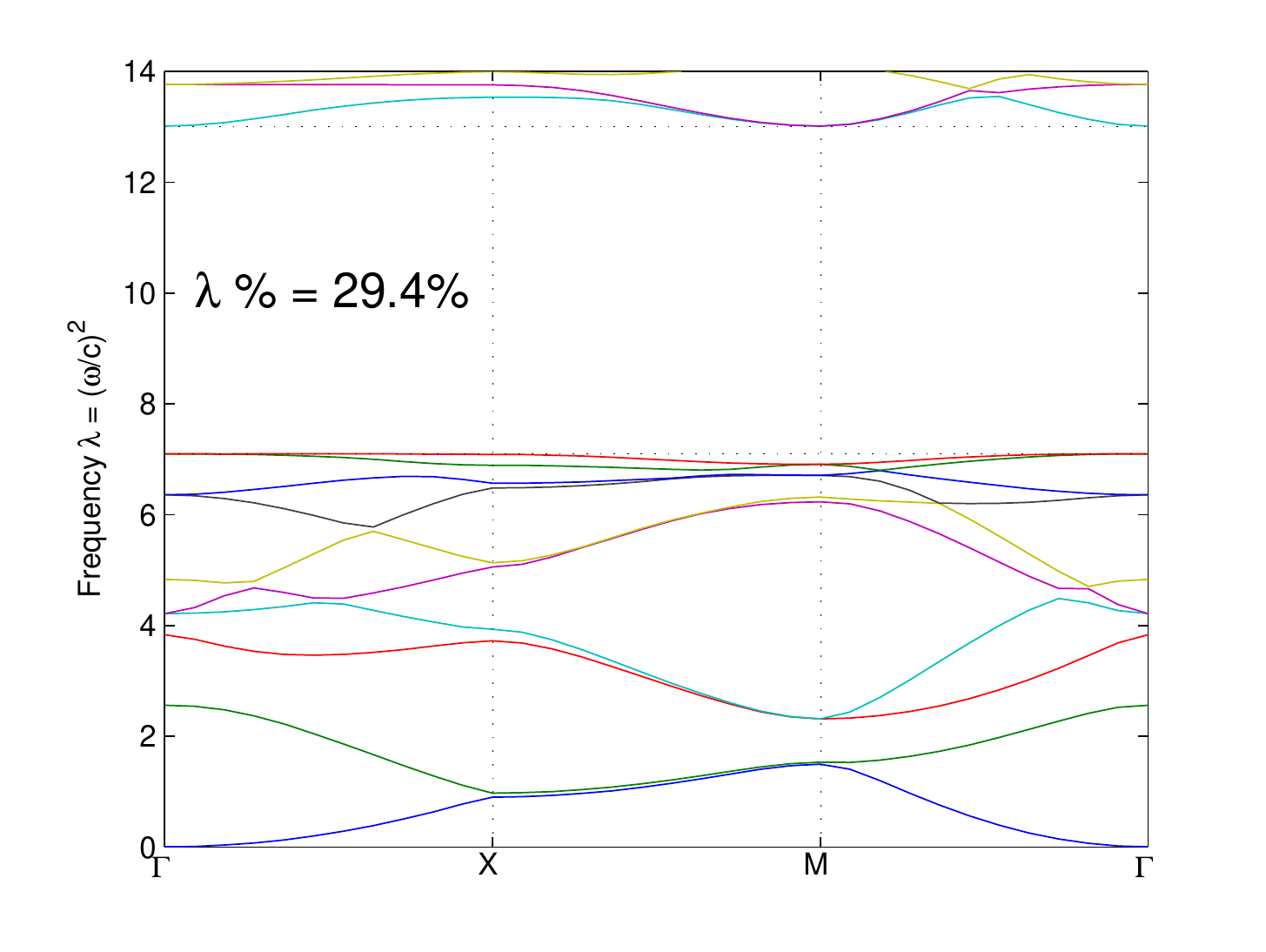}}
\caption{Optimization of band gap between $\lambda^{10}_{\text{TE}}$ and $\lambda^{11}_{\text{TE}}$ in the square lattice.}
\label{figOR_TE1011}
\end{figure}

In Figures \ref{figOR_TM67} through \ref{figOR_TE1011}, we present
only plots of the final optimized crystal structures and the
corresponding band structures for the $6^{\rm th}$ through $10^{\rm
th}$ optimized band gaps for TE and TM polarizations. We see that the optimized TM band gaps are exhibited in isolated high-$\varepsilon$ structures, while the optimized TE band gaps appear in connected high-$\varepsilon$ structures. This observation has also
been pointed out in \cite{joannopoulos2008pcm} (p$75$) \emph{``the TM band gaps are favored in a lattice of isolated high-$\bm
\varepsilon$ regions, and TE band gaps are favored in a connected lattice}'', and observed in \cite{kao2005mbg}
previously. For both TE and TM polarizations, the crystal structures
become more and more complicated as we progress to higher bands. It
would be very difficult to create such structures using physical
intuition alone.  The largest gap-midgap ratio for the TM case is
$43.9\%$ between the seventh and eighth frequency bands, while the
largest ratio for the TE case is $44.1\%$, also between the seventh
and eighth bands. The results presented here are not guaranteed to be
globally optimal, as pointed out in Section 4.2.1.  While most crystal
structures in the TM cases appear similar to those presented in
\cite{kao2005mbg}, we have shown quite different TE structures.  A qualitative
comparison between the two results in the background indicates larger
band gaps (both in absolute value and in the gap-midgap ratio) in our
results.

\section{Conclusions and Future Work}
\label{conclusion}

We have introduced a novel approach, based on reduced eigenspaces and
semidefinite programming, for the optimization of band gaps of
two-dimensional photonic crystals on square lattices.  Our numerical
results convincingly show that the proposed method is very effective
in producing a variety of structures with large band gaps at various
frequency levels in the spectrum.

Since our computational techniques make essential use of the finite
element method, we anticipate that notions of mesh adaptivity can be
easily incorporated into our approach, and thus its computational
efficiency will be improved even further. For example, one can start
with a relatively coarse mesh and converge to a near-optimal solution,
and then judiciously refine the finite element mesh (e.g., refining
elements at the interface of dielectric materials) using the current
optimal solution at the coarser mesh as the new initial configuration.
We intend to explore this approach and report the details and results
in a forthcoming paper.

The main strengths of our proposed approach to solve eigenvalue gap
optimization problem is the fact that SDP-based methods do not require
explicit computation of (sub-)gradients of the objective function (which
are ill-defined in the case of eigenvalue multiplicities), hence
maintaining the regularity of the formulation.  The approach proposed
in this paper can also be readily extended to deal with more general
problems, such as the optimization of photonic crystals in combined TE
and TM fields, optimizing multiple band gaps, dealing with other types
of lattices (e.g. triangular), as well as modeling and optimizing the
design of three-dimensional photonic crystals.

\medskip

\noindent {\bf Acknowledgments.}  We are grateful to Professor Steven
Johnson of the Mathematics Department of MIT for numerous discussions
on this research.  We thank Professor Lim Kian Meng of National
University of Singapore for advising and supporting of H. Men.

\bibliographystyle{plain}
\bibliography{bandgap_ref}

\begin{thebibliography}{10}

\bibitem{alizadeh1995ipm}
F.~Alizadeh.
\newblock {Interior point methods in semidefinite programming with applications
  to combinatorial optimization}.
\newblock {\em SIAM Journal on Optimization}, 5(1):13--51, 1995.

\bibitem{alizadeh1998pdi}
F.~Alizadeh, J.~P.~A. Haeberly, and M.~L. Overton.
\newblock {Primal-dual interior-point methods for semidefinite programming:
  convergence rates, stability and numerical results}.
\newblock {\em SIAM Journal on Optimization}, 8(3):746--768, 1998.

\bibitem{bloch1929qek}
F.~Bloch.
\newblock {{\"U}ber die quantenmechanik der elektronen in kristallgittern}.
\newblock {\em Zeitschrift f{\"u}r Physik A Hadrons and Nuclei},
  52(7):555--600, 1929.

\bibitem{burger87ipt}
M.~Burger, S.~J. Osher, and E.~Yablonovitch.
\newblock {Inverse problem techniques for the design of photonic crystals}.
\newblock {\em IEICE Trans. Electron. E}, 87:258--265.

\bibitem{cances2007fac}
E.~Cances, C.~LeBris, N.~C. Nguyen, Y.~Maday, A.~T. Patera, and G.~S.~H. Pau.
\newblock {Feasibility and competitiveness of a reduced basis approach for
  rapid electronic structure calculations in quantum chemistry}.
\newblock In {\em Proceedings of the Workshop for High-dimensional Partial
  Differential Equations in Science and Engineering (Montreal)}, volume~41,
  pages 15--57, 2007.

\bibitem{charnes1962plf}
A.~Charnes and W.~W. Cooper.
\newblock {Programming with linear functionals}.
\newblock {\em Naval Research Logistics Quarterly}, 9, 1962.

\bibitem{cox2000bso}
S.~J. Cox and D.~C. Dobson.
\newblock {Band structure optimization of two-dimensional photonic crystals in
  H-polarization}.
\newblock {\em Journal of Computational Physics}, 158(2):214--224, 2000.

\bibitem{craven1973dfl}
B.~D. Craven and B.~Mond.
\newblock {The dual of a fractional linear program}.
\newblock {\em Journal of Mathematical Analysis and Applications},
  42(3):507--512, 1973.

\bibitem{doosje2000pbo}
M.~Doosje, B.~J. Hoenders, and J.~Knoester.
\newblock {Photonic bandgap optimization in inverted fcc photonic crystals}.
\newblock {\em Journal of the Optical Society of America B}, 17(4):600--606,
  2000.

\bibitem{fan1995gad}
S.~Fan, J.~D. Joannopoulos, J.~N. Winn, A.~Devenyi, J.~C. Chen, and R.~D.
  Meade.
\newblock {Guided and defect modes in periodic dielectric waveguides}.
\newblock {\em Journal of the Optical Society of America B}, 12(7):1267--1272,
  1995.

\bibitem{fan1998cdf}
S.~Fan, P.~Villeneuve, J.~Joannopoulos, and H.~Haus.
\newblock {Channel drop filters in photonic crystals}.
\newblock {\em Optics Express}, 3(1):4--11, 1998.

\bibitem{floquet1883edl}
G.~Floquet.
\newblock {Sur les equations differentielles lineaires a coefficients
  periodiques}.
\newblock {\em Ann. Ecole Norm. Ser}, 2(12):47--89, 1883.

\bibitem{joannopoulos2008pcm}
J.~D. Joannopoulos, S.~G. Johnson, J.~N. Winn, and R.~D. Meade.
\newblock {\em {Photonic crystals: molding the flow of light}}.
\newblock Princeton university press, 2008.

\bibitem{john1987slp}
S.~John.
\newblock {Strong localization of photons in certain disordered dielectric
  superlattices}.
\newblock {\em Physical Review Letters}, 58(23):2486--2489, 1987.

\bibitem{kao2005mbg}
C.~Y. Kao, S.~Osher, and E.~Yablonovitch.
\newblock {Maximizing band gaps in two-dimensional photonic crystals by using
  level set methods}.
\newblock {\em Applied Physics B: Lasers and Optics}, 81(2):235--244, 2005.

\bibitem{nesterov1994ipp}
Y.~Nesterov and A.~Nemirovskii.
\newblock {Interior-point polynomial algorithms in convex programming}.
\newblock {\em SIAM studies in applied mathematics}, 13, 1994.

\bibitem{pau:046704}
G.~S.~H. Pau.
\newblock Reduced-basis method for band structure calculations.
\newblock {\em Physical Review E (Statistical, Nonlinear, and Soft Matter
  Physics)}, 76(4):046704, 2007.

\bibitem{rayleigh1887mvf}
L.~Rayleigh.
\newblock {On the maintenance of vibrations by forces of double frequency, and
  on the propagation of waves through a medium endowed with a periodic
  structure}.
\newblock {\em Phil. Mag}, 24(147):145--159, 1887.

\bibitem{sigmund2003sdp}
O.~Sigmund and J.~S. Jensen.
\newblock {Systematic design of phononic band-gap materials and structures by
  topology optimization}.
\newblock {\em Philosophical Transactions: Mathematical, Physical and
  Engineering Sciences}, pages 1001--1019, 2003.

\bibitem{SoljacicIb02}
M.~Soljacic, S.~G. Johnson, M.~Ibanescu, Y.~Fink, and J.~D. Joannopoulos.
\newblock Optimal bistable switching in nonlinear photonic crystals.
\newblock {\em Physical Review~E}, 66:055601, 2002.

\bibitem{tutuncu2003ssq}
R.~H. T{\"u}t{\"u}nc{\"u}, K.~C. Toh, and M.~J. Todd.
\newblock {Solving semidefinite-quadratic-linear programs using SDPT3}.
\newblock {\em Mathematical Programming}, 95(2):189--217, 2003.

\bibitem{vandenberghe1996sp}
L.~Vandenberghe and S.~Boyd.
\newblock {Semidefinite programming}.
\newblock {\em SIAM review}, 38(1):49--95, 1996.

\bibitem{weyl1952sp}
H.~Weyl.
\newblock {\em Symmetry}.
\newblock Princeton Univ. Press, 1952.

\bibitem{wolkowicz2000hsp}
H.~Wolkowicz, R.~Saigal, and L.~Vandenberghe.
\newblock {\em {Handbook of semidefinite programming: theory, algorithms, and
  applications}}.
\newblock Kluwer Academic Publishers, 2000.

\bibitem{yablonovitch1987ise}
E.~Yablonovitch.
\newblock {Inhibited spontaneous emission in solid-state physics and
  electronics}.
\newblock {\em Physical Review Letters}, 58(20):2059--2062, 1987.

\bibitem{yang2008obg}
X.~L. Yang, L.~Z. Cai, Y.~R. Wang, C.~S. Feng, G.~Y. Dong, X.~X. Shen, X.~F.
  Meng, and Y.~Hu.
\newblock {Optimization of band gap of photonic crystals fabricated by
  holographic lithography}.
\newblock {\em EPL-Europhysics Letters}, 81(1):14001--14001, 2008.

\bibitem{yanik2005sas}
M.~F. Yanik and S.~Fan.
\newblock Stopping and storing light coherently.
\newblock {\em Phys. Rev. A}, 71(1):013803, Jan 2005.

\end{thebibliography}

\end{document}